\documentclass[acmtog]{acmart}

\usepackage[title]{appendix}
\usepackage{algorithmicx}
\usepackage{algorithm}
\usepackage[noend]{algpseudocode}
\usepackage{amsmath,bm}
\usepackage{wrapfig}
\usepackage{xcolor}
\usepackage{cancel}
\usepackage{bm}
\usepackage{float}
\usepackage{soul}
\usepackage{setspace}
\usepackage{enumitem}
\usepackage{multirow}
\usepackage{blindtext}
\usepackage{xcolor}
\usepackage{textcomp}
\usepackage{booktabs} % For formal tables
\usepackage[most]{tcolorbox}

\setcopyright{acmlicensed}
\acmJournal{TOG}
\acmYear{2023} \acmVolume{42} \acmNumber{6} \acmArticle{215} \acmMonth{12} \acmPrice{15.00}\acmDOI{10.1145/3618374}

\AtBeginDocument{%
  \providecommand\BibTeX{{%
    \normalfont B\kern-0.5em{\scshape i\kern-0.25em b}\kern-0.8em\TeX}}}

\makeatletter
\def\@gobbleappendixname#1\csname thesubsection\endcsname{\Alph{section}.\arabic{subsection}}
\g@addto@macro{\appendix}{\renewcommand{\p@subsection}{\@gobbleappendixname}}
\makeatother

\definecolor{color1}{RGB}{77, 190, 238}
\definecolor{color2}{RGB}{234, 113, 43}
\definecolor{color3}{RGB}{119, 172, 48}
\definecolor{color4}{RGB}{255, 184, 0}
\definecolor{lightbluishgrey}{rgb}{0.76078,0.88235,0.92157}

\newtcolorbox{rounded1}[1][]{enhanced jigsaw,colback=color1,
boxrule=0pt,arc=2mm,auto outer arc,boxsep=0pt,coltext={black},fontupper=\sffamily,
box align=center,
#1}
\newtcolorbox{rounded2}[1][]{enhanced jigsaw,colback=color2,
boxrule=0pt,arc=2mm,auto outer arc,boxsep=0pt,coltext={black},fontupper=\sffamily,
box align=center,
#1}
\newtcolorbox{rounded3}[1][]{enhanced jigsaw,colback=color3,
boxrule=0pt,arc=2mm,auto outer arc,boxsep=0pt,coltext={black},fontupper=\sffamily,
box align=center,
#1}
\newtcolorbox{rounded4}[1][]{enhanced jigsaw,colback=color4,
boxrule=0pt,arc=2mm,auto outer arc,boxsep=0pt,coltext={black},fontupper=\sffamily,
box align=center,
#1}

\def\R{\mathbb{R}}
\def\1{\mathbb{1}}

\def\O{\mathcal{O}}

\newcommand{\ra}[1]{\renewcommand{\arraystretch}{#1}}

\newcommand{\integral}[1]{\overline{#1}}
\newcommand{\quadrature}[1]{\mathring{#1}}
\newcommand{\hybridint}[1]{\mathring{\overline{#1}}}
\newcommand{\hybridquad}[1]{\mathring{#1}_{\textrm{H}}}

\newcommand{\abs}[1]{\ensuremath{ {\lvert #1 \rvert} }}

\begin{document}

\title{%
An Adaptive Fast-Multipole-Accelerated Hybrid Boundary Integral
Equation Method for Accurate Diffusion Curves
}

\author{Seungbae Bang}
\authornote{The publication was written prior to Seungbae Bang joining Amazon.}
\affiliation{
  \institution{University of Toronto}
    \city{Toronto}
  \country{Canada}}
\affiliation{
  \institution{Amazon}
    \city{Sunnyvale}
  \country{USA}}
\email{seungbae@cs.toronto.edu}

\author{Kirill Serkh}
\affiliation{
  \institution{University of Toronto}
  \city{Toronto}
  \country{Canada}}
\email{kserkh@math.toronto.edu}

\author{Oded Stein}
\affiliation{
  \institution{Columbia University}
  \city{New York}
  %\country{USA}
  }
\affiliation{
  \institution{Massachusetts Institute of Technology}
  \city{Cambridge}
  %\country{USA}
  }
\affiliation{
  \institution{University of Southern California}
  \city{Los Angeles}
  \country{USA}}
\email{ostein@usc.edu}

\author{Alec Jacobson}
\affiliation{
  \institution{University of Toronto}
  \city{Toronto}
  %\country{Canada}
  }
\affiliation{
  \institution{Adobe Research}
  \city{Toronto}
  \country{Canada}}
\email{jacobson@cs.toronto.edu}

%%
%% The abstract is a short summary of the work to be presented in the
%% article.
\begin{abstract}
  In theory, diffusion curves promise complex color gradations for
  infinite-resolution vector graphics.
  In practice, existing realizations suffer from poor scaling,
  discretization artifacts, or insufficient support for rich boundary
  conditions.
  Previous applications of the boundary element method to diffusion
  curves have relied on polygonal approximations, which either forfeit the high-order
  smoothness of B\'ezier curves, or, when the polygonal approximation  is extremely 
  detailed, result in large and costly systems of equations that must be solved.
  In this paper, we utilize the boundary integral equation method to
  accurately and efficiently solve the underlying partial differential
  equation. Given a desired resolution and viewport, we then interpolate
  this solution and use the boundary element method to render it.
  We couple this hybrid approach with the fast multipole method on a
  non-uniform quadtree for efficient computation.
  Furthermore, we introduce an adaptive strategy to enable truly scalable
  infinite-resolution diffusion curves.
\end{abstract}




\begin{teaserfigure}
 \centering
 \includegraphics[width=\linewidth]{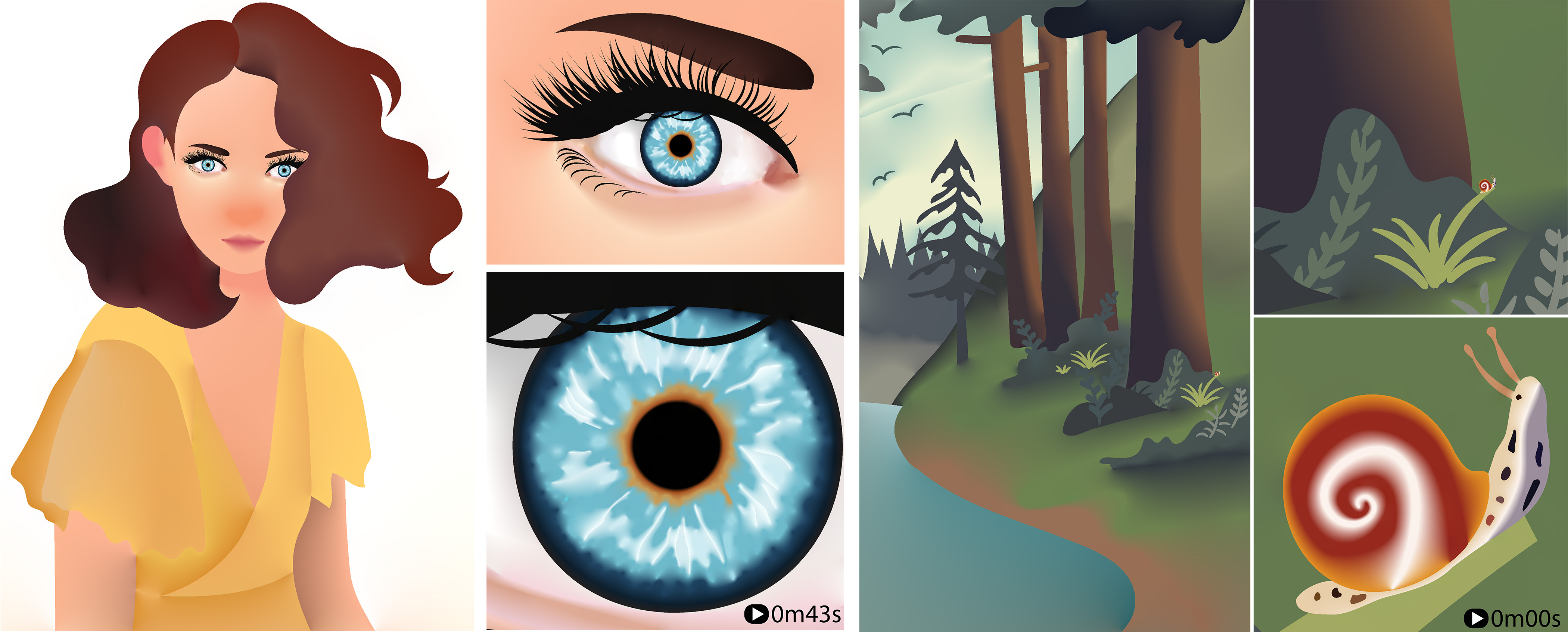}
 \caption{
A diffusion curve image drawn with our method demonstrates that extremely
zoomed-in views maintain high accuracy.
  Please enjoy these images best using a high-resolution digital screen (or
  by printing on billboard).  See the accompanying demo video for full-frame
  zoom-in transition on a given timestamp.
  We obtained the vector graphics image from buysellgraphic.com with purchased commercial license and have modified it to enable extreme zoom-in.
 }
 \label{fig:teaser}
\end{teaserfigure}

\maketitle

\section{Introduction}

% Quickly narrow in on topic

Diffusion curves are primitives for smoothly interpolating color data in
vector graphics images, where the continuous color data is defined to be the
solution to Laplace's equation with boundary values specified along vector
graphics curves.
Laplace's equation is the prototypical elliptic partial differential equation
(PDE), and at first glance it would appear that any numerical method for
elliptic PDEs could potentially be used to solve it, such as finite differences,
finite elements, boundary elements or random walks.
Unfortunately, in practice, diffusion curves present a number of
complications which cause problems in many existing numerical methods.  

Finite difference-based diffusion curve methods 
\cite{orzan2008diffusion,finch2011freeform} rely on lossy rasterization of
boundary data onto a fixed pixel grid, which may either be too dense
(and slow) or too coarse (and inaccurate and aliased) for a desired display
resolution. %
Finite element methods \cite{pang2011fast,JacobsonWS12}
similarly commit to a fixed, albeit adaptive, grid resolution which
simultaneously determines the solution accuracy, solution smoothness, and
boundary curve fidelity.
Both linear elements and isogeometric (curved) elements present
their own respective difficulties.
While popular, linear FEM requires approximating curved B\'ezier curves
by linear segments.
Alternatively, higher-order FEM, with its more complicated functions
spaces \cite{SchneiderHDGPZ18,ilbery2013biharmonic} could be used
on a mesh made of curved elements which conform to boundary curves.
Unfortunately, generating these meshes automatically remains an open
problem with very recent advances \cite{hu2019triwild,MandadC20}.
Furthermore, once meshes are generated, FEM struggles to provide accuracy
near boundary singularities \cite{gopal2019new}.
Unlike many other PDE-based problems in computer graphics, diffusion curves are
rife with both geometric boundary singularities (sharp corners or endpoints of
open curves) and discontinuities in prescribed color values.
Stochastic methods based on random walks, like the recent Walk on Spheres
method \cite{sawhney2020monte}, can overcome some of these difficulties,
however, such methods do not support problems where the boundary conditions
are predominantly Neumann, which are essential to practical applications of
diffusion curves.

An alternative to discretizing the entire image domain is to employ
boundary-only methods, where the color value at every point can be computed from
calculations performed on the boundary alone.  The boundary element method (BEM)
discretizes only the boundary using boundary elements, and can then evaluate the
solution at any point in the domain after a precomputation step which involves
solving an integral equation.  BEM, however, still requires discretization of
the boundary into line segments \cite{van2010high, sun2012diffusion},  
which can lead to resolution problems at the boundary, similar to those
encountered in linear FEM.

We propose a boundary-only method which does not represent the
solution on line segments approximating
the boundary geometry.  Instead, we sample directly from the exact spline
representation of boundary curves using the boundary integral equation
method (BIEM), and solve the associated integral equation in a way that
allows us to color pixels at an arbitrary resolution.  To evaluate the color
data, we interpolate our smooth BIEM solution to a resolution-
and viewport-aware BEM discretization. 
A large part of the calculations required by our method can be precomputed
and, during changes of the viewport, the solution to the BIE only needs to
be re-solved on a sparse set of boundary curves.
We employ the Fast Multipole Method (FMM) to efficiently evaluate the color
data for a large number of curves. In applying the FMM to diffusion
curves, we find that the FMM, as it is typically
presented~\cite{martinsson2019fast} and implemented~\cite{fmm2d}, is not
especially friendly to a graphics audience, and forgoes some
precomputations which we find to be essential in our application.  Thus, we
provide a self-contained presentation of the FMM in the context of diffusion
curves, along with the novel strategies we employ to accelerate our
computations.
Our method, together with our optimized FMM, results in an efficient and
fully adaptive infinite-resolution algorithm for evaluating diffusion
curves, and can be viewed as a hybrid of BEM and BIEM, which
maximally exploits the advantages of both.

\begin{figure}
  \includegraphics[width=0.8\linewidth]{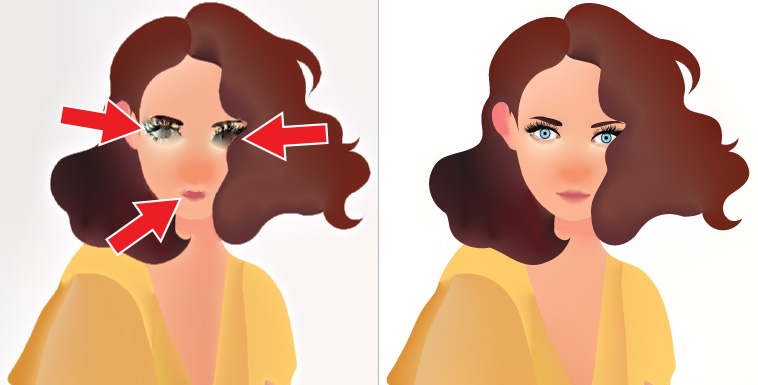}
  %\vspace{-10pt}
  \caption{The original finite-difference method of
  \citet{orzan2008diffusion} exhibits inaccuracies (left), e.g., around the
  eye.  
  Our hybrid, boundary-only method shows accurate solution (right).
  }
\label{fig:compare_FD}
\end{figure}

\section{Related Works}

\subsection{Diffusion Curves}

When Diffusion Curves (DCs) were first introduced,
\citet{orzan2008diffusion} solved Laplace's equation using the
Finite Difference (FD) Method, with follow-up work reformulating the
equation as a constraint problem \cite{bezerra2010diffusion}
on a grid of pixels.  Despite its strengths of simplicity and easy
parallelization, rasterizing the input curve to a pixel domain can lead to
inaccurate results, as shown in Fig.~\ref{fig:compare_FD} (left). 
Follow-up work of \cite{jeschke2009gpu} overcomes this issue by initializing
each pixel to the color of the closest curve point and blending the image
with a Jacobi-like iteration. While resulting images are visually excellent,
they can nonetheless differ slightly from the converged solutions.

\begin{figure}
  \includegraphics[width=\linewidth]{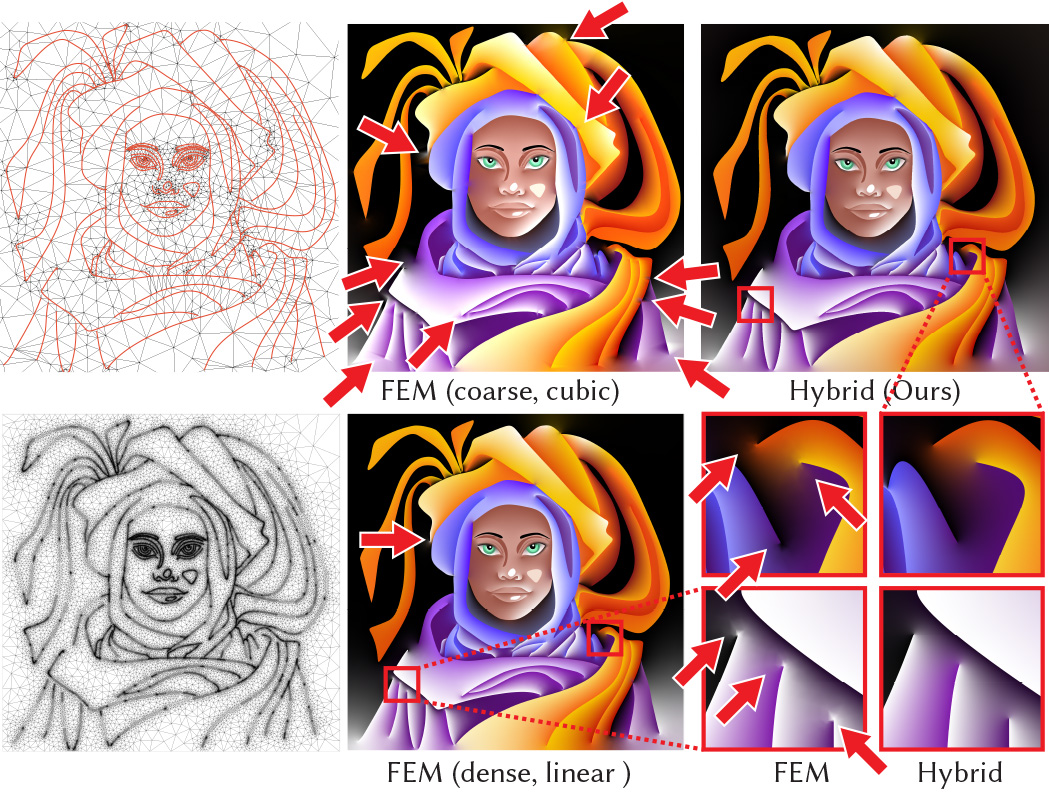}
  \caption{Top row: Bleeding artifacts are pervasive in the FEM results of
  TriWild~\cite{hu2019triwild}, as shown in their Figure 12, reproduced here
  with arrows added. (Used under permission.)
  Bottom row: With much denser triangulation,  with linear
  elements, it still shows unnatural transition of color near end
  point of curves. 
  Our method shows accurate and smooth solution with the same data.}
\label{fig:fem_artifact}
\end{figure}

To overcome some of the problems of the FD method, the Finite
Element Method (FEM) was employed to evaluate DCs
\cite{pang2011fast,takayama2010volumetric} since FEM can more precisely
represent the boundary geometry using constrained triangulation
along curves.  While the boundary can indeed be better represented,
triangulation itself can become burden if the input curves are too numerous
or have complex shapes.
Using the powerful triangulation tool TriWild~\cite{hu2019triwild}, we could not
successfully generate a triangulation of example Fig.~\ref{fig:compare_FD} with
sufficient detail preserved.
Even if triangulation succeeds, FEM still suffers from bleeding artifacts if the
triangulation is not dense enough, as shown in Fig.~\ref{fig:fem_artifact}.
FEM has notoriously poor accuracy near singularities such as re-entrant corners,
which are commonplace in DCs (see, e.g., \cite{gopal2019new}).
While \cite{boye2012vectorial} (Sec.4.3) did present a heuristic method to
circumvent this singularity problem, it does not provide as accurate a
solution as our approach does.

\begin{figure}
  \includegraphics[width=\linewidth]{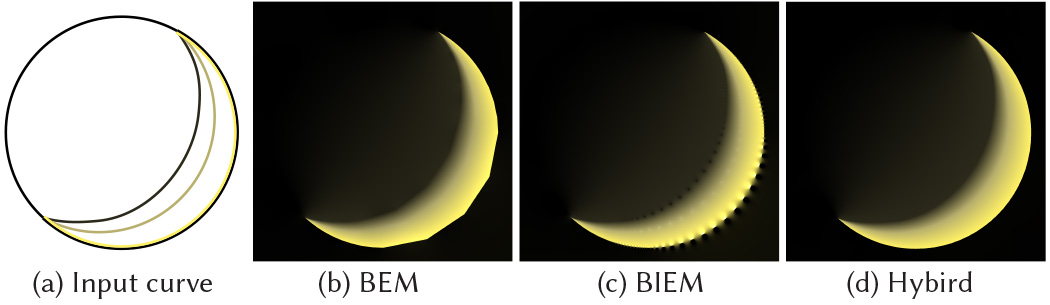}
  \caption{Result comparisons between BEM,  BIEM,  and our Hybrid method.
  BEM suffers from visible polyline, and BIEM shows dotted-looking artifacts
  near the boundary, whereas our method is free from both problems. 
  }
\label{fig:bvp_compare_visual}
\end{figure}

The Boundary Element Method (BEM)~\cite{van2010high, sun2012diffusion} can be
used to avoid triangulation by only discretizing boundary
curves and re-formulating the problem as an integral equation. The evaluation of
color values, which would otherwise be fairly expensive 
with brute force computation, can be accelerated
using the Fast Multipole Method (FMM)~\cite{sun2014fast}.  However, BEM
still suffers from visible polyline discretization, as shown in
Fig.~\ref{fig:bvp_compare_visual} (b).

\begin{figure}
  \includegraphics[width=\linewidth]{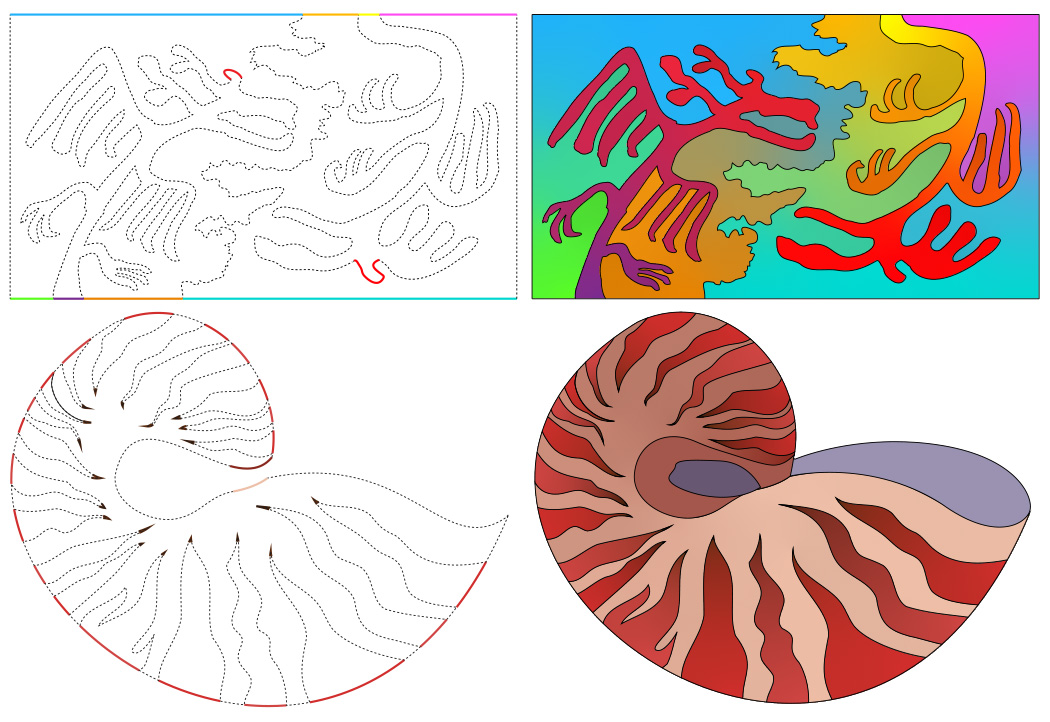}
  \caption{Input diffusion curves, with Dirichlet boundary conditions on
  colored curves and zero Neumann boundary conditions on dotted curves (left),
  and its solution (right), with an example image inspired from
  \cite{hofstadter1979godel} (top) and an example of a Nautilus shell (bottom).
  For the Nautilus shell example, a solid color inner shell region is
  overlayed. Nautilus shell vector graphics art is obtained from buysellgraphic.com with purchased commercial license.  }
\label{fig:neumann}
\end{figure}

Diffusion curves can also be evaluated using stochastic methods.  Stochastic
ray tracing \cite{bowers2011ray} treats the curves as light sources emitting
radiant energy, and determines the color at a pixel by computing the
radiance received at that point.  This method was further combined with FEM
in follow-up work \cite{prevost2015vectorial}.
While stochastic ray
tracing is able to achieve real-time performance using a GPU-based
implementation, it is unable to diffuse colors around corners or obstacles,
resulting in visual differences when compared with diffusion curves
evaluated by other methods.  The fully
meshless Walk on Spheres (WoS) \cite{sawhney2020monte,sawhney2022grid},
on the other hand, does diffuse colors around obstacles. However, WoS has
difficulties with Neumann boundary conditions, and this turns out to be a
major limitation, since such boundary conditions turn out to be exceedingly useful
in practice.
For complicated collections of
input curves, it is difficult to specify Dirichlet boundary conditions on
every single curve. By specifying a zero Neumann boundary condition on a
majority of the input curves, one only needs to specify Dirichlet boundary
conditions on a small subset of curves to create a smooth and natural color
interpolation on the domain,  as shown in Fig.~\ref{fig:neumann}.

Besides the various methods for evaluating diffusion curves, the notion of a
diffusion curve itself has been generalized in several directions.  The typical
definition of the interpolated colors of a diffusion curve is as a harmonic
function; this definition has been generalized to a biharmonic function with FD
\cite{finch2011freeform}, FEM \cite{boye2012vectorial, JacobsonWS12}, and BEM
\cite{ilbery2013biharmonic}.  A blending of two harmonic functions
\cite{jeschke2016generalized} has also been introduced to overcome some of the
unintuitive extrapolation behaviour of biharmonic functions.  Diffusion curves
with harmonic interpolated colors satisfying Laplace's equation have been
generalized to interpolated colors satisfying Poisson's equation \cite{8449116},
where the inhomogeneous term in the Poisson equation was used to provide more
nuanced control over blending and diffusion.  Finally, while
diffusion curves are typically presented as an artistic tool, methods have been
proposed for constructing diffusion curve images from rasterized images
\cite{jeschke2011estimating, xie2014hierarchical,zhao2017inverse}.

Furthermore, it's important to note that, besides diffusion curves, 
there exists a multitude of other approaches to color gradations in
vector graphics representations. 
One such example is the patch-based method \cite{xia2009patch}. However, we will
not cover other approaches in this paper.

\subsection{BEM \& BIEM in Graphics} 

The Boundary Element Method (BEM) reformulates the PDE to be solved as a
boundary integral equation. It then discretizes this integral
equation by approximating the boundary curves by line segments in 2D, or by
approximating the boundary surfaces by triangular elements in 3D. It then
represents the solution to the integral equation as a
piecewise constant function on these line segments or flat surface elements.
BEM was first introduced in the graphics community for real time deformable
objects \cite{james1999artdefo}, followed by ocean wave animation
\cite{keeler2014ocean, schreck2019fundamental}, and surface only liquids
simulation \cite{da2016surface}. BEM can accelerate simulations
while retaining visual accuracy, in those cases where the simulation 
involves only the boundary of the object in question.

The Boundary Integral Equation Method (BIEM) \cite{greengard2009fast} uses
the same integral equation formulations employed by the BEM, with the
difference being that the BIEM represents the curve and the data by
spectrally-accurate quadrature-based discretizations, where by
spectrally-accurate, we mean discretizations for which the approximation
error decays exponentially with the number of degrees of freedom used. This
efficient representation means that a very small number of degrees of
freedom are required to represent the solution to high accuracy.  As far as
we know, there is no work that employs BIEM in the graphics community.  The
BIEM has gained popularity for simulations in mathematical physics due to
its simple quadrature-based integration scheme, its favorable conditioning
properties, and its high accuracy. Although evaluating the solutions of BIEM
close to boundaries presents substantial
challenges~\cite{helsing2008evaluation}, it happens that, in many physics-related
applications, e.g., acoustic scattering, the solution is mainly
desired away from the boundaries.

\subsection{FMM in Graphics}

The Fast Multipole Method (FMM) \cite{greengard1987fast} has been
sporadically explored in the computer graphics community.
\citet{sun2014fast} introduced the FMM for diffusion curves with
a simple uniform quadtree structure.
The FMM has also been used for fast computations of repulsive curves
\cite{yu2021repulsive},  ferrofluids~\cite{huang2019accurate}, 
and fast linking numbers~\cite{qu2021fast}.
Fast summation methods similar to FMM have been employed in graphics, e.g.,
to compute winding numbers \cite{barill2018fast} and to simulate
fluids \cite{zhang2014}.

\section{Overview of Methods}

We begin by formulating a boundary value problem using different discretization approaches: the Boundary Element Method (BEM), the Boundary Integral Equation Method (BIEM), and our Hybrid Method, which combines the strengths of BEM and BIEM. 
We include Neumann boundary conditions in our framework for the diffusion curve problem.

While our proposed Hybrid Method offers notably improved accuracy compared to BEM, BIEM, and previous methods, its computational efficiency is suboptimal if 
applied naively to a complex diffusion curve image. 
To address this, we incorporate the Fast Multipole Method (FMM) for faster computations. The efficiency of the FMM is enhanced by incorporating a non-uniform quadtree approach, departing from the previously used uniform quadtree \cite{sun2014fast}. This efficiency improvement is complemented by the precision gains achieved through quadtree clipping. The density values are obtained through the Generalized Minimum Residual (GMRES) algorithm.

To enable detailed zoom-in into localized parts of a diffusion curve image, we introduce an adaptive strategy that efficiently re-solves local density values without requiring a full re-solve of the entire image.

Finally, we present an anti-aliasing scheme involving weighted integration, leveraging the structure of the non-uniform quadtree.

Our main contributions can be summarized as follows:
% \begin{itemize}[label={}, labelsep=0pt, leftmargin=5pt]
\begin{itemize}
%\item Introduction of BIEM to the computer graphics community.
\item We perform a comprehensive comparison of BEM and BIEM, leading to the development of a hybrid method that effectively utilizes their respective strengths.
\item We adopt the Neumann boundary condition for diffusion curve problems.
\item We implement the FMM using a non-uniform quadtree approach, combined with quadtree clipping, resulting in rapid and accurate computation of diffusion curves.
\item We introduce an adaptive strategy for optimal discretization tailored to the viewport.
\item We introduce an anti-aliasing scheme based on the non-uniform quadtree structure.
\end{itemize}

\section{Boundary Value Problem}
\label{sec:method_bvp}

Before considering the more complicated case of diffusion curves, where
double-sided boundary conditions are specified over a collection of open
curves, we consider the model problem of Laplace's equation on region
$V\subset \R^2$ with a simple, closed boundary $S$:

\begin{equation}
\Delta u = 0 \; \textrm{on} \; V, \quad \text{ subject to one of } \quad 
\begin{cases}
\text{$u=u^*$ on $S$},\\
\text{$\frac{\partial u}{\partial n}=\psi^*$ on $S$}.
\end{cases}
\label{eq:laplace_eqn}
\end{equation}

For simplicity of presentation, we will only consider the Dirichlet boundary
condition until Sec.~\ref{sec:neumann_bc}.

We discuss three approaches to solving this problem: 
the boundary element method (BEM),  the boundary integral
equation method (BIEM), and our newly proposed hybrid of BEM and BIEM.

\subsection{Boundary Integral Equation}

Both the BEM and the BIEM reformulate the underlying PDE over the volume $V$
as boundary integral equations over the boundary $S$. The key idea is to use
a representation involving the free-space Green's function, which ensures
that the candidate solution always satisfies the PDE. The problem is thus
reduced to enforcing the correct boundary conditions on $S$.

\subsubsection{Green's function}

The free space Green's function $G$ is defined to be
the solution to the Laplace equation
\begin{equation}
\Delta G(p,q) = \delta (p,q),
\label{eq:greens_func}
\end{equation}
where $p,q \in \mathbb{R}^2$.  The Dirac delta function $\delta (p,q)$
represents a unit impulse at the source point $p$, and $G(p,q)$ represents
the response at the point $q$ due to that source.

The Green's function for Laplace's equation in two-dimensional Euclidean
space, as well as its directional derivative, are well known to be
\begin{equation}
G(p,q) = -\frac{\log(\| p-q \|)}{2\pi}
\label{eq:green_G}
\end{equation}
and
\begin{equation}
F(p,q) = \frac{\partial G(p,q)}{\partial n(p)} = -\frac{(p-q)\cdot
n(p)}{2\pi \|p-q\|^2},
\label{eq:green_F}
\end{equation}
respectively.  Where $n(p)$ is the normal vector at $p$.

\subsubsection{Integral Equation}
%
%\paragraph{Indirect Approach}
Using the free space Green's function, we can convert the boundary value
problem Eq.~\ref{eq:laplace_eqn} into its Boundary Integral Equation (BIE)
formulation.
Consider the so-called single layer potential, which represents our
candidate solution $u$ as an integral of the Green's function over a
boundary density $\sigma$:
\begin{equation}
u(x) = \int_S G(p,x)\sigma(p)dS(p), \quad \forall x \in V.
\label{eq:int_V}
\end{equation}
Letting $x$ approach the boundary
$S$, we obtain the following BIE, which we can solve for the unknown density
$\sigma(p)$ on the boundary given Dirichlet boundary values $u^*(q)$:
\begin{equation}
u^*(q) = \int_S G(p,q)\sigma(p)dS(p), \quad \forall q \in S.
\label{eq:int_S}
\end{equation}

The process of solving the boundary value problem Eq.~\ref{eq:laplace_eqn}
using its BIE formulation can be broken into two distinct stages.
We call process of solving for the density $\sigma(p)$ using
Eq.~\ref{eq:int_S} the \textbf{solution} stage, and the process of
evaluating our solution $u(x)$ on domain using formula Eq.~\ref{eq:int_V}
the \textbf{evaluation} stage.

It is also possible to represent $u(x)$ using a so-called double-layer
potential, where $G(p,q)$ is replaced by $F(p,q)$.  For simplicity, we
consider only the case of the single-layer potential here.

\subsection{Boundary Element Method}
 \label{sec:bem}
We can apply the boundary element method to discretize BIEs in order to
solve them numerically. Suppose, without any loss of generality, that
the boundary consists of a single curve $S$.  We begin by discretizing $S$
into line segments $\overline{S}_j$.  Then we assume that the density value
$\sigma_j$ is constant on each line segment.  The BIE Eq.~\ref{eq:int_S} can
be expressed as:
 \begin{equation}
 u^*(q) = \sum_{j=1}^{s} \int_{\integral{S}_j} G(p,q)dS(p)\, \sigma_j ,
 \end{equation}
where $s$ is the number of boundary elements, and the integrals
$\int_{\integral{S}_j} G(p,q)dS(p)$
are computed analytically using well-known formulas that depend on
$\integral{S}_j$ being a line segment.  We have $s$ unknowns $\sigma_j$, and
so we need at least $s$ equations to determine a unique solution.
We choose to evaluate $u^*(q)$ at the midpoint of each segment, which we
denote by $q_i$, to
arrive at the system of equations
  \begin{equation}
 u^*(q_i) = \sum_{j=1}^{s} \int_{\integral{S}_j} G(p,q_i)dS(p)\sigma_j,  \quad
 \text{for each line segment $i$}.
 \label{eq:bem_S}
 \end{equation}
In matrix form, this system of equations is
 \begin{equation}
\integral{\mathbf{u}}^*=\integral{\mathbf{G}} \integral{\bm{\sigma}},
\label{eq:bem_mat_S}
 \end{equation}
where $\integral{\mathbf{u}}^*, \integral{\bm{\sigma}} \in \mathbb{R}^{s}$ are
the column vectors of boundary values and density values, and
$\integral{\mathbf{G}}  \in \mathbb{R}^{s \times s}$ is a (dense) matrix with elements
$\integral{\mathbf{G}}_{ij}=\int_{\integral{S}_j} G(p,q_i)dS(p)$.
After having obtained the density values $\integral{\bm{\sigma}}$ on the
boundary $S$, we can evaluate $u(x) \in V$ using the formula
 \begin{equation}
u(x) = \sum_{j=1}^{e} \int_{\integral{S}_j} G(p,x)dS(p) \sigma_j,
\label{eq:bem_V}
 \end{equation}
where $e = s$ is the number of boundary elements.
Formulas for the analytic integration of Green's functions on line segments
are detailed in Appendix~\ref{app:analytic_int_G}.

\subsection{Boundary Integral Equation Method}

The Boundary Integral Equation Method (BIEM) can accurately represent
continuous functions defined on curved boundaries without any lossy
approximations to the boundary geometry, in contrast to how BEM approximates
$S$ with linear segments.  Functions are represented using carefully
chosen discretizations based on quadrature formulas, and are interpolated by
mapping their values at the discretization points to the coefficients of
spectral expansions. The rapid convergence of quadrature-based approximations
means that functions can be represented with minimal loss of accuracy.

\subsubsection{Integral Equation}
 
The integral equation Eq.~\ref{eq:int_S} can be written as a system of
equations by discretizing the boundary data at Gauss-Legendre nodes:
\begin{equation}
u^*_i = \int_S G(p,q_i)\sigma(p)dS(p),
\label{eq:int_biem_S}
\end{equation}
where, without loss of generality, we assume the geometric boundary curve is
given by a function $\gamma(t)\colon [-1,1] \to \mathbb{R}^{2}$, 
$q_i=\gamma(t_i)$ are the sampled Gauss-Legendre quadrature points,
and $u^*_i = u^*(q_i)$.
Evaluating the integrals in Eq.~\ref{eq:int_biem_S} requires some extra care,
since the Green's function $G(p,q)$ has a logarithmic singularity
at the point $p=q$. It turns out that, for each target point $q_i$, it is
possible to construct special-purpose quadrature nodes $t_{ij}$ and weights
$w_{ij}$, for $j=1,2,\ldots,g_i$, such that each integral above is evaluated
accurately:
\begin{equation}
u^*_i  = \sum_{j=1}^{g_i} w_{ij} G(p_{ij},q_i)\sigma(p_{ij}),
\label{eq:biem_S}
\end{equation}
where $p_{ij}=\gamma(t_{ij})$ are the sampled special-purpose quadrature points
(see, for example,~\cite{kolm2001numerical}).

Discretizing the density value $\quadrature{\bm{\sigma}}$ at the
Gauss-Legendre quadrature points $p_j = \gamma(t_j)$, we can approximate the
continuous function $\sigma(p)$ appearing in Eq.~\ref{eq:biem_S} by solving
for the coefficients of its corresponding Legendre expansion:
\begin{equation}
\mathbf{c} = \quadrature{\mathbf{P}}^{-1}  \quadrature{\bm{\sigma}}.
\label{eq:coefs_interp}
\end{equation}
We can then evaluate the density value $\sigma(p)$ by evaluating Legendre
polynomials:
\begin{equation}
\sigma(p) = \sum_{i=1}^{g} c_i P_{i-1} (p).
\label{eq:legendre_poly}
\end{equation}
Ultimately, this procedure can be written in matrix form:
\begin{equation}
\quadrature{\mathbf{u}}^*= \quadrature{\widetilde{\mathbf{G}}} 
\quadrature{\bm{\sigma}},
\label{eq:biem_mat}
\end{equation}
where $\quadrature{\mathbf{u}}^*, \quadrature{\bm{\sigma}} \in
\mathbb{R}^{g}$, and $\quadrature{\widetilde{\mathbf{G}}} \in \mathbb{R}^{g
\times g}$ is a matrix constructed row-by-row by combining the
quadrature Eq.~\ref{eq:biem_S} with the interpolation described
by Eq.~\ref{eq:coefs_interp} and Eq.~\ref{eq:legendre_poly}.

Because we have a continuous function $\sigma(p)$,$\forall p \in S$, on the
boundary, when it comes to evaluating the solution using the
formula Eq.~\ref{eq:int_V}, we are not bound to use the same Gauss-Legendre
quadrature approximation.  Instead, we can evaluate
\begin{equation}
u(x) = \sum_{j=1}^{e} w_j G(p_j,x)\sigma(p_j),
\label{eq:biem_V}
\end{equation}
where $e$ is the number quadrature points for \emph{evaluation},
$p_j$ are the sampled Gauss-Legendre quadrature points
corresponding to the roots of the $e$-th order Legendre polynomial, with
$w_j$ the corresponding quadrature weights.  Note that $e$ does not have to
be equal to the number of quadrature points $g$ at the \emph{solution}
stage.  We call this interpolation process the \textbf{interpolation} stage.

\begin{figure}
  \includegraphics[width=\linewidth]{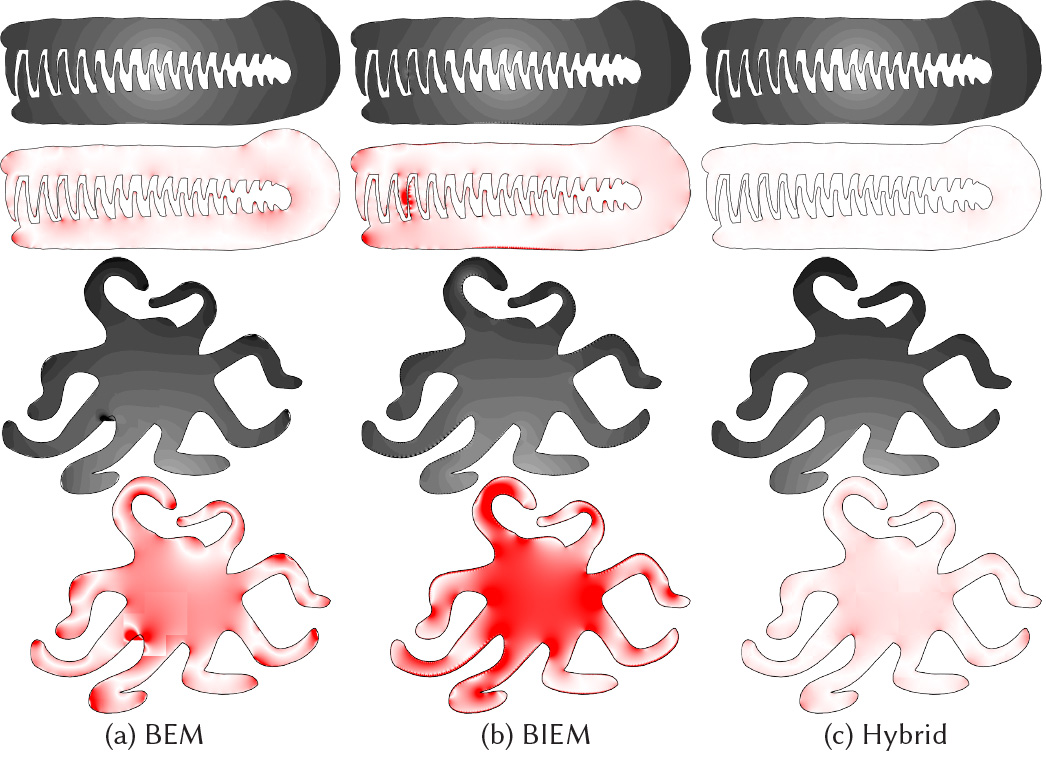}
  \caption{Accuracy comparison between BEM, BIEM, and our Hybrid method.  
  The boundary values are constructed by placing a single source Green's
  function in the middle of the figure.
  The solutions should exactly match the potential induced by
  that Green's function.  The second and last rows show
  the error relative to the ground truth, highlighted in red color.
  }
\label{fig:bvp_compare_accuracy}
\end{figure}

Unfortunately, unlike in the BEM, evaluating the solution using a quadrature
approximation like the one above results in artifacts near the quadrature
points (see the image close to the boundary curves of
Fig.~\ref{fig:bvp_compare_visual}~(c) and
Fig.~\ref{fig:bvp_compare_accuracy}~(b)).  Refining using additional
quadrature points only somewhat alleviates the problem.  Additionally, if
the boundary is composed of multiple different curves, and some of them are
very close to one another, then the kernel $G(p,q)$ can be
close-to-singular, and the evaluation of the integrals in Eq.~\ref{eq:int_S}
by quadrature becomes inaccurate (see near the tip of the moon in
Fig.~\ref{fig:bvp_compare_visual}~(c)).  This situation may be seen as
pathological from the point of view of physical simulations, but it is
commonplace for an artist to create closely positioned diffusion curves as a
technique for achieving high contrast color changes.

The question of how best to evaluate the potential induced by a
continuous density $\sigma(p)$ has been the subject of much recent research
(see, for example,~\cite{helsing2008evaluation,af2021accurate}).
Historically, this has not been a major issue for BIEM, since many
important physical applications of BIEs, e.g., acoustic and electromagnetic
scattering, often do not require the evaluation of the solution close to
boundaries.

\section{Accurate Discretization with Hybrid Method}
\label{sec:hybrid}

Our proposed method combines the advantages of the BEM and BIEM
approaches into a hybrid technique. We start by comparing
these two techniques.

\subsection{Comparison between BEM and BIEM} 
 
For the comparative analysis between BEM and BIEM, we will divide the
diffusion curves algorithm into 3 steps: 
(1) \emph{solution}, (2) \emph{interpolation}, (3) \emph{evaluation}.
BEM uses analytic integration on line segments for \textit{solution} and \textit{evaluation}
but does not have any \emph{interpolation} stage. BIEM uses quadrature-based
integration for \textit{solution} and \textit{evaluation}, and it uses
Legendre polynomial interpolation on density values to populate quadrature
points for \emph{evaluation}. 

BEM has the limitation that the number of degrees of freedom
representing the piecewise constant density $\sigma$ is bounded by the
number of elements in the spatial discretization of the boundary curves.
BIEM is free from this limitation, and the number of degrees of freedom in the
representation of the continous density $\sigma$ is decoupled from the 
 number of
quadrature points $e$ used for \emph{evaluation}. 
On the other hand, BIEM has the limitation that it is inaccurate when curves
are close-to-touching in the \emph{solution} stage, and has artifacts in the
induced potential near the quadrature points in the \emph{evaluation} stage.
BEM, however, is free from both of these problems, since it uses analytic
integration along line segments.

\subsection{Combination of BEM and BIEM}
\label{sec:comb_bem_biem}

We propose to combine these two methods, inheriting the strengths of both.
We discretize both the solution $\sigma$ and the boundary data $u^*$ at
Gauss-Legendre nodes, as in BIEM. However, we also introduce the BEM in two
places. In order to evaluate integrals of the form Eq.~\ref{eq:int_S} in the
\emph{solution} stage, we interpolate the density using
formulas Eq.~\ref{eq:coefs_interp} and Eq.~\ref{eq:legendre_poly} to a BEM-like
approximation, which corrects the shortcoming of BIEM for close-to-touching
curves. Once we have solved for the solution $\quadrature{\bm{\sigma}}$
at the quadrature nodes, we evaluate the
potential by once again interpolating to a BEM-like approximation, which
corrects the shortcoming of BIEM with respect to artifacts in the induced
potential.

We begin by discretizing the boundary data at Gauss-Legendre nodes, leading
to the system of equations Eq.~\ref{eq:int_biem_S}. We then discretize the
boundary curve $S$ into $s$ line segments $\overline{S}_j$. If the density values
$\integral{\bm{\sigma}} \in \mathbb{R}^{s}$
on these
line segments are known, then we can write Eq.~\ref{eq:int_S} as
\begin{equation} u^*(q_i) = \sum_j^{s} \int_{\overline{S}_j}
G(p,q_i)dS_j(p)\sigma_j,
\quad \text{for each quadrature point $i$}.  \label{eq:biem_bem_S}
\end{equation}
In matrix form:
  \begin{equation} \quadrature{\mathbf{u}}^*=
  \hybridint{\mathbf{G}}\integral{\bm{\sigma}}, \end{equation}
Where $\quadrature{\mathbf{u}}^* \in \mathbb{R}^{g}$ are given
boundary values at quadrature points, $\integral{\bm{\sigma}} \in
\mathbb{R}^{s}$ are
density value on line segments of the boundary,
and $\hybridint{\mathbf{G}} \in
\mathbb{R}^{g \times s}$. 

Since we choose to discretize the solution $\sigma$ at Gauss-Legendre
nodes like in the BIEM, we recover the density values $\integral{\bm{\sigma}}$
by using Legendre polynomial interpolation. 
Computing the coefficients of the Legendre expansion of $\sigma$ by
$\mathbf{c} =  \quadrature{\mathbf{P}}^{-1}  \quadrature{\bm{\sigma}}$, we can
evaluate the density value on the midpoint of each line segment $\overline{S}_j$ by the formula
$\integral{\bm{\sigma}}=\integral{\mathbf{P}}  \mathbf{c}$ ,where
$\integral{\mathbf{P}} \in  \mathbb{R}^{s \times g}$ is the Legendre
interpolation matrix constructed by evaluating the Legendre polynomials at
$\integral{\mathbf{t}}$, which is a vector of curve parameter values corresponding
to the midpoints of the line segments $\overline{S}_j$.
Hence, we have the relation $\integral{\bm{\sigma}}=\integral{\mathbf{P}}
\quadrature{ \mathbf{P}}^{-1} \quadrature{\bm{\sigma}}$.

We can thus express our system in matrix form in terms of
$\quadrature{\bm{\sigma}}$ as:
  \begin{equation}
  \quadrature{\mathbf{u}}^*=\underbrace{\hybridint{\mathbf{G}} \,
  \integral{\mathbf{P}}  \quadrature{\mathbf{P}}^{-1}
  }_{\hybridquad{\mathbf{G}}} \quadrature{\bm{\sigma}},
    \label{eq:hybrid_sys}
 \end{equation}
where $\hybridquad{\mathbf{G}} \in \mathbb{R}^{g \times g}$.  In order for
$\hybridquad{\mathbf{G}}$ to have full rank, the number of quadrature points
$g$ must be $\leq$ the number of line segments $s$.  Note that, regardless
of the size of $s$, the dimensionality of the system is $g\times g$.  This
is beneficial for us, as the matrix that needs to be inverted is much
smaller than the corresponding matrix for BEM, $\integral{\mathbf{G}}  \in
\mathbb{R}^{s \times s}$.

Once we solve the system Eq.~\ref{eq:hybrid_sys}, we have, by Legendre
polynomial interpolation, a density value $\sigma(p)$ that can be evaluated
anywhere on the curve.  At the \emph{evaluation} stage, we employ the
BEM-like approach of Eq.~\ref{eq:bem_V}, and now we can use an arbitrary
number of line segments $e$, that is independent both of the number line
segments $s$ used at \emph{solution} stage and the number of quadrature
points $g$ used to represent the solution. Note that we must use arc length
when we integrate over line segments, in order to have consistent
integration lengths between the \emph{solution} and \emph{evaluation} stages
(see the details in Appendix~\ref{app:arc_len_int}). We also use arc length
parametrization when we construct our BEM-like discretization, so that we
have line segments of equal arc length when we subdivide each curve  (see
the details in Appendix~\ref{app:arc_len_param}).

Our method is free from both the visible polyline discretization problem of
BEM for a system of the same size, and also from the artifacts around
quadrature points that are found in BIEM (see
Fig.~\ref{fig:bvp_compare_visual}).  Our method shows the most accurate
results when the number of degrees of freedom in the solution stage and the
evaluation stage are both kept fixed (see
Fig.~\ref{fig:bvp_compare_accuracy}).
In Fig.~\ref{fig:bvp_compare_visual}, we set $s=e=4$ for BEM, $g=4,e=20$ for
BIEM,  and $s=20,g=4,e=20$ for our Hybrid method.  In
Fig.~\ref{fig:bvp_compare_accuracy}, we set $s=e=8$ for BEM, $g=8,e=40$ for
BIEM,  and $s=40,g=8,e=40$ for our Hybrid method.

\begin{figure}
  \includegraphics[width=\linewidth]{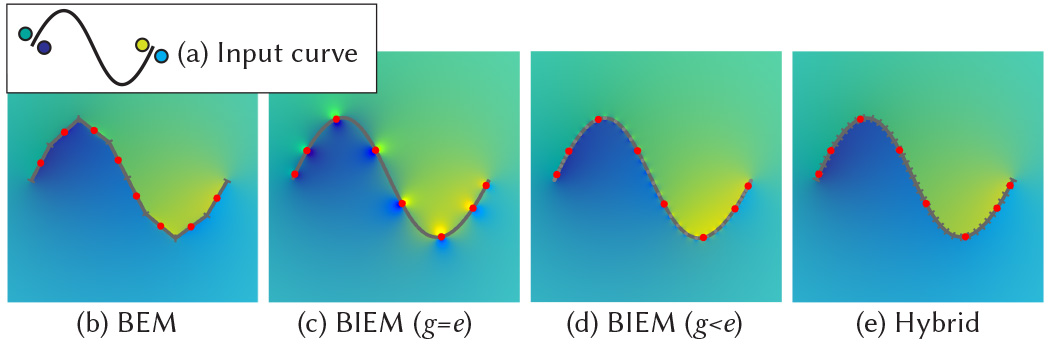}
  \caption{
Diffusion curve with double-sided boundary condition given a single B\'ezier
curve (a).  Solution comparison with BEM with $s=8$ line segments (b), BIEM
with $g=e=8$ (c), BIEM with $g=8, e=40$, and our Hybrid method with
$g=8,s=e=40$ (e).
  }
\label{fig:bvp_comparision_curve}
\end{figure}

\subsection{Double-Sided Boundary Condition}
\label{sec:double_bc}

Up to this point, we have been formulating our equations using the single
layer potential of Eq.~\ref{eq:int_V} for simplicity of presentation.
However, in order to specify two different boundary conditions on each side
of an open curve, we must add to our single layer potential representation
of Eq.~\ref{eq:int_V} a so-called double layer potential, in which the kernel
of Eq.~\ref{eq:int_V} is replaced by $F(p,x)$ from Eq.~\ref{eq:green_F}. Our candidate solution is
thus represented as
\begin{equation}
u(x) = \int_S [G(p,x)\sigma(p)+F(p,x)\mu(p)]dS(p), \forall x \in V.
\label{eq:int_V_dbnd}
\end{equation}

Letting $x$ approach the boundary $S$, we obtain the following BIE:
\begin{equation}
\begin{split}
u^{*}_{+}(q) = \int_S [G(p,q)\sigma(p)+F(p,q)\mu(p)]dS(p)+\frac{1}{2}\mu(q), \forall q \in S, \\
u^{*}_{-}(q) = \int_S [G(p,q)\sigma(p)+F(p,q)\mu(p)]dS(p)-\frac{1}{2}\mu(q), \forall q \in S,
\end{split}
\label{eq:int_S_dbnd_twoside}
\end{equation}
where the $+$,$-$ subscripts indicate on which side the limit is taken.  The
terms $\pm\frac{1}{2}\mu(q)$ come from the well-known ``jump relations''
\cite{martinsson2019fast} of the double-layer potential (see also
Appendix~\ref{app:sing_val}).
Note that the limit process of approaching $x$ to the boundary
requires special attention to deal with the problem of
singularities in both $G(p,x)$ and $F(p,x)$.  (see
Appendix~\ref{app:sing_val} for details).
Subtracting these two equations from one another and adding them to one
another results in the two equations
 \begin{equation}
 \begin{split}
 u^{*}_{+}(q) - u^{*}_{-}(q) = \mu(q), \\
 \frac{1}{2}[u^{*}_{+}(q) + u^{*}_{-}(q)] = \int_S [G(p,q)\sigma(p) + F(p,q)\mu(p)] dS(p).
 \end{split}
 \end{equation}
We thus have two equations, which we can solve for the two unknown density functions
$\sigma(p)$ and $\mu(p)$. In fact, we see that the value of $\mu(p)$ is given
explicitly as the jump $u^{*}_{+}(q) - u^{*}_{-}(q)$.

When BEM is used, the double boundary
condition can be expressed in matrix form as:
  \begin{equation}
  \frac{1}{2}( \overline{\mathbf{u}}^{*}_{+} +
  \overline{\mathbf{u}}^{*}_{-} )
  =\overline{\mathbf{G}}\overline{\bm{\sigma}} +
  \overline{\mathbf{F}}\,\overline{\bm{\mu}}.
  \label{eq:db_mat_solve_bem}
  \end{equation}
Likewise, when our hybrid method combining BEM and BIEM is used, the double
boundary condition can be expressed as:
  \begin{equation}
  \frac{1}{2}( \quadrature{\mathbf{u}}^{*}_{+} + \quadrature{\mathbf{u}}^{*}_{-} ) =\hybridquad{\mathbf{G}}\quadrature{\bm{\sigma}} + \hybridquad{\mathbf{F}}\quadrature{\bm{\mu}},
  \label{eq:db_mat_solve}
  \end{equation}
where $\hybridquad{\mathbf{F}}=\hybridint{\mathbf{F}} \,
\integral{\mathbf{P}}  \quadrature{\mathbf{P}}^{-1}$, and
$\quadrature{\bm{\mu}}=\quadrature{\mathbf{u}}^{*}_{+}-\quadrature{\mathbf{u}}^{*}_{-}$.

This double-sided boundary condition is precisely the condition
required to specify the colors on each side of a diffusion curve.

Fig.~\ref{fig:bvp_comparision_curve} shows a single diffusion curve example
solved with different methods. 
We used a discretization size of $s=8$ for BEM (b), $g=e=8$ for BIEM (c), 
$g=8,e=40$ for BIEM (d), and $g=8,s=e=40$ for our Hybrid method (e).
Note the visible polyline on BEM, and the 
dotted-looking artifact on BIEM even for large discretization size (d),
whereas our Hybrid method is free from both limitations.

\begin{figure}
  \includegraphics[width=\linewidth]{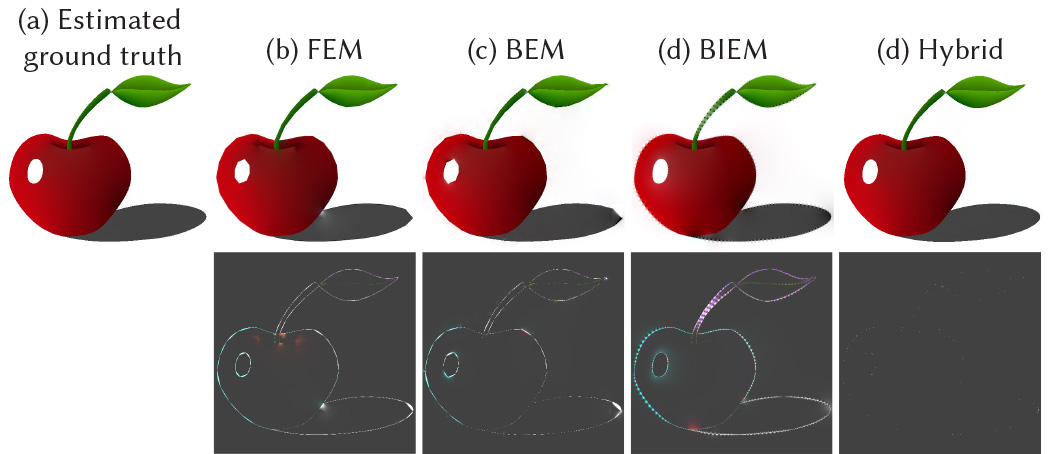}
  \caption{
  Diffusion curve comparison with FEM, BEM, BIEM and our Hybrid method (top),
  and its error (bottom) between estimated ground truth, which was achieved
  by running a dense discretization of BEM ($s=100$ line segments per curve).
  }
\label{fig:compare_methods}
\end{figure}

Fig.~\ref{fig:compare_methods} shows a diffusion curve example of a cherry,
with a comparison between FEM, BEM, BIEM and our Hybrid method.  We used a
discretization size of $s=4$ for BEM and FEM,  $g=4,e=20$ for BIEM, and
$g=4,s=e=20$ for our Hybrid method.

 \subsection{Neumann Boundary Condition}
\label{sec:neumann_bc}
We can also formulate and solve boundary integral equations for Neumann
boundary conditions.  When a domain $V$ is bounded by a simple closed curve
$S$, we can construct a BIE for a Neumann boundary condition on $S$, as
follows.  Using Green's third identity, we have the integral representation
\begin{equation}
u(x) = \int_S [G(p,x)\psi(p)+F(p,x)u(p)]dS(p), \forall x \in V,
\label{eq:int_V_neumann}
\end{equation}
where $\psi(p)=\frac{\partial u(p)}{\partial n}$ denotes the Neumann
boundary values.
Assuming that the Neumann boundary condition is given on the boundary as
$\psi=\psi^{*}$, we let $x$ approach the boundary $S$, and obtain the
following BIE:
\begin{equation}
u(q) = \int_S [G(p,q)\psi^*(p)+F(p,q)u(p)]dS(p)+\frac{1}{2}u(q), \forall q \in S \\
\label{eq:int_S_dbnd}
\end{equation}
We can solve this equation for the unknown Dirichlet value on boundary.

In general, if we are given a boundary condition which 
specifies a combination of Dirichlet boundary values $u^*_d$ on some parts of
the boundary,
and Neumann boundary values $\psi^*_n$ on other parts, 
then we can solve the following matrix system for $u_n$ and $\psi_d$:
\begin{equation}
\frac{1}{2}
\begin{bmatrix}
\quadrature{\mathbf{u}}_{d}^* \\ \quadrature{\mathbf{u}}_{n}
\end{bmatrix}=
\begin{bmatrix}
(\hybridquad{\mathbf{G}})_{dd} & (\hybridquad{\mathbf{G}})_{dn} \\
(\hybridquad{\mathbf{G}})_{nd} & (\hybridquad{\mathbf{G}})_{nn} 
\end{bmatrix}
\begin{bmatrix}
\quadrature{\bm{\psi}}_{d} \\ \quadrature{\bm{\psi}}_{n}^{*}
\end{bmatrix}+
\begin{bmatrix}
(\hybridquad{\mathbf{F}})_{dd} & (\hybridquad{\mathbf{F}})_{dn} \\
(\hybridquad{\mathbf{F}})_{nd} & (\hybridquad{\mathbf{F}})_{nn} 
\end{bmatrix}
\begin{bmatrix}
\quadrature{\mathbf{u}}_{d}^* \\ \quadrature{\mathbf{u}}_{n}
\end{bmatrix},
\label{eq:neumann_mat}
\end{equation}
where the subscripts $d$,$n$ denote the parts of the boundary on which Dirichlet
and Neumann boundary conditions are specificed, respectively.
Fig.~\ref{fig:neumann} shows the resulting image when a combination of Dirichlet
boundary conditions and Neumann boundary conditions are specified.  

Double-sided Neumann boundary conditions on open curves also admit BIE
formulations, and can similarly be handled by our Hybrid Method, with the
key difference that the BIEs involving double-sided Neumann boundary
conditions require the further introduction of an additional hypersingular
kernel $H(p,q)$ (see Appendix~A.1 of~\cite{liu2009fast}).

\subsection{Shortcomings of Brute Force Evaluation}

The discussion up to this point provides us with a new method to solve
for diffusion curves, one which outperforms the standard BEM in accuracy
as demonstrated in Fig.~\ref{fig:bvp_compare_accuracy}, and also shows much 
better performance (see Table~\ref{tab:timing_hybrid}) because the system matrix
to be solved becomes much smaller.

\begin{figure}[h]
  \includegraphics[width=\linewidth]{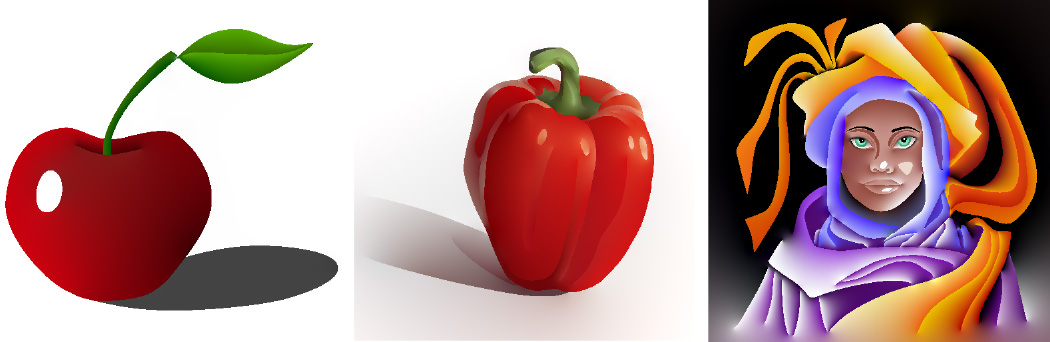}
  \caption{
   Diffusion curve results generated with our Hybrid Method. 
  All images are in a resolution of $512 \times 512$.  
  Examples were taken from \cite{orzan2008diffusion}.}
\label{fig:results_hybrid}
\end{figure} 
% !TEX root = winding.tex
\begin{table}[h!]
\caption{Computation time comparison between BEM solve and Hybrid solve for 
Fig.~\ref{fig:results_hybrid} examples.
Once solved, the evaluation step becomes identical if 
the number of segments are set to be equal.   }

\centering
\ra{1.2}
%\vspace*{-0.1cm}
\setlength{\tabcolsep}{5.5pt}
\rowcolors{2}{lightbluishgrey}{white}
\begin{tabular}{l r r r r r r r r r r r r r r r r r r}
\toprule
\rowcolor{white}
	&  & BEM & Hybrid & \\
\rowcolor{white}
   & curves & solve & solve &eval \\
\midrule
  cherry   												& 32    & 0.10s & 0.008s & 13.9s  \\%0.00287s \\
  red pepper                    						& 109    & 1.47s  & 0.079s & 64.7s \\%0.383s   \\
  person with purple cloak                         & 326    & 32.7s & 0.831s & 326.7s \\%0.124s   \\
\bottomrule
\end{tabular}

\label{tab:timing_hybrid}
  % \vspace*{-0.9cm}
\end{table}

However, as shown in Table~\ref{tab:timing_hybrid}, brute force computation with our
Hybrid method still requires an extremely heavy calculation (especially 
for the evaluation stage,  because the number of pixels is much larger).
 We see then that it is essential to use a fast summation method, such as the Fast
Multipole Method (FMM), to achieve reasonable rendering speeds.

\section{Fast Solution and Evaluation}

The Fast Multipole Method (FMM) is a technique which can be used to rapidly
evaluate the potential induced by a collection of $N$ source curves $S$ at
a set of $M$ targets $q$:
  \begin{align}
u^q = \int_{S} G(p,q) \sigma(p) dS(p).
\label{eq:naive_eval}
  \end{align}
When these potentials are evaluated naively by brute force, the cost grows
as $\mathcal{O}(NM)$. With the Fast Multipole Method, this cost is reduced
to $\mathcal{O}(KN+KM)$, where the factor $K$ grows logarithmically with the
desired accuracy. The key idea behind the FMM is the observation that the
potential induced by a collection of sources, when evaluated at a
well-separated target, can be represented to high accuracy by a expansion
containing only $K$ terms, where the number of terms is \textit{independent}
of the complexity of the source distribution.

We include a complete  description of
the Fast Multipole Method we implemented to accommodate our use cases
in Appendix~\ref{sec:fmmtop}.
We consider this a reproducibility contribution to the graphics community.
Readers unfamiliar with the Fast Multipole Method are \emph{strongly encouraged}
to first read our appendix. FMM is a fairly complicated method, and though
many books and previous descriptions exist, we have made a special effort to
write a self-contained introduction in terms graphics readers will hopefully
better understand.
Nevertheless, readers
who are willing to treat the FMM as a black box can jump right into
Sec.~\ref{sec:adapinf} with no loss of continuity.
We denote the  evaluation of the single-layer integral
operator using the FMM by:
  \begin{align}
u^q = \mathbf{FMM}_{G}(S,\sigma,q),
\label{eq:fmm_eval}
  \end{align}
where $S$ is a collection of source curves, $\sigma$ is a density, and $q$ is
a collection of target points.

\subsection{Non-Uniform Quadtree}

\begin{figure}[h]
 \includegraphics[width=\linewidth]{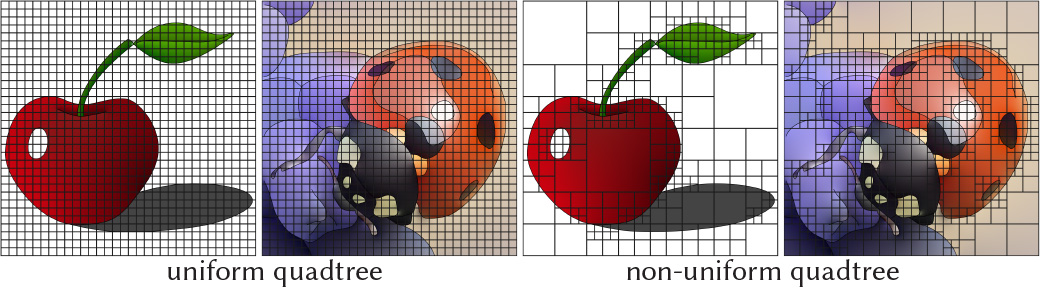}
\caption{Comparison of a uniform quadtree and a non-uniform quadtree}
\label{fig:quadtree_comparison}
\end{figure}

Diffusion curves with the FMM using a uniform quadtree \cite{sun2014fast} 
becomes inefficient for large domains as it requires a lot of memory for cell
allocation as well as significant computation time. 
Figure~\ref{fig:quadtree_comparison} compares the use of
uniform (perfect) and non-uniform (sparse) quadtrees.
We can see that the non-uniform quadtree
is much more efficient, allocating a denser quadtree only in regions 
requiring it.  

Please refer to \ref{sec:non_uniform_quadtree} for a detailed discussion of the
technical differences between uniform and non-uniform quadtrees.

\subsection{Quadtree Clipping}

In our discussion of the Fast Multipole Method, we assumed that the curves
$S$ were discretized with $M$ degrees of freedom. Suppose that $S$ is
discretized with $M$ line segments, and that the integrals over $S$ in the
\textit{target-from-source}, \textit{outgoing-from-source}, and
\textit{incoming-from-source} formulas are computed using BEM, as
described in Section~\ref{sec:multipole_exp}.

At the quadtree construction stage, all cells are subdivided until each cell
contains fewer than $b$ degrees of freedom, which in this case means fewer
than $b$ line segments. Since the degree of freedom associated with a BEM
line segment is located at its midpoint, there will be line segments which
span multiple leaf cells in the quadtree, but since the midpoint of such a
line segment only belongs to a single cell, this segment is handled by only
a single one of the \textit{target-from-source},
\textit{outgoing-from-source}, and \textit{incoming-from-source} formulas.
This can result in invalid computations, if, e.g., a part of the line
segment that should be handled by the \textit{target-from-source} formula is
handled by the \textit{incoming-from-source} formula. Such a situation can
occur when computing the \textit{incoming-from-source} terms from the
\textit{bigger separated list}, when a line segment with its midpoint in a
cell belonging to the \textit{bigger separated list} has an endpoint in a
cell close to or even inside the target cell. The part of the curve near the
endpoint should be computed using the \textit{target-from-source} formula,
but will instead be computed incorrectly by the
\textit{incoming-from-source} formula. In practice, such a situation will
indeed occasionally occur, since we are allowing for dramatically different
scales and curve sizes in our problem.

\begin{wrapfigure}[6]{r}{1.0in}
  \includegraphics[width=\linewidth,trim={6mm 0mm 0mm 4mm}]{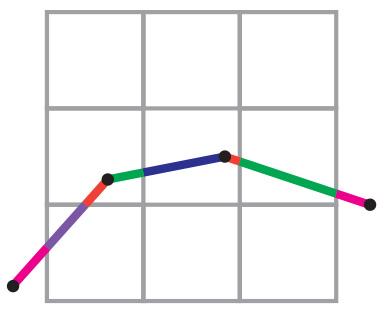}
  \label{fig:quad_clip}
\end{wrapfigure} 
This situation is completely remedied by clipping each BEM segment into
multiple segments using the quadtree, so that each resulting segment is
contained entirely within a single leaf cell (see inset). This ensures that
each part of a BEM segment is handled by the correct formula.

\subsection{FMM for normal derivative of Green's function}

In previous sections,  we described the Fast Multipole Method for rapidly
evaluating the single-layer potentials Eq.~\ref{eq:int_V}.  Suppose instead
that we would like to evaluate the double-layer potential with density $\mu$
over the collection of all curves $S$ at all query points $q$:
  \begin{align}
u^q = \int_{S} F(p,q) \mu(p)
dS(p).
  \end{align}
It turns out that the Fast Multipole Method for evaluating double-layer
potentials is essentially identical to the one presented for single-layer
potentials. The only parts of the method that are changed are the
\textit{target-from-source} (see Appendix~\ref{app:analytic_int_F} and
~\ref{app:sing_int_F}), \textit{outgoing-from-source} (see Appendix~\ref{app:int_out_exp_f}), 
and \textit{incoming-from-source} (see Appendix~\ref{app:int_inc_exp_f})
operators. We denote this FMM evaluation of the 
double-layer integral operator by:
  \begin{align}
u^q = \mathbf{FMM}_{F}(S,\mu,q).
  \end{align}

\subsection{Precomputations}
\label{sec:pre_fmm}

When the Fast Multipole Method is used to evaluate the potential produced by
several different density functions $\sigma$ over a single discretization
$\overline{S}$ of a set of curves $S$ and a single set of target points $q$, a
large number of computations can be reused between evaluations. Such
quantities are independent of the density function, and can be precomputed
and stored, in order to accelerate the evaluation of the FMM. The first
quantity that can be stored is the quadtree over the discretized source
curves $\overline{S}$, which we denote by $\mathcal{Q}$. 
%Given a discretization $\overline{S}$ of the source curves $S$, 
We denote the function constructing the quadtree over that discretization by
  \begin{align}
\mathcal{Q} = \mathbf{quadtree}(\overline{S}).
  \end{align}
The next set of quantities which can be precomputed are the various terms
that appear in the operators used by the Fast Multipole Method. We denote
the collection of precomputed quantities associated with these operators by
$\mathcal{P}_G,\mathcal{P}_F$, corresponding to the single-layer and
double-layer FMM respectively, and denote the function constructing these
quantities by
  \begin{align}
\mathcal{P}_G,\mathcal{P}_F = \mathbf{pre\_FMM}(\overline{S},q,\mathcal{Q}).
  \end{align}
We describe the precise quantities which are precomputed for each one of the
FMM operators in Appendix~\ref{app:pre_fmm}. 

When these precomputed quantities are available, we can accelerate the FMM
by skipping the associated calculations. We indicate that precomputed
quantities are used by providing the precomputed quantities as additional
arguments to the FMM, writing
  \begin{align}
u^q = \mathbf{FMM}_{G}(S,\sigma,q,\mathcal{Q},\mathcal{P}_G).
  \end{align}
%
%%%If only the precomputed quadtree is used, we write
%%%%
%%%  \begin{align}
%%%u^q = \mathbf{FMM}_{G}(S,\sigma,q,\mathcal{Q}).
%%%  \end{align}
%%%%

\subsection{BEM + FMM}

In this section, we describe how to combine the BEM, described in
Section~\ref{sec:bem}, with the FMM, described in Section~\ref{sec:fmm}, in
order to rapidly solve and render diffusion curves.

\subsubsection{BEM + FMM for solving for unknown density}
\label{sec:fmm_solve}

Suppose that we would like to evaluate a single-layer potential $\sigma$ on
a collection of curves $S$ using the boundary element method. We denote the
density and curves, discretized into $M$ boundary elements, by
$\overline{\bm\sigma}$ and $\overline{S}$, respectively. Using the BEM described in
Section~\ref{sec:bem},  we can evaluate the \textit{target-from-source}
Eq.~\ref{eq:t-f-s}, \textit{outgoing-from-source} Eq.~\ref{eq:o-f-s},  and
\textit{incoming-from-source} Eq.~\ref{eq:i-f-s} operators using the
discretized density $\overline{\bm\sigma}$. If the target points
$\overline{{q}}$ are chosen to be same as the source points (the
midpoints of the BEM segments, also called the collocation points), then the
FMM 
  \begin{align}
\overline{\bm{u}} =
\mathbf{FMM}_G(\overline{S},\overline{\bm{\sigma}},\overline{{q}})
  \end{align}
provides an $\mathcal{O}(KM)$ algorithm (recalling that $K$ is the number of
terms in each expansion) for evaluating the matrix-vector product
  \begin{align}
\overline{\bm{u}} = \integral{\bm{G}}\overline{\bm{\sigma}},
  \end{align}
described in Section~\ref{sec:bem}.

If we are given a desired potential $\overline{\bm{u}}^*$ at the BEM
collocation points, 
then we can solve Eq.~\ref{eq:bem_mat_S} for the unknown density
$\overline{\bm\sigma}$ by a direct solver for linear systems, which will
have cost $\mathcal{O}(M^3)$, which is usually prohibitively large. To use
the FMM given by Eq.~\ref{eq:fmm_eval} to solve the linear system, we must
use a so-called iterative method, which requires only a fast method for
evaluating the product of the matrix $\overline{\bm G}$ with a vector. If the
number of iterations required by the iterative method is small, then the
cost will be proportional to the cost of the~FMM, $\mathcal{O}(M)$.

One such iterative method is the Generalized Minimum Residuals Method, or
GMRES, which solves a linear system $\bm A\bm x^*=\bm y^*$ for a possibly nonsymmetric
matrix $\bm A$, and which minimizes the residual $\|{\bm A\bm x-\bm y^*}\|$, where $\bm x
\approx \bm x^*$ is the approximate solution computed by GMRES.
To indicate that the GMRES method is used to solve the linear system 
$\bm A\bm x=\bm y$ to within an error of $\epsilon$ in the residual, where $f(\cdot)$ is
a function approximating the matrix-vector product $f(\bm x) \approx \bm A\bm x$, and
where $\bm x_0$ is the initial guess for the solution, we write
  \begin{align}
\bm x = \mathbf{GMRES}(f(\cdot),\bm x_0,\bm y,\epsilon).
  \end{align}

Thus, to solve for the unknown density $\overline{\bm\sigma}$ in
Eq.~\ref{eq:bem_mat_S}  using the FMM, we compute
  \begin{align}
\overline{\bm\sigma} =
\mathbf{GMRES}(\mathbf{FMM}_G(\overline{S},(\cdot),\overline{q}),
\overline{\bm \sigma}_0,\overline{\bm u}^*,\epsilon),
  \end{align}
using a random initial guess $\overline{\bm\sigma}_0$. 
(we typically choose $\overline{\bm\sigma}_0=\overrightarrow{1}$)

\subsubsection{BEM + FMM for double-sided boundary condition}
\label{sec:bem_fmm_double_bc}

If a double sided boundary condition is given, like the one described in
Section~\ref{sec:double_bc},  we need to solve for unknown density
$\overline{\bm\sigma}$ in Eq.~\ref{eq:db_mat_solve_bem}.  In other words,
given a discretized curve $\overline{S}$ and boundary conditions
$\overline{\bm u}^*_+$ and $\overline{\bm u}^*_-$ on each side, we must solve for
$\overline{\bm\sigma}$ in
\begin{align}
\frac{1}{2}( \overline{\mathbf{u}}^{*}_{+} +
\overline{\mathbf{u}}^{*}_{-} )
=\overline{\mathbf{G}}\overline{\bm{\sigma}} +
\overline{\mathbf{F}}\,\overline{\bm{\mu}},
\end{align}
where $\overline{\bm\mu}=\overline{\bm u}^*_+ - \overline{\bm u}^*_-$.

Using the FMM, we can rapidly solve for $\overline{\bm\sigma}$, as follows. 
First, we compute a right hand side vector $\bm{b}$ by the computation
  \begin{align}
\mathbf{b} = \frac{1}{2}
( \overline{\mathbf{u}}^{*}_{+} +
\overline{\mathbf{u}}^{*}_{-} )
- \mathbf{FMM}_F(\overline{S},\overline{\bm\mu},\overline{q}).
  \end{align}
Next, we solve for $\overline{\bm\sigma}$ using GMRES:
  \begin{align}
\overline{\bm\sigma} =
\mathbf{GMRES}(\mathbf{FMM}_G(\overline{S},(\cdot),\overline{q}),
\overline{\bm\sigma}_0,\mathbf{b},\epsilon),
  \end{align}
using a random initial guess $\overline{\bm\sigma}_0$. The total cost of this
computation will be $\mathcal{O}(M)$, where $M$ is the number of BEM segments used.

To improve things further, we can precompute the quantities needed by the FMM,
as described in Section~\ref{sec:pre_fmm}. The full algorithm for solving for a
double-sided boundary condition using the FMM and BEM is described in 
Algorithm~\ref{alg:fmm_sol}.

\begin{algorithm}[H]
\caption{FMM + BEM for Solving}
\begin{flushleft}
 \hspace*{\algorithmicindent} \textbf{Input:} source curves $\overline{S}$,
 collocation points $\overline{q}$, boundary values
 $\overline{\bm u}^{*}_{+},\overline{\bm u}^{*}_{-}$, initial guess for density
 $\overline{\bm \sigma}_0$ \\
 \hspace*{\algorithmicindent} \textbf{Output:} density values
 $\overline{\bm\sigma}$, quadtree $\mathcal{Q}$, precomputed values
 $\mathcal{P}_G$, $\mathcal{P}_F$
 \end{flushleft}
\begin{algorithmic}[1]
\State $\mathcal{Q} = \mathbf{quadtree}(\overline{S})$
\State $\mathcal{P}_G,\mathcal{P}_F = \mathbf{pre\_FMM}(\overline{S},\overline{q},\mathcal{Q})$
\State Set $\overline{\bm\mu}=\overline{\bm u}^*_+ - \overline{\bm u}^*_-$
\State $\mathbf{b} = \frac{1}{2}
( \overline{\bm{u}}^{*}_{+} +
\overline{\bm{u}}^{*}_{-} )
- \mathbf{FMM}_F(\overline{S},\overline{\bm\mu},\overline{q},\mathcal{Q},\mathcal{P}_F)$
\State $\overline{\bm\sigma} =
\mathbf{GMRES}(\mathbf{FMM}_G(\overline{S},(\cdot),\overline{q},\mathcal{Q},\mathcal{P}_G),
\overline{\bm\sigma}_0,\mathbf{b},\epsilon)$
\end{algorithmic}
\label{alg:fmm_sol}
\end{algorithm}

%%%Let's denote the above algorithm as simply:
%%%%
%%%\begin{equation}
%%%\bm{\sigma}=\mathbf{solve\_FMM}(\mathbf{p},\mathbf{q},\mathbf{n},\mathbf{c},\bm{\mu}^{*})
%%%\end{equation}

\subsubsection{Diffusion Curve with BEM + FMM}

The overall algorithm for using the BEM and FMM to compute pixel values
$u^q$ at all pixels $q$ on 2D domain, given a set of discretized
diffusion curves $\overline{S}$ with collocation points $\overline{\bm q}$, and
a double-sided boundary condition $\overline{\bm u}^{*}_{+}$ and
$\overline{\bm u}^{*}_{-}$, is as follows:

\begin{algorithm}[H]
\caption{Diffusion Curve with FMM + BEM}
\begin{flushleft}
 \hspace*{\algorithmicindent} \textbf{Input:} source curves $\overline{S}$, collocation
 points $\overline{q}$, pixel targets $q$, boundary values
 $\overline{\bm u}^{*}_{+},\overline{\bm u}^{*}_{-}$, initial guess for density
 $\overline{\bm \sigma}_0$ \\
 \hspace*{\algorithmicindent} \textbf{Output:} target pixel values: $u^q$ 
 \end{flushleft}
\begin{algorithmic}[1]
\State Use Algorithm~\ref{alg:fmm_sol} to solve for the density of the single layer
$\overline{\bm \sigma}$ at the 
collocation points $\overline{q}$, using the inputs $(\overline{S},\overline{q},
\overline{\bm u}^{*}_{+},\overline{\bm u}^{*}_{-})$
\State Set $\overline{\bm \mu}=\overline{\bm u}^*_+ - \overline{\bm u}^*_-$
\State
$u^q=\mathbf{FMM}_G(\overline{S},\overline{\bm \sigma},q,\mathcal{Q},\mathcal{P}_G)
+ \mathbf{FMM}_F(\overline{S},\overline{\bm \mu},q,\mathcal{Q},\mathcal{P}_F)$
\end{algorithmic}
\label{alg:diff_curv_fmm}
\end{algorithm}

\begin{figure}
  \includegraphics[width=\linewidth]{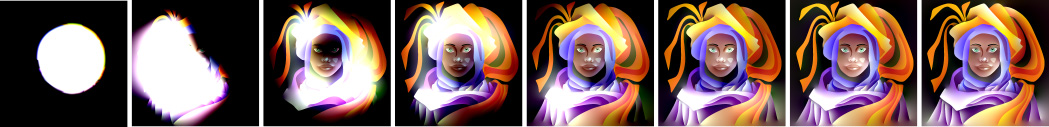}
  \caption{Visualization of the first few iterations of GMRES}
\label{fig:gmres_iter}
\end{figure}

\subsection{Hybrid Method + FMM}

In this section, we describe how to combine our Hybrid Method, described in
Section~\ref{sec:hybrid}, with the FMM, described in Section~\ref{sec:fmm},
in order to rapidly solve and render diffusion curves.

\subsubsection{Hybrid Method + FMM for double-sided boundary condition}

When the FMM is combined with BEM for solving for an unknown density with
a double-sided boundary condition (see
Section~\ref{sec:bem_fmm_double_bc}),  the target points $\overline{q}$ coincide 
with the collocation points on the discretized curves $\overline{S}$:
  \begin{align}
\overline{\bm u} = \mathbf{FMM}_G(\overline{S},\overline{\bm \sigma},\overline{q}).
  \end{align}
On the other hand, when the Hybrid Method is used, the density $\sigma$ is
represented at $g$ Gauss-Legendre nodes $\quadrature{q}$, and this
discretized density is written as $\quadrature{\bm \sigma}$ (see
Section~\ref{sec:hybrid}).
The potential at these quadrature nodes $\quadrature{q}$ is evaluated using
the BEM, where the curve $S$ is descretized into $s$ line segments
$\overline{S}$, with $s > g$. The density $\quadrature{\bm \sigma}$ is
interpolated to the $s$ BEM collocation points by the formula
$\integral{\bm{\sigma}}=\integral{\mathbf{P}}\quadrature{\mathbf{P}}^{-1}
\quadrature{\bm{\sigma}}$. The potential created by the 
density $\quadrature{\bm\sigma}$  can thus be
evaluated at the points $\quadrature{q}$ using the calculation
  \begin{align}
\quadrature{\bm u} =
\mathbf{FMM}_G(\overline{S},\overline{\bm P}\quadrature{\bm P}^{-1}\quadrature{\bm\sigma},\quadrature{q}).
  \end{align}
This is an $\mathcal{O}(s)$ algorithm for evaluating the matrix-vector product
\begin{equation*}
  \quadrature{\mathbf{u}}=\underbrace{\hybridint{\mathbf{G}} \, \integral{\mathbf{P}}  \quadrature{\mathbf{P}}^{-1} }_{\hybridquad{\mathbf{G}}} \quadrature{\bm{\sigma}},
\end{equation*}
where $\hybridint{\mathbf{G}} \in \mathbb{R}^{g \times s}$ and
$\hybridquad{\mathbf{G}} \in \mathbb{R}^{g \times g}$.

The FMM + BEM method for solving for an unknown density with double-sided 
boundary conditions can thus be
reformulated using the Hybrid Method, as follows:

\begin{algorithm}[H]
\caption{FMM + Hybrid Method for Solving}
\begin{flushleft}
 \hspace*{\algorithmicindent} \textbf{Input:} source curves $\overline{S}$,
 quadrature nodes $\quadrature{q}$, boundary values
 $\quadrature{\bm u}^{*}_{+},\quadrature{\bm u}^{*}_{-}$, initial guess for density
 $\quadrature{\bm \sigma}_0$ \\
 \hspace*{\algorithmicindent} \textbf{Output:} density values $\quadrature{\bm \sigma}$,
 quadtree $\mathcal{Q}$, precomputed values $\mathcal{P}_G$, $\mathcal{P}_F$
 \end{flushleft}
\begin{algorithmic}[1]
\State $\mathcal{Q} = \mathbf{quadtree}(\overline{S})$
\State $\mathcal{P}_G,\mathcal{P}_F = \mathbf{pre\_FMM}(\overline{S},\quadrature{q},\mathcal{Q})$
\State Set $\quadrature{\bm\mu}=\quadrature{\bm u}^*_+ - \quadrature{\bm u}^*_-$
\State $\mathbf{b} = \frac{1}{2}
( \quadrature{\mathbf{u}}^{*}_{+} +
\quadrature{\mathbf{u}}^{*}_{-} )
- \mathbf{FMM}_F(\overline{S},\integral{\mathbf{P}}\quadrature{\mathbf{P}}^{-1}\quadrature{\bm \mu},\quadrature{q},\mathcal{Q},\mathcal{P}_F)$
\State $\quadrature{\bm\sigma} =
\mathbf{GMRES}(\mathbf{FMM}_G(\overline{S},\integral{\mathbf{P}}\quadrature{\mathbf{P}}^{-1}(\cdot),\quadrature{q},\mathcal{Q},\mathcal{P}_G),
\quadrature{\bm\sigma}_0,\mathbf{b},\epsilon)$
\end{algorithmic}
\label{alg:fmm_sol_hybrid}
\end{algorithm}

With the Hybrid Method, the dimensionality of the unknown density is much
smaller than for the standard BEM, and the number of BEM segments $s$ used
can also be smaller than for the standard BEM. Hence, the FMM evaluation
step above is much faster.  Furthermore, since the number of degrees of
freedom in the
discretization is smaller, the number of iterations of GMRES is much smaller
as well. This is because we are solving an integral equation involving a
single-layer potential, and the condition number of such a system grows with
the number of degrees of freedom in the discretization.  Fig.~\ref{fig:gmres_iter}
visualizes the first few iterations of GMRES.
%%%The number of iterations needed by GMRES is well-known to
%%%grow with the condition number of the linear system. We can see the performance
%%%gains in the Hybrid Method over BEM in Table~???.

\subsubsection{Diffusion Curve with Hybrid Method + FMM}

The overall algorithm for using the Hybrid Method, together with the FMM,
to compute pixel values
$u^q$ at all pixels $q$ on 2D domain, given a set of discretized
diffusion curves $\overline{S}$ with quadrature points $\quadrature{q}$, and
a double-sided boundary condition $\quadrature{\bm u}^{*}_{+}$ and
$\quadrature{\bm u}^{*}_{-}$, is as follows:

\begin{algorithm}[H]
\caption{Diffusion Curve with FMM + Hybrid Method}
\begin{flushleft}
 \hspace*{\algorithmicindent} \textbf{Input:} source curves $\overline{S}$,
 quadrature
 points $\quadrature{q}$, pixel targets $q$, boundary values
 $\quadrature{\bm u}^{*}_{+},\quadrature{\bm u}^{*}_{-}$, initial guess for density
 $\quadrature{\bm\sigma}_0$ \\
 \hspace*{\algorithmicindent} \textbf{Output:} target pixel values: $u^q$ 
 \end{flushleft}
\begin{algorithmic}[1]
\State Use Algorithm~\ref{alg:fmm_sol_hybrid} to solve for the density of
the single layer $\quadrature{\bm\sigma}$ at the quadrature nodes
$\quadrature{q}$, using the inputs $(\overline{S},\quadrature{q},
\quadrature{\bm u}^{*}_{+},\quadrature{\bm u}^{*}_{-})$
\State Set $\quadrature{\bm\mu}=\quadrature{\bm u}^*_+ - \quadrature{\bm u}^*_-$
\State
$u^q=\mathbf{FMM}_G(\overline{S},\integral{\mathbf{P}}\quadrature{\mathbf{P}}^{-1}\quadrature{\bm\sigma},q,\mathcal{Q},\mathcal{P}_G)
+ \mathbf{FMM}_F(\overline{S},\integral{\mathbf{P}}\quadrature{\mathbf{P}}^{-1}\quadrature{\bm\mu},q,\mathcal{Q},\mathcal{P}_F)$
\end{algorithmic}
\label{alg:diff_curv_fmm_hybrid}
\end{algorithm}

\subsection{Need for Adaptive Subdivision}

\begin{figure}[h]
  \includegraphics[width=\linewidth]{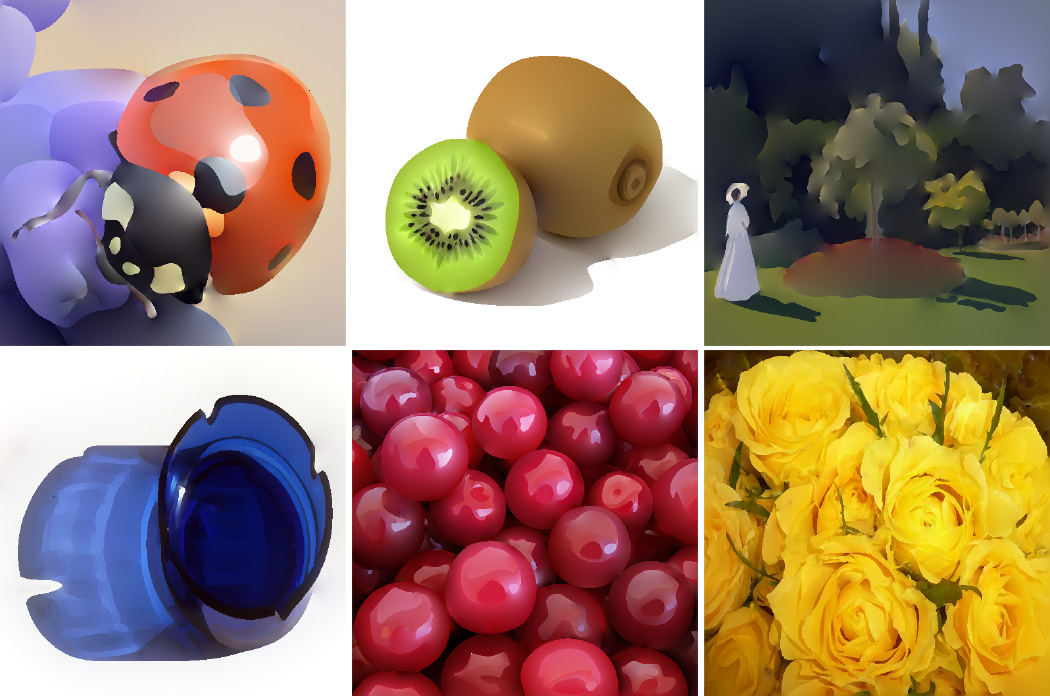}
  \caption{ Diffusion curve results generated with our Hybrid Method + FMM. 
  All images are in a resolution of $512 \times 512$.
Lady bug,  blue glass, and yellow roses examples were taken from
\cite{orzan2008diffusion};  kiwi and Monet examples were taken from
\cite{jeschke2011estimating};  cherries example was taken from
\cite{sun2014fast}.
}   

\label{fig:results_fmm}
\end{figure}

The discussion up to this point provides us with an FMM to solve for
diffusion curves with the Hybrid Method that we introduced in
Sec.~\ref{sec:hybrid}.  Fig.~\ref{fig:results_fmm} shows the 
results of several examples
generated with the Hybrid Method + FMM and its corresponding computation 
time in Table~\ref{tab:timing_fmm}

% !TEX root = winding.tex
\begin{table}[h!]
\caption{Computation time comparison between brute force and FMM for above 
Fig.~\ref{fig:results_fmm} examples
(N\textbackslash A:
cherries and yellow roses examples hung in the evaluation stage).}
\centering
\ra{1.2}
%\vspace*{-0.1cm}
\setlength{\tabcolsep}{5.5pt}
\rowcolors{2}{lightbluishgrey}{white}
\begin{tabular}{l r r r r r r r r r r r r r r r r r r}
\toprule
\rowcolor{white}
	&  & \multicolumn{2}{c}{Brute} & \multicolumn{2}{c}{FMM} \\
\rowcolor{white}
   & curves & solve & eval & solve &eval \\
\midrule
  ladybug   												& 151    	& 0.14s & 133s & 0.33s & 0.41s  \\
  kiwi                    										& 330    & 0.82s  & 355s & 0.59s & 0.81s \\
  monet                         								& 423    & 1.43s & 429s & 0.53s & 0.58s \\
  blue glass                         							& 525    & 2.42s & 519s & 0.53s & 0.68s \\
  cherries                         							& 1110    & 17.26s & N\textbackslash A & 3.55s & 0.95s \\
  yellow roses                         						& 4632    & 1126s &  N\textbackslash A & 10.49s & 1.8s \\
\bottomrule
\end{tabular}
\label{tab:timing_fmm}
  % \vspace*{-0.9cm}
\end{table}

\begin{figure}
  \includegraphics[width=\linewidth]{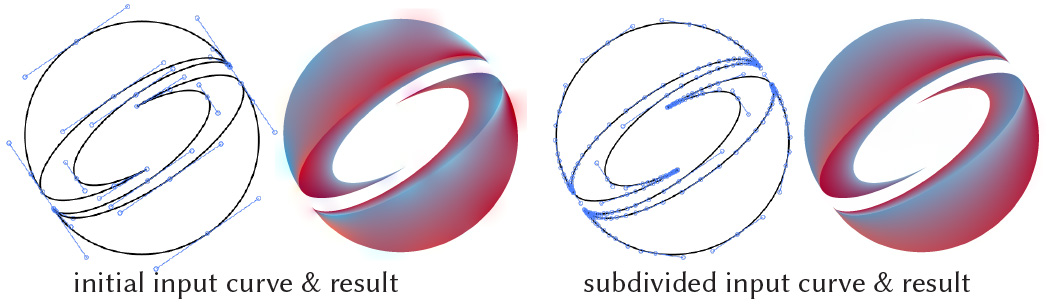}
  \caption{An initial input curve with its resulting color values (left),
  and the same curve subdivided, with its resulting color values (right).}
\label{fig:adap_sub}
\end{figure}

If the given discretization is not sufficient to represent the given
diffusion curve, as shown in left of Fig.~\ref{fig:adap_sub}, or if the user
performs an extreme zoom-in of small region, as demonstrated in
Fig.~\ref{fig:teaser},  we need to employ an adaptive subdivision strategy
to ensure accurate evaluation.

\section{Adaptive Strategy for Infinite Resolution }
  \label{sec:adapinf}

To create an infinite resolution image representation, we adopt the
following adaptive strategy. In our Hybrid Method (see Section~\ref{sec:hybrid}), there
are three discretization parameters that we have control over: $g$, the
number of quadrature nodes; $s$, the number of BEM line segments used to
discretize the source curves during the solution step; and $e$, the number
of BEM segments used to discretize the curves in the evaluation step. 

On the one hand, the number of quadrature nodes $g$ must be large enough so
that the single layer density $\sigma$ solving the double-sided boundary
value problem is well-represented. In other words, we would like $g$ to be
large enough so the Legendre expansions associated with the quadrature nodes
accurately represent $\sigma$. The number of BEM segments $s$ in the solving
step should also be large enough to resolve the potential at all of the
quadrature nodes, and, in practice, this is achieved when $s$ is a constant
multiple of $g$ (we discuss this in the sequel). The number of
discretization nodes $g$ and $s$ thus depend solely on the smoothness of the
solution $\sigma$, and the desired accuracy of approximation. 

On the other hand, the number of BEM segments $e$ used in the evaluation
stage must be large enough so that the layer potentials in the final image
look smooth and continuous. This depends on the particular part of the image
the user is requesting to view, called the \textbf{viewport}, as well as the
pixel resolution that the user is requesting.

The key idea behind our adaptive strategy is to select the optimal number
of quadrature nodes $q$ and BEM segments for the solution stage $s$
\textit{separately} from the number of evaluation BEM segments $e$.  This is
made possible by the interpolation formula described in
Section~\ref{sec:comb_bem_biem}, which
allows us to interpolate a density known at Gauss-Legendre nodes to its
values at arbitrarily many BEM segments. Since the optimal number of
quadrature nodes $q$ will be significantly smaller than the number of
evaluation BEM segments $e$, this will result is a dramatic reduction of
cost in the solution stage. Furthermore, since the density can be interpolated
to any number of BEM segments, this means that it is possible to change the 
viewport and pixel density without re-solving the equations for the density,
which gives us an efficient infinite resolution zoom.

In Section~\ref{sec:adap_sub}, we describe the optimal strategy for
selecting the number of quadrature nodes $q$, and their locations. In
Section~\ref{sec:adap_line}, we describe a strategy for selecting the number
of BEM segments $s$ and $e$. Finally, in Section~\ref{sec:update_viewport}, we describe a
strategy for updating our discretization based on the user's requested
viewport and pixel resolution.

\subsection{Adaptive Subdivision}
\label{sec:adap_sub}

The number of quadrature points needed to represent the density $\sigma$ on
any particular curve will need to be larger if the associated double sided
boundary condition is complicated, if the curve has an intricate geometry
(e.g.  a shape with high curvature), or if the curve has many other curves
nearby.  The relationship between the number of required quadrature nodes
and this a priori information can be somewhat complicated, so instead of using
this a priori information directly, we opt
to use a simpler condition depending only on the density itself.  

If we are able to quickly determine whether or not a particular set of
quadrature points has adequately represented the density, then, by
repeatedly adding more degrees of freedom to any underresolved curve, we can
achieve an optimal discretization of $\sigma$. There are two ways of adding
degrees of freedom to a curve: by choosing a higher-order Gauss-Legendre
quadrature, or by subdividing the curve into two sub-curves, and keeping the
order of the quadrature on each curve the same. Of the two approaches, the
latter is more stable, since increasing the order of the quadrature can
cause the quadrature points to cluster excessively, leading to an increase
in the condition number of the Green's function matrix. In our examples, we
subdivide the curve while keeping the number of quadrature points on each
curve equal to $4$.

To determine whether a panel on a curve needs subdivision, we examine the
size of the highest order expansion coefficient of the Legendre polynomial
expansion of the density.  A large value in the
highest order coefficient means that the density function requires more
degrees of freedom to be properly resolved. Hence, we subdivide the curve if
the highest order coefficient of the Legendre polynomial expansion is larger than
a threshold, which is set globally independent of the domain. This approach
works for any smooth density, however, if the true density is singular (near
a corner, for example), then this strategy can result in an infinite number
subdivisions if we don't provide any constraint on the smallest size of the
resulting curves. Thus, we ensure that further subdivision isn't performed
if the length of the curve is below a second threshold, which is dependent on
the size of the viewport domain.  The overall algorithm can be described as
follows:

\begin{algorithm}
\caption{Adaptive Subdivision}
\begin{algorithmic}[1]
	\Loop \, over every curve:
		\If{length of curve < $\epsilon_1$ (depends on size of pixel) }
		\State skip
		\EndIf
		\If{highest order coefficient of Legendre polynomial expansion > $\epsilon_2$ (global)}
		\State subdivide curve
		\EndIf
	\EndLoop
\end{algorithmic}
\label{alg:adap_sub}
\end{algorithm}

If a collection of curves is determined to need subdivision, then their
density values need to be re-solved.  Fig.~\ref{fig:adap_sub} shows the
initial bleeding artifact being fixed with our adaptive subdivision
algorithm, combined with a re-solving for the density values (note that
subdivision is only applied on a sharp corner of the curve with an abrupt
color change). We accelerate the process of re-solving for the subdivided
densities by the techniques described in the following subsections.

\subsubsection{Warm start of density value}

When re-solving for density values, it is necessary to re-run GMRES with the
FMM (see Section~\ref{sec:fmm_solve}). However, instead of re-solving for
the density values $\quadrature{\bm \sigma}$ from scratch, we can instead
specify an initial guess $\quadrature{\bm \sigma}_0$ from the previous solution,
by interpolating the old density to the quadrature nodes on the subdivided
curves.  This significantly reduces the number of iterations required by
GMRES.

\begin{figure}[h]
  \includegraphics[width=\linewidth]{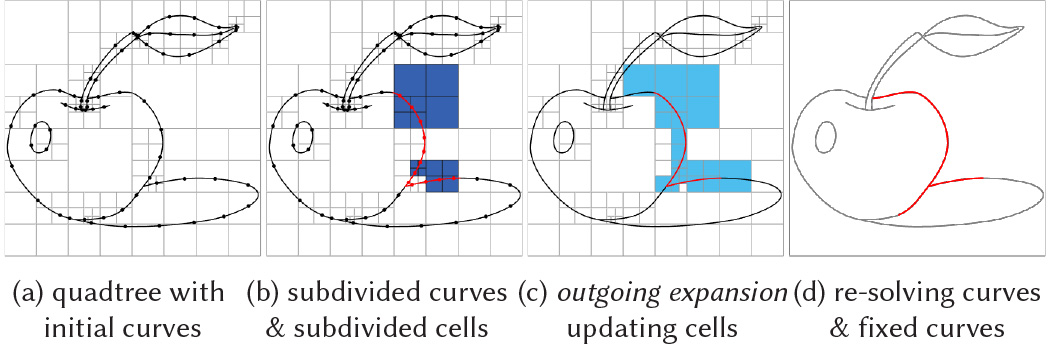}
  \caption{An initial quadtree is generated from an initial curve
  (a). Subdivided curves appear, causing some of cells to further
  subdivide to satisfy the condition on the minimum number of segments per
  cell (b)
  (Section~\ref{sec:update_quadtree}).  For every cell that includes
  subdivided curves, we update its outgoing expansion (c)
  (Section~\ref{sec:update_expansion}).  Re-solving curves are determined by
  perturbing the density value on each subdivided curve (d)
  (Section~\ref{sec:update_density}).  }
\label{fig:adap_quad}
\end{figure}

\subsubsection{Local update of quadtree}
\label{sec:update_quadtree}
Instead of re-constructing the quadtree every time a curve panel is
subdivided, we can update the already pre-constructed quadtree
$\mathcal{Q}$, so that we can re-use the pre-computed expansions $\mathcal{P}_G$
(see Section~\ref{sec:pre_fmm}). Recall that the quadtree is
constructed by repeatedly subdividing cells, until each cell contains fewer
than $b$ degrees of freedom.  To update the quadtree, we check the
prescribed condition on the maximum number of BEM line segments in each
cell and, whenever it is violated after curve subdivision, the containing cell is
subdivided until the prescribed condition is satisfied  (see Fig.
~\ref{fig:adap_quad} (b)). Note that, whenever the quadtree is updated, the
cell relationships described in the various lists (the neighbor list, the
interaction list, etc.) also need to be
updated. These lists can also be locally updated by only
considering newly generated cells or any cells that contain newly generated
cells in their relations. We denote the local update of the quadtree by 
  \begin{align}
\mathcal{Q}' = \mathbf{update\_quadtree}(\overline{S}',\mathcal{Q}),
  \end{align}
where $\mathcal{Q}$ is the quadtree constructed for the old set of discretized
curves $\overline{S}$, and where $\overline{S}'$ and $\overline{\mathcal Q}'$ denote
the updated curves and updated quadtree, respectively.

\subsubsection{Local update of expansions}
\label{sec:update_expansion}
The pre-computed expansion information from the pre-computation step
$\mathcal{P}_G=\mathbf{pre\_FMM}_{G}(\overline{S},\quadrature{q},\mathcal{Q})$
can be also locally updated while preserving some of the previously computed
information. We denote the operator updating the precomputed information by
  \begin{align}
\mathcal{P}_G'=\mathbf{update\_pre\_FMM}_{G}(\overline{S}',\quadrature{q}',
\mathcal{Q}',\mathcal{P}_G),
  \end{align}
where $\overline{S}'$, $\quadrature{q}'$, and $\mathcal{Q}'$ are the updated
discretized curves, quadrature nodes, and quadtree respectively.
%%%

Determining which cells have precomputations that need to be updated depends
on the particular precomputation being considered.  For the
\emph{outgoing-from-source} operator, only the leaf cells that have a change
in their contained line segments need to be updated (see Fig. ~\ref{fig:adap_quad} (c)).  For the
\emph{outgoing-from-outgoing} operator, cells that have a child cell that
has been updated need to be updated accordingly.  For the
\emph{incoming-from-outgoing} operator, the cells that have had a change in
their \emph{interaction list} need to be updated.  For the
\emph{incoming-from-source} operator, the cells that have any change on
their \emph{bigger separated list} need to updated.  For the
\emph{incoming-from-incoming} operator, essentially all cells will need
updating, since the upward and downward process of the FMM will eventually
touch every cell.  Similarly, the \emph{target-from-incoming} operator will
require updates on every cell. For the \emph{target-from-outgoing} operator,
the cells that have any change on their \emph{seperated list} need to be updated,
and for the \emph{target-from-source} operator, those cells that have any
change on their \emph{adjacency list}, or any change on their own cells,
need to be updated. All newly generated cells will also obviously need to
constuct every pre-computation for each of the FMM operators.

\subsubsection{Local re-solve of density value}
\label{sec:update_density}

When some curves are subdivided, the density does not necessarily have to be
updated everywhere. The density on curves with subdivided panels must
naturally be updated, as well as the density on nearby curves which might be
affected by this subdivision. Curves that are far away may not need updating
at all.

We use the following method to determine which curves are in need of
re-solving.  We run the FMM separately on each of the newly subdivided
curves $S_\alpha \in \overline{S}'$, using a density $\overline{\bm\rho}$ that
is given the value one on the subdivided $S_\alpha$, and is given the
value zero on all other curves: 
  \begin{align}
\quadrature{\bm v} = \mathbf{FMM}_G(\overline{S}',
\integral{\bm\rho},\quadrature{q}',
\mathcal{Q}',\mathcal{P}_G')
  \label{eq:resolve_fmm}
  \end{align}
where $\overline{S}'$, $\quadrature{q}'$, $\mathcal{Q}'$, and
$\mathcal{P}_G'$ are the updated discretized curves, quadrature nodes,
quadtree, and precomputations, respectively.
We use the resulting output potential $\quadrature{\bm v}$ to
determine which curves are influenced most by the subdivision.  
Intuitively, this can be
seen as perturbing the density value with a unit charge on each subdivided
curve to check for its effect on other curves. 
We then determine which other curves need to be updated, by the following
procedure. First, we determine how much the induced potential changes over
the subdivided curve $S_\alpha$ by computing the quantity
$\max\quadrature{\bm v}_{S_\alpha} - \min{ \quadrature{\bm v}_{S_\alpha}}$.
We then find all curves $S_\beta$ in need of re-solving by checking if the
change in the 
induced potential $\max\quadrature{\bm v}_{S_\beta} - \min{
\quadrature{\bm v}_{S_\beta}}$ is greater than $0.9$ times the change
over $S_\alpha$. In other words, a curve $S_\beta$ is determined to be in
need of re-solving if
  \begin{align}
\max\quadrature{\bm v}_{S_\beta} - \min{ \quadrature{\bm v}_{S_\beta}}
> 0.9 (\max\quadrature{\bm v}_{S_\alpha} - \min{ \quadrature{\bm v}_{S_\alpha}}).
  \end{align}
We label all such curves as the \textbf{re-solving curves} (see Fig. ~\ref{fig:adap_quad} (d)).  
One of the reasons that this scheme works so well in practice is that the mean-value
theorem tells us that the change in the induced potential over the curve is
related to the integral of the derivative of the potential. 
The potential induced by the kernel $G(p,q)$ is proportional to
$\log(\|p-q\|)$ and decays very slowly, and does not provide a good measure
for identifying nearby re-solving curves.  On the other hand, its derivative
$F(p,q)$ is proportional to $1/\|p-q\|$ and decays much more quickly, and
provides an excellent measure for identifying nearby curves. 

The extra step of potential evaluation for each subdivided curve $S_\alpha$
needed to identify the \textit{re-solving curves} saves a great amount of
computation in the end, because it reduces the dimensionality of the linear
system we must solve for the new density value, which needs to be computed
using several steps of FMM evaluation within GMRES. In fact, it is not
necessary to compute the full FMM in Eq.~\ref{eq:resolve_fmm}, since, as we
describe later in this section, a more efficient local FMM can be performed
instead.

After we have determined the \emph{re-solving curves}, we only need to
re-solve for the density on the \emph{re-solving curves}, while leaving the
density on the remaining curves, which we call the \textbf{constrained
curves}, unchanged.  We denote the set of \textit{re-solving curves} by
$r\subset S$, and the set of \textit{constrained curves} by $c \subset S$.
We then denote the quadrature nodes on the \emph{resolving curves} by
$\quadrature{q}_r'$, and the quadrature nodes on the \emph{constrained
curves} by $\quadrature{q}_c$. Likewise, we denote the BEM segments
on the \emph{re-solving curves} by $\integral{q}_r'$ and on the
\emph{constrained curves} by $\integral{q}_c$ (omitting the $(\cdot)'$ from
the constained curves, since the nodes and BEM segments on
those curves are unchanged).
We represent the above procedure for determining which curves are
\textit{re-solving curves} and which are \textit{constrained curves} by
  \begin{align}
r,c = \mathbf{label\_curves}(\overline{S}',\quadrature{q}',
\mathcal{Q}',\mathcal{P}_G'),
  \end{align}
where $\overline{S}'$, $\quadrature{q}'$, $\mathcal{Q}'$, and
$\mathcal{P}_G'$ are the updated discretized curves, quadrature nodes,
quadtree, and precomputations, respectively.

\subsubsection{Accelerating the local re-solve}

Since the density is unchanged on all of the constrained curves, we can
solve a much smaller linear system than we would if we were solving for the
density from scratch.  Let's first formulate this problem with a matrix
system so that we have clear idea what we want to achieve, and then
reformuate the problem using the FMM.
Recall that, using our Hybrid Method, the unknown density
$\quadrature{\bm\sigma}'$ is the solution to the linear system
  \begin{align}
\hybridquad{\mathbf{G}}\quadrature{\bm\sigma}' = \quadrature{\bm b}^{*'},
  \end{align}
for some right hand side ${\bm b^*}'$, where $\hybridquad{\mathbf{G}} =
\hybridint{\mathbf{G}} \, \integral{\mathbf{P}}
\quadrature{\mathbf{P}}^{-1}\in \R^{g\times g}$, and $\hybridint{\mathbf{G}}
\in \R^{g \times s}$ (see Section~\ref{sec:comb_bem_biem}).  Solving for the potential
only at the quadrature nodes $\quadrature{q}_r'$ on the \emph{re-solving
curves}, while constraining the density at the quadrature nodes
$\quadrature{q}_c$ on the \textit{constrained curves}, gives us the linear
system
\begin{equation}
\begin{bmatrix}
(\hybridquad{\mathbf{G}})_{rr} & (\hybridquad{\mathbf{G}})_{rc} \\
\end{bmatrix}
\begin{bmatrix}
\quadrature{\bm{\sigma}}_r' \\ \quadrature{\bm{\sigma}}_c
\end{bmatrix} = \quadrature{\bm b}^{*'}_{r}
\end{equation}
Then, setting $\quadrature{\bm{\sigma}}_c$ to the previously solved
value and placing it on right hand side: 
\begin{equation}
(\hybridquad{\mathbf{G}})_{rr} \quadrature{\bm{\sigma}}_r' =
\quadrature{\bm b}^{*'}_{r} - (\hybridquad{\mathbf{G}})_{cr}
\quadrature{\bm{\sigma}}_{c}.
\label{eq:local_resolve_mat}
\end{equation}
To solve this equation with the FMM, we compute:
\begin{align}
\quadrature{\bm b}^{\#'}_r =
\quadrature{\bm b}^{*'}_{r} - 
\mathbf{FMM}_G(\overline{S}_c,
\integral{\mathbf{P}}\quadrature{\mathbf{P}}^{-1}\quadrature{\bm \sigma}_c,
\quadrature{q}_r'),
\end{align}
and then solve for $\quadrature{\bm \sigma}_r'$ using GMRES:
  \begin{align}
\quadrature{\bm \sigma}_r' =
\mathbf{GMRES}(\mathbf{FMM}_G(\overline{S}_r',
\integral{\mathbf{P}}\quadrature{\mathbf{P}}^{-1}(\cdot),\quadrature{q}_r'),
(\quadrature{\bm \sigma}_r')_0,{\quadrature{\bm b}^{\#'}}_r,\epsilon),
  \end{align}
where $(\quadrature{\bm\sigma}_r')_0$ is the initial guess for
$\quadrature{\bm\sigma}_r'$ at the updated quadrature nodes $\quadrature{q}'_r$.

While we can use the FMM to compute the potentials created by densities on
the \emph{constrained curves} $S_c$  and the \emph{resolving curves} $S_r'$
from scratch, we can dramatically accelerate the calculations by using the
fact that the $S_c$ and $S_r'$ are subsets of the updated curves $S'$, for
which we already have an updated quadtree $\mathcal{Q}'$ and updated
precomputations $\mathcal{P}_G'$. 
We thus define the following modified FMM, for performing local calculations that 
take advantage of
precomputations on a larger set of curves and targets.

Suppose that $\mathcal{Q}$ and $\mathcal{P}_G$ are the quadtree and
precomputations for the FMM with source curves $S$ and targets $q$.  We
define the FMM
  \begin{align}
u^q = \mathbf{local\_FMM}_{G}(\widetilde{S},\sigma,\widetilde{q},\mathcal{Q},
\mathcal{P}_G)
  \end{align}
for computing the potential created by a subset $\widetilde{S} \subset {S}$
of source curves at a subset $\widetilde{q} \subset q$ of target points, by
modifying the \emph{outgoing-from-source}, \emph{incoming-from-source},
\emph{target-from-incoming}, \emph{target-from-outgoing}, and
\emph{target-from-source} operators (see Section~\ref{sec:multipole_exp}), as follows.  First,
the two operators used to construct outgoing expansions and incoming
expansions from sources are modified to only use the source curves
$\widetilde{S}$.  The last three operators used to compute potentials at the
targets are modified to only compute the potentials at the points
$\widetilde{q}$. Note that the operators for translating expansions can remain
unchanged from Algorithm~\ref{alg:fmm}. 

Finally, using this local FMM, we can find the solution to the system
Eq.~\ref{eq:local_resolve_mat} for the density $\quadrature{\bm\sigma}_r'$ on
the \emph{re-solving curves}, by
first computing
\begin{align}
\quadrature{\bm b}^{\#'}_r =
\quadrature{\bm b}^{*'}_{r} - 
\mathbf{local\_FMM}_G(\overline{S}_c,
\integral{\mathbf{P}}\quadrature{\mathbf{P}}^{-1}\quadrature{\bm\sigma}_c,
\quadrature{q}_r',\mathcal{Q}',\mathcal{P}_G'),
\end{align}
and then solving for $\quadrature{\bm\sigma}_r'$ using GMRES:
  \begin{align}
\quadrature{\bm\sigma}_r' =
\mathbf{GMRES}(\mathbf{local\_FMM}_G(\overline{S}_r',
\integral{\mathbf{P}}\quadrature{\mathbf{P}}^{-1}(\cdot),\quadrature{q}_r',
\mathcal{Q}',\mathcal{P}_G'),
(\quadrature{\bm\sigma}_r')_0,{\quadrature{\bm b}^{\#'}}_r,\epsilon),
  \end{align}
where $(\quadrature{\bm\sigma}_r')_0$ is the initial guess for
$\quadrature{\bm\sigma}_r'$ at the updated quadrature nodes
$\quadrature{q}'_r$.

\subsubsection{Local update for adaptive subdivision}

The overall algorithm for locally re-solving for the density on the
\emph{re-solving curves} can be described as follows: 

\begin{algorithm}[H]
\caption{Local Re-solve after Curve Subdivision}
\begin{flushleft}
 \hspace*{\algorithmicindent} \textbf{Input:} source curves $\overline{S}$,
 quadrature
 nodes $\quadrature{q}$, right hand side $\quadrature{\bm b}^*$, previous
 density $\quadrature{\bm \sigma}$, sources curves after subdivision
 $\overline{S}'$, quadrature nodes after subdivision $\quadrature{q}'$ \\
 \hspace*{\algorithmicindent} \textbf{Output:} updated density value
 $\quadrature{\bm\sigma}'$
\end{flushleft}
\begin{algorithmic}[1]
\State
$\mathcal{Q}' = \mathbf{update\_quadtree}(\overline{S}',\mathcal{Q})$
%%%\Comment{Local update of quadtree}

\State 
$\mathcal{P}_G'=\mathbf{update\_pre\_FMM}_{G}(\overline{S}',\quadrature{q}',
\mathcal{Q}',\mathcal{P}_G)$

\State 
$r,c = \mathbf{label\_curves}(\overline{S}',\quadrature{q}',
\mathcal{Q}',\mathcal{P}_G')$
%%%\Comment{Determining re-solving curves}

\State
$(\quadrature{\bm\sigma}^{'}_r)_0=\mathbf{legendre\_interpolation}(\quadrature{\bm\sigma}_r)$
%%%\Comment{warm start of density value}

\State
$\quadrature{\bm b}^{*'}_r=\mathbf{legendre\_interpolation}(\quadrature{\bm b}^*_r)$
%%%(interpolation of RHS)

\State
$\quadrature{\bm b}^{\#'}_r =
\quadrature{\bm b}^{*'}_{r} - 
\mathbf{local\_FMM}_G(\overline{S}_c,
\integral{\mathbf{P}}\quadrature{\mathbf{P}}^{-1}\quadrature{\bm\sigma}_c,
\quadrature{q}_r',\mathcal{Q}',\mathcal{P}_G')$

\State
$\quadrature{\bm\sigma}_r' =
\mathbf{GMRES}(\mathbf{local\_FMM}_G(\overline{S}_r',
\integral{\mathbf{P}}\quadrature{\mathbf{P}}^{-1}(\cdot),\quadrature{q}_r',
\mathcal{Q}',\mathcal{P}_G'),
(\quadrature{\bm \sigma}_r')_0,{\quadrature{\bm b}^{\#'}}_r,\epsilon)$

\end{algorithmic}
\label{alg:local_update}
\end{algorithm}

\begin{figure}
  \includegraphics[width=\linewidth]{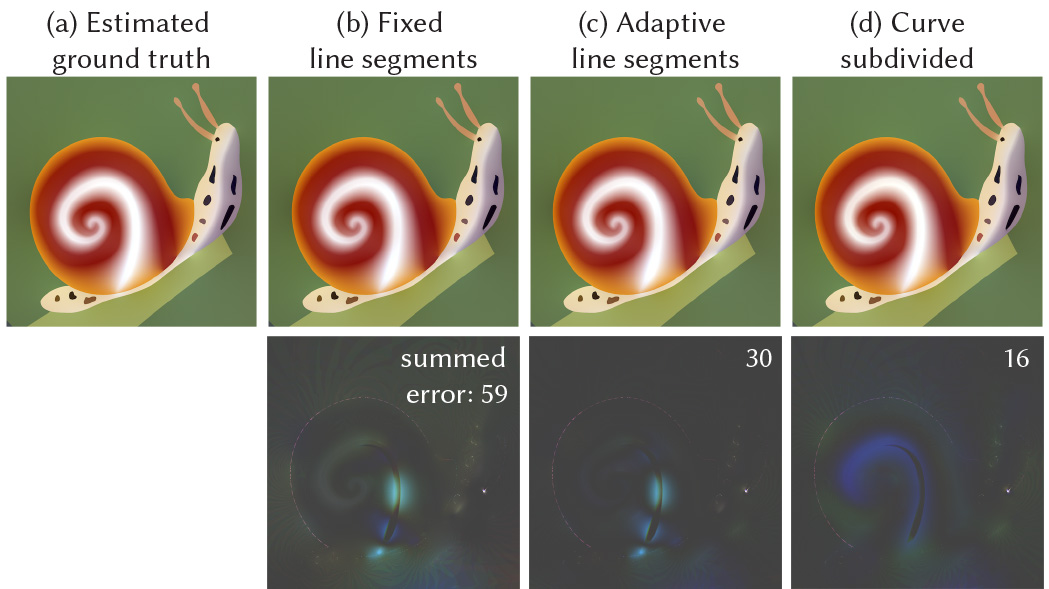}
  \caption{
  Result comparison of snail example with fixed line segments,
  adaptive line segments (Sec.~\ref{sec:adap_line}), and adaptive curve 
  subdivision (Sec.~\ref{sec:adap_sub}) (top) and its error (bottom).
  }
\label{fig:compare_subdivision}
\end{figure}

\subsection{Adaptive BEM Line Segments}
\label{sec:adap_line}

The number of line segments $s$ at the solving stage can be set to be a small
multiple of the number of quadrature points $g$, so if the number of
quadrature points on  each panel of curve is fixed, then the number of line
segments $s$ for the solution stage is fixed as well. In our examples, we
set $s=5g$.

The number of line segments $e$ at the evaluation stage, on the other hand,
needs to be determined to be just fine enough so that user does not perceive
a discretized poly-line, but not so fine as to result in a burdensome
computation. Such a choice of $e$ is best determined by the size of the
pixel and by the user's viewport.

Interestingly, it turns out we do not have to make $e$ depend on the pixel
size directly. If
we set number of the line segments at evaluation stage $e$ to be proportional 
to the arc-length of the curve plus a constant with the following simple equation,
  \begin{align}
e = (\text{length of curve})/10 + s,
  \end{align}
then the result looks perfectly smooth when zooming in, without requiring
an excessive number of calculations. The reason for this is that our adaptive
subdivision algorithm for the density, Algorithm~\ref{sec:adap_sub}, uses the pixel
size as a threshold for subdivision. When the pixel size becomes smaller upon
zooming in, the large curves are subdivided, which causes the total
number of line segments $s$ on that curve to increase.

Note also that this formula ensures that $e \ge s$, since if $e$ is smaller
than $s$, the potential created by the curve will not approximate the
potential we solved for in the solution stage.

\subsection{Diffusion Curve with Adaptive Subdivision}

The overall algorithm for computing a diffusion curve with adaptive
subdivision is described by the following algorithm:

\begin{algorithm}
\caption{Diffusion Curve with Adaptive Subdivision}
\begin{algorithmic}[1]
  \State Solve for the density at the quadrature points using
  Algorithm~\ref{alg:fmm_sol_hybrid} (with a fixed a priori
  discretization)
  \While{ running the adaptive subdivision Algorithm~\ref{alg:adap_sub}
  results in subdivided curves}
  \State Update the density with a local re-solve using
  Algorithm~\ref{alg:local_update} \EndWhile
	\State Use the FMM to evaluate the single and double-layer potentials at 
  the pixel targets on the viewport domain, choosing the BEM line segments
  as described in Section~\ref{sec:adap_line}.
\end{algorithmic}
\label{alg:diff_curv_adap_sub}
\end{algorithm}

Fig.~\ref{fig:compare_subdivision} compares the fixed line segments,
adaptive line segments, and adaptive curve subdivision (top),
showing the resulting error (bottom).
Note that the error has been amplified for clear visualization.

\begin{figure}
  \includegraphics[width=\linewidth]{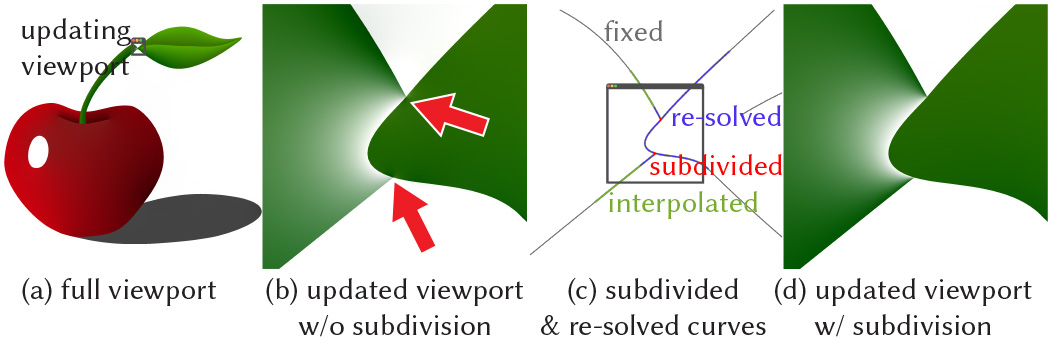}
  \caption{Starting from the full viewport (a), when the user updates the
  viewport with an extreme zoom-in, artifacts are visible without adaptive
  subdivision (b).  Assigning each curve with a suitable label for
  efficiency (c) and updating the viewport with adaptive subdivision gives a
  plausible result (d).  }
\label{fig:updated_viewport}
\end{figure}

\subsection{Updated Viewport}
\label{sec:update_viewport}

Suppose that the user is exploring the domain with their viewport, so that
pixel values need to be re-computed whenever the viewport changes.
Re-evaluating pixel values can be done very quickly with
Eq.~\ref{eq:fmm_eval}, however, whenever the discretization of curves into
BEM line segments must change (if the user zooms in, for example), re-solving
for the density values will require a heavy re-computation if a BEM-only
algorithm is used. This is the moment our hybrid method shines, since we can
simply construct interpolated density values for any re-discretized curves
using Legendre polynomial interpolation as in Eq.~\ref{eq:legendre_poly}.
This process of interpolation will work until the 
viewport domain is zoomed in on such a small region, that the adaptive subdivision
process of Algorithm~\ref{alg:adap_sub} needs further subdivision due to the
smaller threshold resulting from the reduced pixel size (see line 2 from
Algorithm~\ref{alg:adap_sub}).

\begin{figure}[t]
  \includegraphics[width=\linewidth]{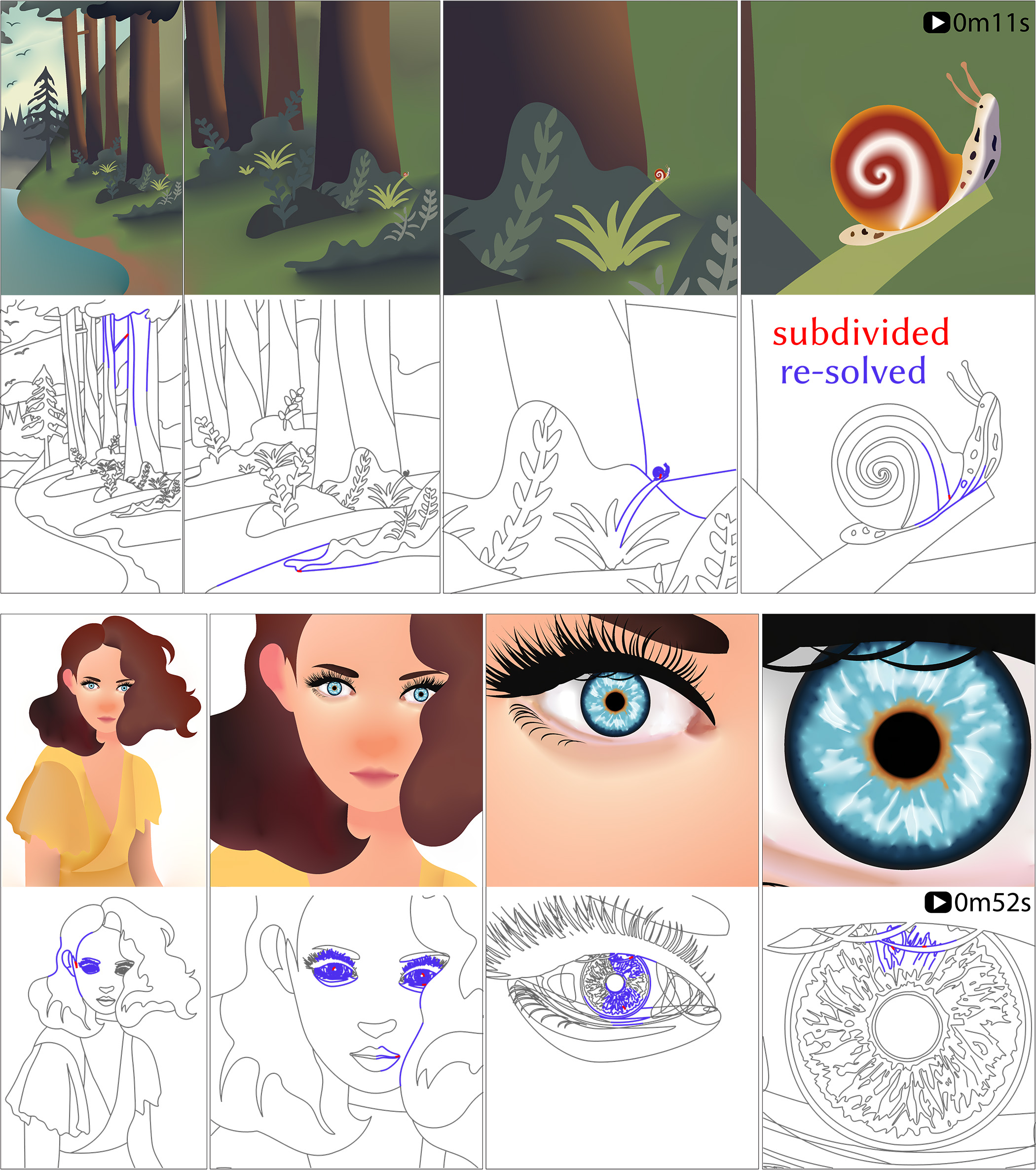}
  \caption{Visualization our of adaptive subdivision procedure while the
  viewport is zooming in.}
\label{fig:adap_viz}
\end{figure}

To account for a change in viewport, we assign to each curve one of the
following three labels: \emph{fixed curve}, \emph{interpolating curve}, or
\emph{re-solving curve} (see Figure~\ref{fig:updated_viewport}). All curves
that are fully outside of the domain of the viewport are labelled as
\textbf{fixed curves}.
Such curves do not require re-discretization, and can retain their BEM
discretization.  The curves for which the density must be re-solved
according to the algorithm described in Sec.~\ref{sec:update_density}
are labelled as \textbf{re-solving curves}, and include both subdivided curves
and neighboring curves.  Finally, all of the remaining curves the intersect
the viewport are labelled as \textbf{interpolating-curves},  and are
re-discretized with smaller BEM line segments using only Legendre
polynomial interpolation, with no need for re-solving. Note that, when
labeling curves to account for the user's viewport, the \emph{constrained
curves} of Section~\ref{sec:update_density} can be either \emph{fixed
curves} or \emph{interpolating curves}, depending on whether or not they
intersect the viewport.  Fig.~\ref{fig:adap_viz} visualizes
subdivided and local re-solving curves determined by its viewport domain
while a user is zooming-in. Table.~\ref{tab:timing_adap} is a comparison
of computation time with global re-solve and with
our local re-solve.

% !TEX root = winding.tex
\begin{table}[t!]
\caption{Computation time comparison between the initial solve of the curve,
the global re-solve, and the local re-solve with zoom-in for the example from
Fig.~\ref{fig:adap_viz}.  }
\centering
%\ra{1.2}
%\vspace*{-0.1cm}
\setlength{\tabcolsep}{5.5pt}
\rowcolors{2}{lightbluishgrey}{white}
\begin{tabular}{l r r r r r r r r r r r r r r r r r r}
\toprule
\rowcolor{white}
	& Initial  & Initial & Global & Local \\
\rowcolor{white}
   & curves & solve & re-solve & re-solve  \\
\midrule
  snail in forest   												& 1493     & 9.48s & 5.35s & 0.95s  \\
  lady with blue eyes                    						& 3126    & 17.47s & 11.29s  & 1.12s  \\
\bottomrule
\end{tabular}
\label{tab:timing_adap}
\end{table}

\section{Anti Aliasing }

Anti aliasing can be smartly handled using the quadtree structure we built
for the FMM.  Instead of sampling color values at the mid point of each
pixel, we can compute a better estimate of the color values by using
area-weighted integration.  This can be easily achieved by recursively
computing pixel values as weighted sums of pixel values from child cells.

Pixels that are inside cells which are bigger than the pixel size
will not benefit from the above strategy.  Since we don't require every
pixel to be assigned smaller child cells, the color strategy for such pixels
reduces to naive multi sampling.  With the key insight that
the aliasing happens
mostly near boundary curves, we modify the quadtree construction condition
so that, if a cell includes any boundary curves, then it will be subdivided
until the leaf cell size becomes smaller then the pixel size. This gives a
very efficient way of anti aliasing.  One limitation of this method is that
the total number of pixels is restricted to be $2^n$, so that the pixels
will always align with the quadtree.

\begin{figure}[h]
  \includegraphics[width=\linewidth]{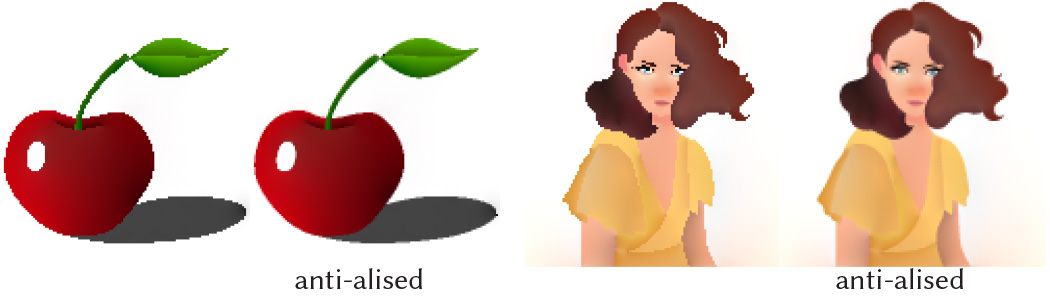}
  \caption{Our method with naive computation of pixel values (left) as well
  as anti-aliased pixel values (right) on a domain with $128 \times 128$ pixels.}
\label{fig:results_fmm_aa}
\end{figure} 

\section{Results}
\label{sec:results}

We demonstrate that our method can accurately compute diffusion curves for a
complex set of input curves with drastically differing scales and sizes of
details, as demonstrated in Fig.~\ref{fig:teaser}. Our algorithm renders the
initial diffusion curves by an adaptive method. The diffusion curves then
retain their accuracy when the viewport is zoomed into the figure, by an
efficient adaptive algorithm that involves re-solving for the density only
on a small subset of curves, as shown in Fig.~\ref{fig:adap_viz}.  Our
adaptive technique is facilitated by our Hybrid Method, which combines the
BIEM and BEM, and allows the density to be accurately interpolated on those
curves appearing in the viewport which are not re-solved, as shown in
Fig.~\ref{fig:updated_viewport}.

\begin{figure}
 \includegraphics[width=\linewidth]{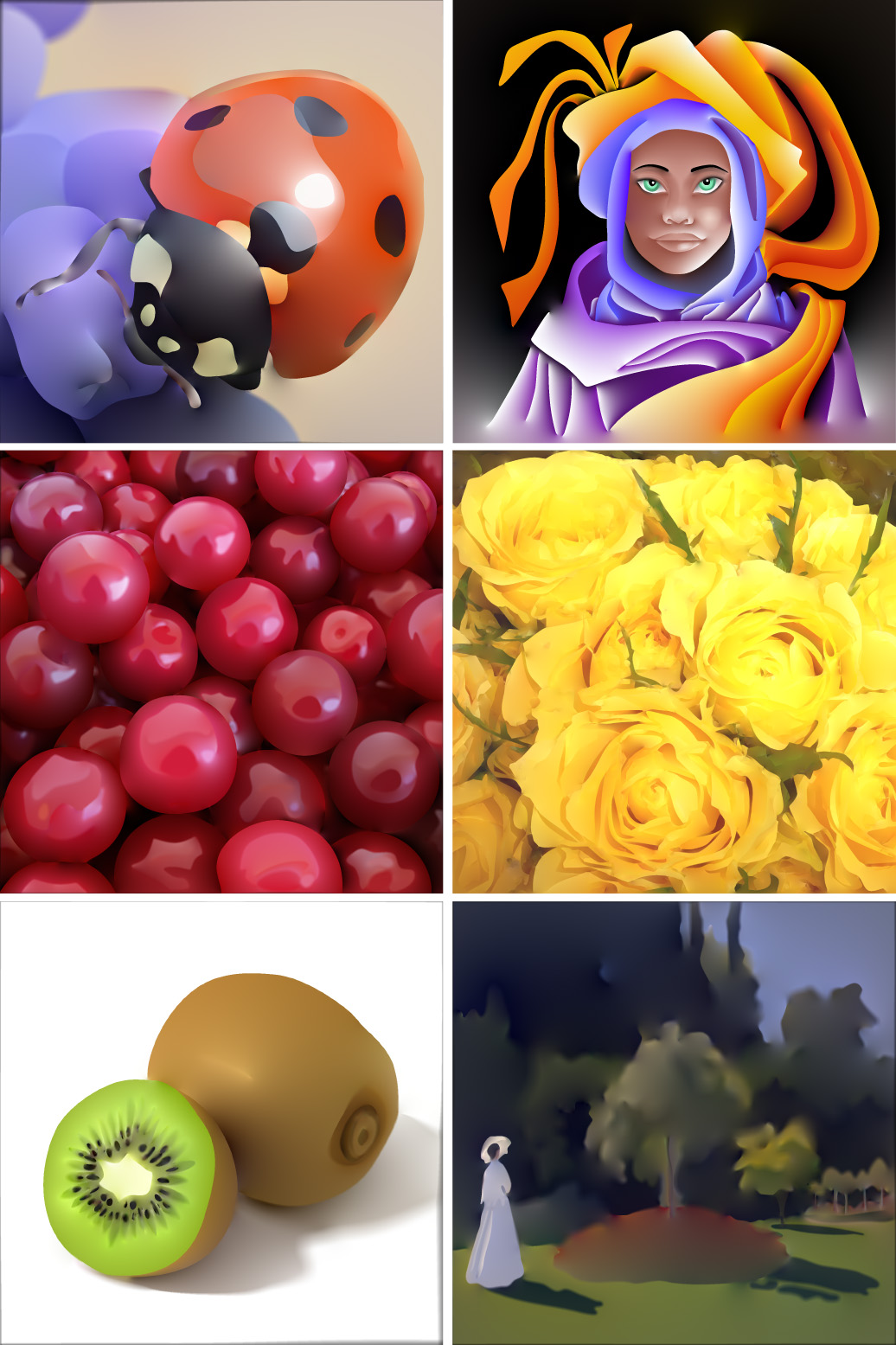}
\caption{
4K resolution results generated with our method.
The blurred scalar field is computed using diffusion curves applied as a
postprocessing step.
This figure is best viewed on a high-resolution digital screen.
Lady bug,  person with purple cloak, and yellow roses examples were taken from
\cite{orzan2008diffusion};  kiwi and  Monet examples were taken from
\cite{jeschke2011estimating};  cherries example was taken from
\cite{sun2014fast}.  }
\label{fig:results_gallery}
\end{figure}

Our method can also be used to generate high resolution images with existing
diffusion curve data from \cite{orzan2008diffusion},  \cite{jeschke2011estimating},
 and \cite{liu2009fast} as shown in Fig.~\ref{fig:results_gallery}.  

We report the computation time in Tables~\ref{tab:timing_hybrid}, ~\ref{tab:timing_fmm},
and~\ref{tab:timing_adap}, but note that our implementation is not fully
optimized, and has a lot of room for faster computation.

\begin{figure}
 \includegraphics[width=\linewidth]{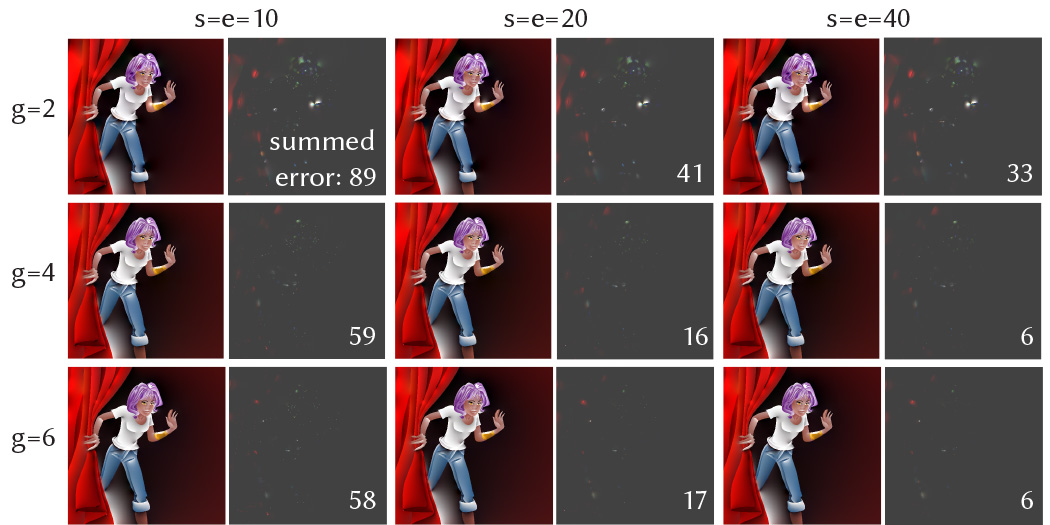}
\caption{
Comparison of different sets of discretization parameters,
with highlighted color difference
between ground truth, and the summed error indicated in white text.
This behind curtain example was taken from \cite{orzan2008diffusion}.}
\label{fig:compare_discretization}
\end{figure}

To demonstrate the effect of setting different values of the number of
solving BEM segments $s$, Gauss-Legendre nodes $g$, and evaluation BEM
segments $e$, we compared the error for different sets of parameters as
visualized in Fig.~\ref{fig:compare_discretization}.  The ground truth here
is derived by running BEM with a very large number of BEM segments $e=200$
for each curve.  As is clear from the figure, the error gets smaller if we
increase the discretization parameters, but the difference is not too large
if the numbers are already high enough.  Empirically,  we found that
$s=20,g=4,e=20$ works best in balancing accuracy and performance.  
Hence, for all the examples in this paper with fixed resolution (in other words,
for all of the examples besides
Fig. ~\ref{fig:teaser} and Fig.~\ref{fig:adap_viz}), we set $s=20,g=4,e=20$,  
except for the example in Fig.~\ref{fig:bvp_compare_accuracy},  where  we set
$s=40,g=8,e=40$ for our Hybrid Method.

\subsection{Implementation}

We implemented the main algorithm of our method in C++ with \cite{libigl},
and additionally used MATLAB and GPTOOLBOX \cite{gptoolbox} for development
and experimentation. 
We note that our implementation is not fully optimized, and has a lot of
room for improvement. Importantly, our code runs almost entirely on the CPU,
and could be accelerated dramatically by an optimized GPU implementation.
All the timings were computed on a MacBook Pro laptop with an Intel 2.4GHz
Quad-Core i9 Processor and 16GB RAM.  Please check the accompanying code to
try out the examples.

\subsection{Comparison with other Methods}

The Finite Difference (FD) method
\cite{orzan2008diffusion,bezerra2010diffusion,finch2011freeform} has its
strength in its easy parallelization and speed when combined with Multigrid
solvers.  The main advantage of our method compared to the FD method, is
that our method correctly constructs diffusion curves around tiny features,
while the FD method flattens down all of its curves to a pixel size, which
loses a lot of detailed information. The resulting pixel-level rasterization
errors are further amplified by the diffusion process.  In Fig.
\ref{fig:compare_FD} (left), for example, all of the important detail is
lost around eye of the woman. Compare this result to Fig.~\ref{fig:teaser}
(left). Furthermore, even when no tiny features are present, the
rasterization process is not stable to small translations, resulting in
strobing artifacts when curves are repositioned.

A solution to this rasterization problem was proposed
in~\cite{jeschke2009gpu}.  There, rather than rasterizing just the boundary
curves, the authors propose rasterizing the entire image by initializing
each pixel to the color of the closest curve point. This initial image is
then smoothed with a Jacobi-like iteration scheme,  in which each pixel is averaged
with four axis-aligned pixels lying on a circle, chosen to be sufficiently
small so as to not intersect any boundary curves. Since averaging over only
axis-aligned pixels instead of the entire circle can create
mach-banding-like artifacts, the stencil size is decreased linearly at each
step, either immediately, or after performing half of the total number of iterations with the stencil at full
size. The iterations on the final stencil are equivalent to classical Jacobi
iterations, with boundary constraints enforced by fixed color data near the
curves, performed on a good initial guess produced by blending the original image rasterization on the larger stencils. Provided an appropriate shrinking
strategy is used, this method can produce a visually excellent diffusion
curve image in real-time, with only slight differences with the fully
converged image.  This method is, however, limited to Dirichlet boundary
conditions, and does not take into account curves outside of the image
domain. Our method can potentially be extended to Neumann boundary
conditions, and handles curves outside of the image domain in a natural way
via the FMM.

The Finite Element Method (FEM)
\cite{takayama2010volumetric,pang2011fast,boye2012vectorial} is the most
widely used method in computer graphics, due to its easy and intuitive
implementation with fast Cholesky solvers. However, FEMs suffer from the
problem of triangulating the domain.  Triangulation of a complex set of
curves is itself very difficult problem, and is the subject of current
research \cite{hu2019triwild}. We have attempted to use TriWild
\cite{hu2019triwild} with the input curves for Fig.~\ref{fig:compare_FD}, but
TriWild discarded all of the highly detailed features, and failed if we
tried to preserve these features. Even when the triangulation succeeds, the
FEM exhibits a bleeding artifact which can be seen in Fig.
\ref{fig:fem_artifact}. 
Using Triangle \cite{shewchuk2005triangle} led to successful triangulation
but resulted in more than 10 million triangles given data from
Fig.~\ref{fig:compare_FD}. 
Even with dense triangulations, FEM still shows
nonsmooth results near curve endpoints, caused by singularities there.
When singularities are present in the PDE solution, 
both increasing only the polynomial order on elements 
of fixed size (called p-refinement) or fixing the polynomial 
order and refining the mesh (called h-refinement) are known to 
provide only modest improvements in solution accuracy. 
However, if a carefully chosen combination of mesh refinement 
and polynomials of varying degree is used, then the FEM can be 
made to converge exponentially fast (this process is called 
hp-refinement or the hp-FEM). Such methods, while often effective, 
can be extremely challenging to implement in a fully automatic fashion
(see~\cite{Gopal_2019}).
A heuristic solution for dealing with singularities, proposed
by~\cite{boye2012vectorial}, is to linearly blend colors around the vertex
of a triangular element lying on a curve endpoint. While visually quite
satisfactory, it is worth noting that this approximation is nonetheless very
different from the true power-type singularities present at such points.
These difficulties, namely the triangulation problem, the bleeding artifact
problem,  and the singularity problem, are all completely absent in our method.
There is no need for triangulation since our method is a boundary-only
method, and the bleeding artifact is resolved completely by adaptive
subdivision, as seen in Fig.~\ref{fig:results_hybrid}.

The pros and cons of the Boundary Element Method (BEM) and Boundary Integral
Equation Method (BIEM) are described in Sec.~\ref{sec:method_bvp}. As 
we discuss in that section, the boundary element method can produce accurate
results when many BEM line segments are used, but results in an extremely large
linear system to solve. On the other hand, the BIEM has a highly efficient 
representation of the solution, and results in a small linear system, but creates
artifacts when the density is evaluated near the curves. Our Hybrid Method
takes the advantages of both, namely, it retains the efficient representation
of BIEM while also obtaining the visual quality of BEM, as shown in
Fig.~\ref{fig:bvp_compare_visual}.

The Walk on Spheres (WoS) method~\cite{sawhney2020monte} has its strength in
its simplicity and its robustness to input data and geometry, but does not
generalize efficiently to Neumann boundary conditions, which is a
significant limitation, since such boundary conditions are so useful in
practice. The two main issues that arise are that WoS can require extremely
long walks, and that the sphere sizes near Neumann boundaries become
very small, significantly impacting performance.  The recently proposed Walk
on Stars (WoSt) method~\cite{sawhney2023walk} overcomes this second issue by
replacing these small spheres with the boundaries of much larger star-shaped
domains.  Nonetheless, when the boundary is predominantly Neumann (like the
boundaries in Fig.~\ref{fig:neumann}), WoSt will still can take very long
walks, since, like WoS, it must reflect back into the domain from Neumann
boundaries and can only terminate at Dirichlet boundaries. This issue can be
somewhat ameliorated by caching solution values on the boundary of the
domain~\cite{miller2023boundary}, but since the solution values must still
be generated by, for example, WoSt, the problem persists.  Our method, on
the other hand, generalizes to Neumann boundary conditions without
difficulty, as shown in Fig~\ref{fig:neumann}.

\section{Limitations and Future Works}

Our proposed method advances diffusion curve representation to a 
high-accuracy level, accommodating fine multiscale features and 
allowing precise zooming and panning while maintaining accuracy.

Despite its
many desirable features, our method still has some limitations and room for
future improvement.

Our Adaptive Strategy was constructed with the assumption that diffusion
curves mostly remain static.  It will not work efficiently for animated
diffusion curves, as it will require large portion of the curve to be
re-computed every step. 

Our requirement of ensuring accurate computations can become burdensome
computationally, because messy or wild curve data will exhibit a lot of
intersecting and overlapping curves, which will require heavy adaptive
subdivision to resolve.  We developed a pre-processing step to deal with
ill-posed curves, but it is difficult to distinguish between an intentional
ill-posed curve placed by an artist and an unintended ill-posed curve. It
will be useful to have a version of our algorithm with softer and less
stringent accuracy requirements, which would allow for wilder curve data.

We demonstrate Neumann boundary conditions in the example of Fig.~\ref{fig:neumann}.
However, our current implementation only supports one-sided Neumann boundary
conditions on closed curves.  These examples were generated by subdividing a
region into disconnected closed subregions.  A more general and powerful
double-sided Neumann boundary condition is possible, but will require the
introduction of a hyper-singular kernel, as described in Sec~2.3 of
\cite{liu2009fast}.

Our current implementation is not fully optimized, and does not reach
real-time computation speeds.  We observed that a
GPU-accelerated implementation of brute force computation leads to
a speedup of more than 100 times.  A GPU-accelerated version of
FMM and adaptive re-computation would put on wings on our method
and make it truly practical. We leave this extension to our future work.

%\section*{Acknowledgement}
\begin{acks}
This work was supported by the National Research Foundation, Korea
(NRF-2020R1A6A3A0303841311), 
the Swiss National Science Foundation's Early Postdoc.Mobility fellowship,  
the NSERC Discovery Grants RGPIN-2020-06022 and DGECR-2020-00356, 
and NSERC Discovery Grant RGPIN-2022-04680, the Ontario Early Research Award program, the Canada Research Chairs Program, a Sloan Research Fellowship, the DSI Catalyst Grant program and gifts by Adobe Inc.

We express deep gratitude to Professor Eitan Grinspun for leading in-depth discussions during the early stages of the research.
We thank Silvia Sell\'an and Otman Benchekroun for their help in conducting experiments and performing proofreading.
\end{acks}

\bibliographystyle{ACM-Reference-Format}
\bibliography{reference}

\begin{appendices}
\section{Arc Length Line Segment Integration}
\label{app:arc_len_int}

Unlike classical BEM, our hybrid method allows for different numbers of BEM
line segments at the solution and evaluation stages. In order for us to use
different discretizations at each stage, it turns out to be essential that
the BEM integration is performed in a way that accounts for the true length
of the curve, since the length of a polyline approximation can change if the
number of segments changes.  Thus, if we naively integrate the line
segments, the result will be different if the number of line
segments is different. Hence, we employ arc length integration to make sure
that different numbers of line segments do not give a different result.

\section{Arc Length Parametrization}
\label{app:arc_len_param}
Unfortunately,  there is no closed form solution for the arc length
parametrization of a cubic B\'ezier curve.  We employ an iterative method to
find the arc length parametrization.  Suppose that we have a cubic Bezier
curve
\begin{align}
B: [0,1] \to \mathbb{R}^2, \\
B(t) = (1-t)^3 C_0 + 3t(1-t)^2 C_1 + 3t^2(1-t) C_2 + t^3 C_3,
\end{align}
where $C_0, C_1, C_2, C_3 \in \mathbb{R}^2$ are the control points.
Then the cumulative length of the Bezier curve up to some parameter $t=s$ 
may be written as an integral:
\begin{equation}
L(s)=\int^s_0 \| B'(t) \| dt.
\end{equation}
We apply Newton's root finding method to find $s$ for a given 
target arc length $T$ such that $T=L(s)$. Given some guess $s_0$, we can
improve this guess by the update formula
\begin{equation}
s \gets s_0 -(L(s)-T)/L'(s).
\end{equation}
However, this requires evaluation of $L(s)$ and $L'(s)$. To compute these
quantities, we use Gauss-Legendre quadrature to approximate a definite
integral over the finite interval $[0,s]$.

\section{Singular Green's Function Integrals}
\label{app:sing_val}

The Green's function $G(p,q)$ and its normal derivative $F(p,q)$ both have
singularities at $p=q$, where the kernel is infinite.
To perform accurate BEM integration involving these kernels, it is important
to carefully consider how these singularities influence the integrals.

Suppose that $x$ is a target point that approaches the boundary.  We divide
the boundary $S$ into two parts: $S-S_{\epsilon}$ and $S_{\epsilon}$, where
$S_{\epsilon}$ is a small segment with arc length $2\epsilon$ centered around
the point to which $x$ will approach.

Let us first consider single layer potential, which involves the integral of
the Green's function $G(p,q)$:
\begin{equation}
\begin{split}
\int_S G(p,x)\sigma(p)dS(p) = \\
\lim_{\epsilon \to 0} \int_{S-S_{\epsilon}} G(p,x)\sigma(p)dS(p) + \lim_{\substack{\epsilon \to 0  \\ x \to S}} \int_{S_{\epsilon}} G(p,x)\sigma(p)dS(p) 
\end{split}
\label{eq:G_sing}
\end{equation}
The limit of the second integral on the right-hand side of Eq.~\ref{eq:G_sing}
turns out to be
\begin{equation}
\lim_{\substack{\epsilon \to 0  \\ x \to S}} \int_{S_{\epsilon}} G(p,x)\sigma(p)dS(p)=0.
\end{equation}
Analytic integration of the limit of the first integral on the right-hand
side of Eq.~\ref{eq:G_sing} over a BEM line segment is described in
Appendix~\ref{app:sing_int_G}.

Let us now consider double layer potential, which involves the integral of
the normal derivative of Green's function, $F(p,q)$: 
\begin{equation}
\begin{split}
\int_S F(p,x)\mu(p)dS(p) = \\
\lim_{\epsilon \to 0} \int_{S-S_{\epsilon}} F(p,x)\mu(p)dS(p) + \lim_{\substack{\epsilon \to 0  \\ x \to S}} \int_{S_{\epsilon}} F(p,x)\mu(p)dS(p) 
\end{split}
\label{eq:F_sing}
\end{equation}
The limit of the second integral on the right-hand side of
Eq.~\ref{eq:F_sing} turns out to be: 
\begin{equation}
\lim_{\substack{\epsilon \to 0  \\ x \to S}} \int_{S_{\epsilon}} F(p,x)\mu(p)dS(p)= \pm \frac{1}{2} \mu(x) , x \in S.
\end{equation}
The sign of the one-half $\mu(x)$ depends on which side of the curve the
point $x$ approaches.

\subsection{Singular Integrals in our Hybrid Method}

For standard BEM, the query points are located on the midpoints of each line
segment.  Thus, for each query point,  integration over the singularity
happens on the one line segment that is located on the query
point, which corresponds to the diagonal entries of $\overline{G}$ and
$\overline{F}$. 

For our Hybrid Method, the query point can be located at any place along
the input curve.
If evaluation points (which are also the quadrature points) are positioned
near the intersection of two BEM line segments, this situation should be
considered as a singular integration for both line segments, and evaluated
by invoking Eq.~\ref{eq:sing_G} and Eq.~\ref{eq:sing_F} in
Appendix~\ref{app:analytic_int}.

\section{Analytic Integration}

In this section, we provide various analytical formulas for evaluating
the integrals of the Green's function kernels $G(p,q)$ and $F(p,q)$ over
the line segment $\overline{p_1 p_2}$.

\label{app:analytic_int}

\subsection{Integration of Green's function}
\label{app:analytic_int_G}
The integral of a Green's function on a line segment can be evaluated
analytically as follows:
\begin{equation}
\begin{split}
-\int_S G(p,q)dS = \int_S \frac{\log(|p-q|)}{2\pi}dp \\
=\frac{|a|}{2\pi} \int_{0}^{1} \log(|b+at|)dt = \frac{|a|}{2\pi} \int_{0}^{1} \log(\sqrt{a^2t^2+2abt+b^2})dt,
\end{split}
\end{equation}
where $a=p_{2}-p_{1}$, $b=p_{1}-q$, with $p_{1}, p_{2}$ being the endpoints
of the line segment. Denoting $a^2=\zeta, 2ab=\eta, b^2=\xi$, we get the
following anti-derivative: 
\begin{equation}
\begin{split}
\frac{|a|}{4\pi} \int_{0}^{1} \log(\zeta t^2+\eta t+\xi) dt = \\
\frac{|a|}{4\pi} \left[ \log(\zeta t^2+\eta t+\xi) (t+\frac{\eta}{\zeta})-\tan^{-1}\left( \frac{2\zeta t + \eta}{\sqrt{4\xi \zeta-\eta^2}} \right)\left(\frac{2\zeta + \eta}{\zeta \sqrt{4\xi \zeta-\eta^2}} \right) \right]_0^1 =\\
\frac{|a|}{4\pi}   \left[ \log(\zeta+\eta+\xi)(1+\frac{\eta}{\zeta})-\log(\xi)\frac{\eta}{\zeta} \right] + \\
\frac{|a|}{4\pi} \left[ \tan^{-1} \left( \frac{\eta}{\sqrt{4\xi \zeta-\eta^2}} \right) - \tan^{-1} \left( \frac{2\zeta+\eta}{\sqrt{4\xi \zeta-\eta^2}}\right) \right]  \left( \frac{2\zeta+\eta}{\sqrt{4\xi \zeta-\eta^2}} \right).
\end{split}
\label{eq:G_analytic}
\end{equation}

\subsection{Integration of Normal Derivative of Green's function}
\label{app:analytic_int_F}
The integral of the normal derivative of the Green's function on a line
segment can be evaluated analytically as follows:
\begin{equation}
\begin{split}
-\int_S F(p,q)dS= \int_S \frac{(p-q) \cdot n}{2\pi |p-q|^2} dp\\
\frac{1}{2\pi} \int_{0}^{1} \frac{(at+b) \cdot n}{|at+b|^2} dt = \frac{b
\cdot a^{\perp}}{2\pi}  \int_{0}^{1} \frac{1}{a^2+2abt+b^2} dt,
\end{split}
\end{equation}
where the numerator in the last integral becomes constant because
$n=a^{\perp}$. Denoting $2b \cdot a^{\perp}=\nu$, we get the
following anti-derivative:
\begin{equation}
\begin{split}
\frac{\nu}{2\pi} \int_{0}^{1} \frac{1}{\zeta t^2 +\eta t +\xi} dt \\
\frac{\nu}{2\pi} \left[ \tan^{-1}\left(\frac{2\zeta t+\eta}{\sqrt{4\xi \zeta-\eta^2}}\right) \left(\frac{\nu}{\sqrt{4\xi \zeta-\eta^2}}\right) \right]_0^1 = \\
\frac{\nu}{2\pi} \left[ \tan^{-1} \left( \frac{2\zeta+\eta}{\sqrt{4\xi \zeta-\eta^2}} \right) \tan^{-1} \left( \frac{\eta}{\sqrt{4\xi \zeta-\eta^2}} \right) \right] \left( \frac{\nu}{\sqrt{4\xi \zeta-\eta^2}}  \right).
\end{split}
\label{eq:F_analytic}
\end{equation}

\subsection{Singular Integration of Green's function}
\label{app:sing_int_G}
Singular integration happens when $q$ lies on the line segment $\overline{p_1
p_2}$. Letting $b=a t^*$, the integral can be expressed simply as: 
\begin{equation}
\begin{split}
\frac{|a|}{2\pi} \int_{0}^{1} \log(|a(t-t^*)|)dt = \\
\frac{|a|}{2\pi} \lim_{\epsilon \to 0} \left( \int_{0}^{t^*-\epsilon} \log(|a(t-t^*)|)dt + \int_{t^*+\epsilon}^{1} \log(|a(t-t^*)|)dt \right) = \\
\frac{|a|}{4\pi}  \lim_{\epsilon \to 0} \left[-2t+(t-t^*)\log(a^2(t-t^*)^2)\right]_{0}^{t^*-\epsilon}  + \\
\frac{|a|}{4\pi}  \lim_{\epsilon \to 0} \left[-2t+(t-t^*)\log(a^2(t-t^*)^2)\right]_{t^*+\epsilon}^{1} =\\
\frac{1}{4\pi} \left( 2|a|-|a|t^* \log(|a|^2 {t^*}^2) +(|a|t^*-|a|) \log(|a|^2(1-t^*)^2) \right)= \\
\frac{1}{2\pi} (s-e_1 \log e_1 - e_2 \log e_2),
\end{split}
\label{eq:sing_G}
\end{equation}
where $s=|a|$, $e_1 = \overline{qp_1}=|a|t^*$, and $e_2 =
\overline{qp_2}=|a|(1-t^*)$.

\subsection{Singular Integration of Normal Derivative of Green's function}
\label{app:sing_int_F}
Similar to the above approach, we can express the singular integral as:
\begin{equation}
\begin{split}
\frac{|a|}{2\pi} \int_{0}^{1} \frac{a(t-t^*) \cdot n}{|a(t-t^*)|^2}dt.
\end{split}
\label{eq:sing_F}
\end{equation}
The above equation becomes $0$ since $(t-t^*) \cdot n=0$.

\section{Fast Multiple Method}
\label{sec:fmmtop}

Consider the evaluation of the BEM integrals in Eq.~\ref{eq:bem_S} over $m$
diffusion curves, for a total of $N=ms$ BEM line segments.  Directly
evaluating the BEM integrals in Eq.~\ref{eq:bem_S} or Eq.~\ref{eq:bem_V} at
$M$ target points requires $\mathcal{O}(NM)$ operations.  Greengard and
Roklin~\cite{greengard1987fast} demonstrated that the task could be done in
$\mathcal{O}(N+M)$ operations in finite precision by introducing the Fast
Multipole Method (FMM), which is asymptotically even better than the
$\mathcal{O}(N\log(N) + M\log(M))$ operations required by the Barnes-Hut
algorithm~\cite{barnes1986hierarchical,pfalzner1997many}.  In this section,
we provide an introduction to the FMM for the complex-valued kernel
$G(p,q)=-\frac{1}{2\pi}\log(q-p)$, whose real part is the Green's function
defined in Eq.~\ref{eq:greens_func}, where $p$ and $q$ are treated as points
in the complex plane.

\subsection{Multipole expansion}
\label{sec:multipole_exp}

The key idea behind the Fast Multipole Method is the observation that, if a
potential is induced by a source involving a Green's function, then the
potential away from the source can be approximated to high accuracy by a
finite sum involving certain basis functions, where the number of terms
involved in the sum is independent of the complexity of the source
distribution. These sums approximating the induced potential are called
\textbf{expansions}.

\subsubsection{Outgoing expansion}
\begin{wrapfigure}[11]{r}{1.2in}
  \includegraphics[width=\linewidth,trim={6mm 0mm 0mm 5mm}]{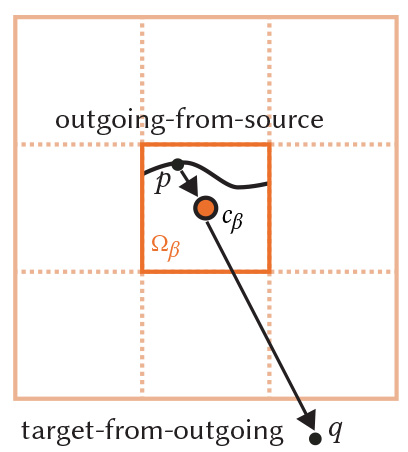}
  \label{fig:ou_exp}
\end{wrapfigure} 
Suppose we have the source curve $S_\beta$ contained in a cell
$\Omega_{\beta}$, and that  we have a query point located far from this
cell, separated by at least one cell width (see inset).  In this case, we
say that a query point is \textbf{well-separated} from cell
$\Omega_{\beta}$. The \textbf{outgoing expansion} is the expansion which
represents the potential induced by this source curve, at all well-separated
target points. The word outgoing reflects the fact that this expansion is directed
outwards from the source curve, and is valid for all target points far away.

Letting $c_\beta$ be the point at the center of the cell $\Omega_\beta$,  
we observe that the Green's function kernel $G(p,q)$ can be written as 
  \begin{align}
&G(p,q)=-\frac{1}{2\pi} \log(q-p) \\
&= -\frac{1}{2\pi} [\log(q - c_\beta) + \log(1 -\frac{p-c_\beta}{q-c_\beta})].
  \end{align} 

Using the Taylor series expansion for $\log$,
\begin{equation}
\log(1-z)=-\sum_{k=1}^{\infty} \frac{z^k}{k} \quad \text{for} \quad |z|<1,
\end{equation}
we see that
  \begin{align}
G(p,q) = -\frac{1}{2\pi} [\log(q-c_\beta) - \sum_{k=1}^\infty
\frac{1}{k} \bigl( \frac{p-c_\beta}{q-c_\beta} \bigr)^k ].
  \end{align}
Since $p$ is in $\Omega_\beta$, and $q$ is at least one cell width away,
it's not difficult to show that that $\abs{p-c_\beta}/\abs{q-c_\beta} \le
\sqrt{2}/3 < 1$. Therefore, the terms in the infinite sum above decay exponentially
fast at the rate $(\sqrt{2}/3)^k$. This means that we can truncate the series to
only $K$ terms, where $K$ is independent of the source curve $S_\beta$ and 
the query point $q$, and depends only on the desired accuracy.

Applying this observation to the Green's function kernel, we have
\begin{equation}
G(p,q) = \frac{1}{2\pi} \sum_{k=0}^{K} O_k(q-c_{\beta})I_k(p-c_{\beta}),
\label{eq:G_sep}
\end{equation}
where
\begin{equation}
I_k(p-c_{\beta})=\frac{(p-c_{\beta})^k}{k!} \quad \text{for} \quad k \geq 0,
\label{eq:Ik0}
\end{equation}
and
\begin{equation}
O_k(q-c_{\beta})=\frac{(k-1)!}{(q-c_{\beta})^k} \quad \text{for} \quad k \geq 1 \quad \text{and} \quad O_0(q-c_{\beta})=-\log(q-c_{\beta}).
\label{eq:Ok0}
\end{equation}
Thus, to evaluate the potential $u^q$ created by the density $\sigma$ on the curve
$S_\beta$, 
  \begin{align}
u^q = \underbrace{\int_{S_\beta} G(p,q) \sigma(p)
dS(p),}_{\textrm{target-from-source}}
\label{eq:t-f-s}
  \end{align}
all that is needed are the coefficients $\hat{\bm{\sigma}}^\beta$,
\begin{equation}
\hat{\bm{\sigma}}^\beta_k
= \underbrace{\int_{S_{\beta}} I_k(p-c_{\beta})\sigma(p)
dS(p),}_{\textrm{outgoing-from-source}}
\label{eq:o-f-s}
\end{equation}
for $k=0,1,\ldots,K$, from which the potential ${u}^q$ is evaluted by the formula
\begin{equation}
\begin{aligned}
u^q = 
\underbrace{\frac{1}{2\pi}  \sum_{k=0}^{K} O_k(q-c_{\beta})
\hat{\bm{\sigma}}^\beta_k.}_{\textrm{target-from-outgoing}}
\end{aligned}
\label{eq:t-f-o}
\end{equation}

Suppose that the integrals in the \textit{target-from-source} and
\textit{outgoing-from-source} formulas require $M$ degrees of freedom to
evaluate. 
To create the \textit{outgoing expansion}, we expend $KM=\mathcal{O}(M)$
operations to evaluate the \textit{outgoing-from-source} formulas.  If we
then evaluate the potential at $N$ well-separated target points $q$ using
the \emph{target-from-outgoing} formula, the cost will be only
$KN=\mathcal{O}(N)$. This is much faster than evaluating the potential
directly using the \textit{target-from-source} formula, since that would
require $MN$ operations.

\subsubsection{Incoming expansion}
\begin{wrapfigure}[11]{r}{1.2in}
  \includegraphics[width=\linewidth,trim={6mm 0mm 0mm 0mm}]{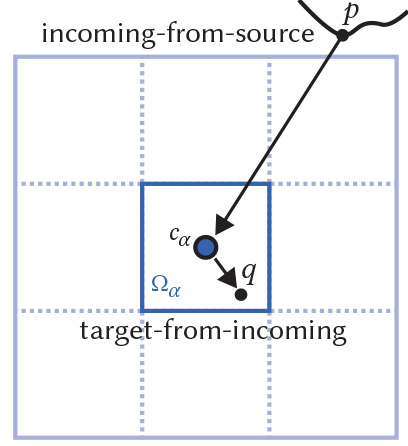}
  \label{fig:in_exp}
\end{wrapfigure} 
There is another alternative way of accelerating the evaluation of the
potential induced by an integral of a Green's function over a source curve.
 Suppose that we have a query point $q$
contained in a cell $\Omega_\alpha$, and that there is a source curve $S$
that is separated from $\Omega_\alpha$ by at least on cell width (see inset). 
We can construct an \textbf{incoming expansion} on $\Omega_\alpha$ 
for the potential created by this
source curve, which is directed inwards, in the sense that it is valid only
for target points $q$ inside $\Omega_\alpha$.

We can express the Green's function kernel as:
\begin{equation}
G(p,q) = \frac{1}{2\pi} \sum_{k=0}^{K} O_k(p-c_{\alpha})I_k(q-c_{\alpha}),
\end{equation}
where $q$ is in $\Omega_\alpha$ and $p$ is at least one cell-width away.
Thus, to evaluate the potential $u^q$ generated by the density $\sigma$ on the 
curve $S$ (see Eq.~\ref{eq:t-f-s}),
all that is needed are the coefficients $\hat{\mathbf{u}}^\alpha$,
\begin{equation}
\hat{\mathbf{u}}^\alpha_k = \underbrace{\int_{S}
O_k(p-c_{\alpha})\sigma(p)dS(p)}_{\textrm{incoming-from-source}},
\label{eq:i-f-s}
\end{equation}
for $k=0,1,\ldots,K$, from which the potential $u^q$ is evaluated by the
formula
\begin{equation}
\begin{aligned}
u^q = \underbrace{\frac{1}{2\pi}\sum_{k=0}^{K} I_k(q-c_{\alpha})
\hat{\mathbf{u}}^\alpha_k.}_{\textrm{target-from-incoming}}
\end{aligned}
\label{eq:t-f-i}
\end{equation}
Suppose, like before, that the integrals in the \textit{target-from-source} and
\textit{incoming-from-source} formulas require $M$ degrees of freedom to
evaluate. 
To create an \textit{incoming expansion}, we expend $KM=\mathcal{O}(M)$
operations to evaluate the \textit{incoming-from-source} formulas.  If we
then evaluate the potential at $N$ target points $q$ inside $\Omega_\alpha$
using the \emph{target-from-incoming} formula, then the cost will be only
$KN=\mathcal{O}(N)$. This is much faster than evaluating the potential
directly using the \textit{target-from-source} formula, since that would
require $MN$ operations.

\subsubsection{Incoming-from-Outgoing}

Suppose now that we have $N$ query points contained in a cell
$\Omega_{\alpha}$ with center $c_{\alpha}$, and that there are $m$ cells
containing source curves, all well-separated from the cell
$\Omega_\alpha$. If we construct \textit{outgoing-from-source} expansions
for each source cell, then the cost of evaluating the potential using the
\textit{target-from-outgoing} expansions will be $m K$ for each query point
$q$ in $\Omega_\alpha$, for a total evaluation cost of $NmK=\mathcal{O}(Nm)$.
However, if we could construct a single incoming
expansion on $\Omega_\alpha$ from all $m$ outgoing expansions, then the cost of
evaluating the
potential using the \textit{target-from-incoming} expansion would be only
$K$ for each query
point $q$ in $\Omega_\alpha$, for a total evaluation cost of
$NK=\mathcal{O}(N)$.
We can
construct an incoming expansion from an outgoing expansion, as follows.

Suppose that $\hat{\bm{\sigma}}^\beta$ is an outgoing expansion for the cell
$\Omega_\beta$, where $\Omega_\beta$ is well-separated from $\Omega_\alpha$.
Recall the \textit{target-from-outgoing} formula
\begin{equation}
\begin{aligned}
&u^q = 
\frac{1}{2\pi}  \sum_{k=0}^{K} O_k(q-c_{\beta})
\hat{\bm{\sigma}}^\beta_k \\
&= \frac{1}{2\pi}  \sum_{k=0}^{K} O_k[(q-c_{\alpha}) + (c_\alpha - c_\beta)]
\hat{\bm{\sigma}}^\beta_k,
\end{aligned}
\end{equation}
where $u^q$ is the potential at the query point $q$ in $\Omega_\alpha$.
If there are two terms inside the $O_{k}$ function, we can separate them
using the formula
\begin{equation}
O_k(z_1+z_2) = \sum_{l=0}^{K} (-1)^l O_{k+l}(z_1)I_l(z_2).
\label{eq:Ok}
\end{equation}
Applying Eq.~\ref{eq:Ok} and exchanging the order of summation, we have
\begin{equation}
\begin{aligned}
u^q = \sum_{k=0}^{K} I_k(q-c_{\alpha})
\hat{\mathbf{u}}^\alpha_k,
\end{aligned}
\end{equation}
where
\begin{equation}
\hat{\mathbf{u}}^\alpha_l = \underbrace{(-1)^l \sum_{k=0}^{K} O_{l+k}
(c_{\alpha}-c_{\beta})
\hat{\bm{\sigma}}^\beta_k. }_{\textrm{incoming-from-outgoing}}
\label{eq:i-f-o}
\end{equation}
Thus, we can turn an \textit{outgoing expansion} into an \textit{incoming
expansion} with $K^2$ operations. The cost of turning all outgoing
expansions into incoming expansions in the previous example is thus
$K^2m=\mathcal{O}(m)$. Recalling that the cost of evaluating the incoming
expansion is $\mathcal{O}(N)$, we have a total cost of $\mathcal{O}(m+N)$,
which is much better than the $\mathcal{O}(Nm)$ cost of evaluating all of
the \textit{outgoing expansions} naively.
%
 %\vspace{-20pt}
\begin{figure}[h]
  \includegraphics[width=\linewidth]{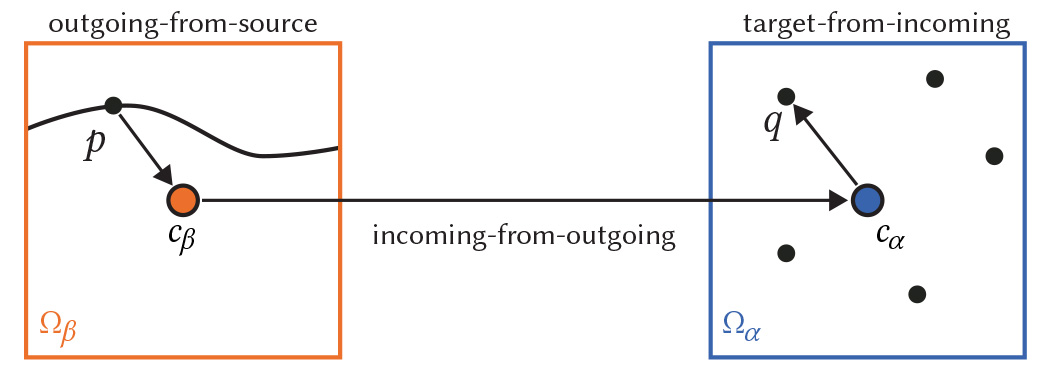}
  %\caption{.}
\label{fig:multipole_exp}
\end{figure}
 %\vspace{-10pt}
 
\subsection{A Single-level Method}
\label{sec:single_level}
Suppose that the computational domain contains a collection of diffusion
curves, and that all of their integrals can be discretized using a total of
$M$ degrees of freedom. Suppose further that we would like to evaluate the
potential induced by these curves at $N$ points. Suppose finally that the
computational domain containing all curves and evaluation points is divided
into $m$ boxes or cells. The three expansions described in the previous
section can be used to accelerate the evaluation of the potential.

First, we define the following relations betweens cells in our computational
domain.

\begin{minipage}{.67\linewidth}
\begin{itemize}[label={}, labelsep=0pt, leftmargin=0pt]
\item  The \textbf{neighbor list} $\mathcal{L}^\textrm{nei}_{\tau}$ of the
cell $\Omega_\tau$ (dark blue cell inset) is the set of all boxes that
directly touch $\Omega_\tau$ (light blue cells inset).
\item  The \textbf{well-separated list} $\mathcal{L}^\textrm{sep}_\tau$ of 
the cell $\Omega_\tau$ is the set of all boxes that do not touch 
$\Omega_\tau$ (orange cells inset).
\end{itemize}
\end{minipage}
\begin{minipage}{.28\linewidth}
\includegraphics[width=\linewidth]{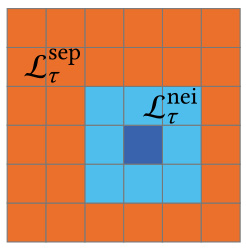}
\end{minipage}

The algorithm proceeds as follows.  An outgoing expansion is constructed for
every box using the \textit{outgoing-from-source} formula Eq.~\ref{eq:o-f-s}.
Next, for
each box $\Omega_\tau$, all outgoing expansions in the well-separated list
$\mathcal{L}^\textrm{sep}_\tau$ are turned into incoming expansions using
the \textit{incoming-from-outgoing} formula Eq.~\ref{eq:i-f-o}. 
Finally, for each box
$\Omega_\tau$, all incoming expansions are evaluated at the target points
using the \textit{target-from-incoming} formula Eq.~\ref{eq:t-f-i}, and then
added to the potentials produced by all sources in the \textit{neighbor
list} $\mathcal{L}^\textrm{nei}_\tau$ using the \textit{target-from-source}
formula Eq.~\ref{eq:t-f-s}. If $m$ is chosen to be $m=(NM)^\frac{1}{3}$, then
we show in Appendix~\ref{app:single_level} that the overall cost is
$\mathcal{O}((NM)^\frac{2}{3})$, which is substantially better than the
$\mathcal{O}(NM)$ cost of computing the potential at all $N$ points naively.

\subsection{Moving Between Levels in the Fast Multipole Method}

The single-level method described in the previous section can be further 
improved by defining a multi-level hierarchy of boxes on the computational
domain, and taking advantage of the fact that the number of terms $K$ in the
incoming and outgoing expansions stays constant across levels in this
hierarchy. The resulting  method is called the Fast Multipole Method (FMM).

The primary geometric data structure used by the FMM is the
\textbf{quadtree}, which is constructed recursively by starting from a
single cell containing the whole domain, and splitting each cell into 4
equal-sized smaller cells until each cell contains a small number of degrees
of freedom, or the algorithm reaches a prescribed maximum depth.

\subsubsection{Relations between Cells} 
\label{sec:cell_relations}
To describe the multi-level scheme, we must first introduce several additional
relations between cells in the hierarchy, besides the \textit{neighbor list}
$\mathcal{L}^\textrm{nei}_\tau$ and 
\textit{well-separated list} $\mathcal{L}^\textrm{sep}_\tau$ 
of the box $\Omega_\tau$, introduced in the
previous section.

\begin{itemize}[label={}, labelsep=0pt, leftmargin=5pt]
\item The \textbf{parent} of a cell $\Omega_\tau$ is the cell on the next
coarsest level that contains $\Omega_{\tau}$.
\item The \textbf{child list} of a cell $\Omega_\tau$ is the set
$\mathcal{L}^\textrm{child}_\tau$ containing the 4 cells whose parent is
$\Omega_{\tau}$.
\item The \textbf{leaf} cells are the cells that don't have any children.
\item The \textbf{interaction list} of a cell $\Omega_\tau$ is the 
set $\mathcal{L}^\textrm{int}_\tau$ containing all cells $\Omega_\sigma$
such that: 1) $\Omega_\tau$ and $\Omega_\sigma$ are on the same level 
and well-separated to each other,  
2) the parents of $\Omega_\tau$ and $\Omega_\sigma$ is not well-separated
\end{itemize}

Note that the \textit{interaction list} is a subset of the
\textit{well-separated list} defined in the previous section, since two boxes
on the same level cannot touch without their parents also touching.

To proceed, we need a way of moving information between levels in the
multi-level hierarchy, in the form of outgoing and incoming expansions.
The formulas relating expansions between different levels of the quadtree
are also called translation operators.

\vspace{10pt}

\subsubsection{Outgoing-from-Outgoing}
\begin{wrapfigure}[9]{r}{0.9in}
  \includegraphics[width=\linewidth,trim={6mm 0mm 0mm 0mm}]{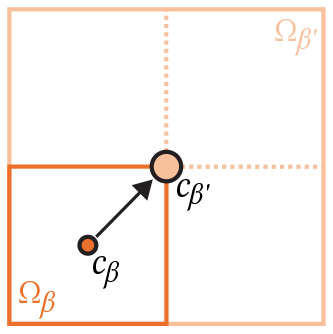}
  %\label{fig:o-f-o}
\end{wrapfigure} 
Consider a parent cell $\Omega_{\beta^{'}}$ containing 4 child cells.
Suppose that $\Omega_\beta$ is a child of $\Omega_{\beta'}$, and that we
have already computed the outgoing expansion $\hat{\bm{\sigma}}^\beta$ on
$\Omega_\beta$. Then, we can transfer the outgoing expansions from the child
to the parent, as follows.  First, recall the \textit{target-from-outgoing}
formula
\begin{equation}
\begin{aligned}
&u^q = 
\frac{1}{2\pi}  \sum_{k=0}^{K} O_k(q-c_{\beta})
\hat{\bm{\sigma}}^\beta_k \\
&= \frac{1}{2\pi}  \sum_{k=0}^{K} O_k[(q-c_{\beta'})+(c_{\beta'}-c_\beta)]
\hat{\bm{\sigma}}^\beta_k.
\end{aligned}
\end{equation}
Applying Eq.~\ref{eq:Ok} and exchanging the order of summation, we have
\begin{equation}
\begin{aligned}
u^q = 
\frac{1}{2\pi}  \sum_{k=0}^{K} O_k(q-c_{\beta'})
\hat{\bm{\sigma}}^{\beta'}_k
\end{aligned}
\end{equation}
where
\begin{equation}
\hat{\bm{\sigma}}^{\beta'}_k = \underbrace{\sum_{l=0}^{k}
I_{k-l}(c_{\beta}-c_{\beta^{'}}) \hat{\bm{\sigma}}^{\beta}_l.
}_{\textrm{outgoing-from-outgoing}}
\label{eq:o-f-o}
\end{equation}
This \textit{outgoing-from-outgoing} formula transfers an outgoing expansion
from the child cell $\Omega_\beta$ to the parent cell $\Omega_{\beta'}$ 
%(see Figure~\ref{fig:exp_trans}).

\subsubsection{Incoming-from-Incoming}
\begin{wrapfigure}[9]{r}{0.9in}
  \includegraphics[width=\linewidth,trim={6mm 0mm 0mm 3mm}]{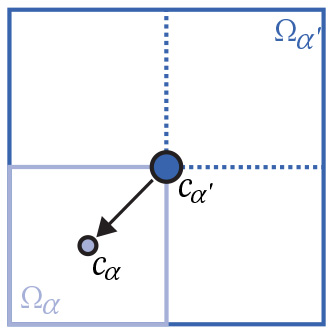}
  %\label{fig:o-f-o}
\end{wrapfigure} 
Likewise,  an incoming expansion can be transferred from a parent cell to a
child cell. Suppose that $\Omega_\alpha$ is a child of $\Omega_{\alpha'}$,
and that we have already computed the incoming expansion of
$\hat{\bm{u}}^{\alpha'}$ of $\Omega_{\alpha'}$. We can transfer the incoming 
expansion from the parent to the child, as follows.
Recall the $\textit{target-from-incoming}$ formula
\begin{equation}
\begin{aligned}
&u^q = \frac{1}{2\pi} \sum_{l=0}^{K}  I_l(q-c_{\alpha'}) \hat{\bm{u}}^{\alpha'}_l \\
&= \frac{1}{2\pi} \sum_{l=0}^{K} 
I_l[(q-c_{\alpha})+(c_{\alpha}-c_{\alpha'})] \hat{\bm{u}}^{\alpha'}_l.
\end{aligned}
\end{equation}
Similarly to Eq.~\ref{eq:Ok}, when there are two terms inside the $I_l$ function, 
we can separate them:
\begin{equation}
I_l(z_1+z_2) = \sum_{k=0}^{l} I_{l-k}(z_1)I_k(z_2).
\label{eq:Ik}
\end{equation}
Applying Eq.~\ref{eq:Ik}, we obtain the following equation:
\begin{align}
&u^q = \frac{1}{2\pi} \sum_{l=0}^{K} I_l(q-c_{\alpha})
\hat{\bm{u}}^{\alpha}_l,
\end{align}
where,
\begin{equation}
\hat{\bm{u}}^{\alpha}_l = \underbrace{\sum_{k=l}^{K}
I_{k-l}(c_{\alpha}-c_{\alpha'})
\hat{\bm{u}}_k^{\alpha'}}_{\textrm{incoming-from-incoming}}.
\label{eq:i-f-i}
\end{equation}
This \emph{incoming-from-incoming} formula transfers an incoming expansion from
the parent cell $\Omega_{\alpha'}$ to its child cell $\Omega_{\alpha}$.

\subsection{The Fast Multipole Method on a Uniform Quadtree}
\label{sec:fmm_uni_quadtree}
A uniform quadtree is a quadtree in which all leaf cells occur at the same
level of the tree.  On a uniform quadtree, the Fast Multipole Method
proceeds as follows.  Suppose that the computational domain contains a collection
of diffusion curves which can be discretized using $M$ degrees of freedom. Suppose
further that we would like to evaluate the potential induced by these
curves at $N$ points. Suppose finally that we construct our quadtree by subdividing
until fewer than $b$ degrees of freedom are contained in each leaf box.
To proceed, first, all outgoing expansions are formed on all of
the \textit{leaf} cells on the lowest level in the hierarchy using the
\textit{outgoing-from-source} formula Eq.~\ref{eq:o-f-s}. 
Then, outgoing expansions are
formed on coarser grids, going from finest to coarsest, by merging outgoing
expansions using the \textit{outgoing-from-outgoing} formula
Eq.~\ref{eq:o-f-o}.  

Once all
outgoing expansions on all levels have been constructed, incoming expansions
are formed, going from the coarsest grid to the finest. In particular,
incoming expansions are first formed on all cells on the coarsest grid by
applying the \textit{incoming-from-outgoing} formula Eq.~\ref{eq:i-f-o}
to all outgoing
expansions in the well-separated list $\mathcal{L}^\text{sep}_\tau$ of each
cell $\Omega_\tau$. Then, on each successive level, the incoming expansion of
each cell $\Omega_\tau$ is formed by first transferring the incoming
expansion from the parent cell $\Omega_{\tau'}$ using the
\textit{incoming-from-incoming} formula Eq.~\ref{eq:i-f-i}, and then by
turning all outgoing expansions in the \textit{interaction list}
$\mathcal{L}_\tau^\text{int}$ into incoming expansions, using the
\textit{incoming-from-outgoing} formula Eq.~\ref{eq:i-f-o}. 

Finally, for each cell
$\Omega_\tau$, the incoming expansions are evaluated at the target points
using the \textit{target-from-incoming} formula Eq.~\ref{eq:t-f-i}, 
and then these values
are added to the potentials from the sources in the \textit{neighbor list}
$\mathcal{L}_\tau^\text{nei}$, which are evaluated directly using the
\textit{target-from-source} formula Eq.~\ref{eq:t-f-s}.

In Appendix~\ref{app:fmm_cost}, we show that if we choose $b = K$, then the
total cost is $\O(KM + KN)$, so it is linear in both the number of degrees of
freedom in the sources $M$ and the number of targets $N$.

\subsection{Extension to a Non-uniform Quadtree}
\label{sec:non_uniform_quadtree}

If we allow leaf cells to occur on different levels of the quadtree,
then the quadtree is called non-uniform.  For such a quadtree, the method 
described in Section~\ref{sec:fmm_uni_quadtree} can fail,  
since there may be cells which end up unaccounted for. 

\subsubsection{Relations between Cells}

We will need to introduce two more relations between cells, in addition to
the ones introduced already in Section~\ref{sec:single_level}
and~\ref{sec:cell_relations}.

\begin{itemize}[label={}, labelsep=0pt, leftmargin=5pt]
\item The \textbf{smaller separated list}
$\mathcal{L}^\text{small}_\alpha$ of the leaf cell $\Omega_\alpha$ is the
set of cells $\Omega_\beta$ that are smaller than the cell
$\Omega_{\alpha}$, such that $\Omega_\alpha$ is in the well-separated list
$\mathcal{L}^\text{sep}_\beta$ of $\Omega_\beta$, and $\Omega_\alpha$ is not
in the well-separated list $\mathcal{L}^\text{sep}_{\beta'}$ of the parent
cell $\Omega_{\beta'}$ of $\Omega_\beta$ (see Fig.~\ref{fig:quad_list} (d)).
\item The \textbf{bigger separated list} $\mathcal{L}^\text{big}_\alpha$ is the
dual of the \textit{smaller separated list}, in the sense that a cell
$\Omega_\beta$ is in $\mathcal{L}^\text{big}_\alpha$ if and only if
$\Omega_\alpha$ is in $\mathcal{L}^\text{small}_\beta$.  It is not too hard
to show that the \textit{bigger separated list}
$\mathcal{L}^\text{big}_\alpha$ of the cell $\Omega_\alpha$ is the set of
leaf cells $\Omega_\beta$ that are larger than the cell $\Omega_{\alpha}$,
such that that $\Omega_\beta$ is  in the well-separated list
$\mathcal{L}^\text{sep}_\alpha$ of $\Omega_\alpha$, and $\Omega_\beta$ is
not in the well-separated list $\mathcal{L}^\text{sep}_{\alpha'}$ of the
parent cell $\Omega_{\alpha'}$ of $\Omega_\alpha$ (see
Fig.~\ref{fig:quad_list} (e)).
\end{itemize}

In the FMM literature, the lists $\mathcal{L}^\text{small}_\tau$ and
$\mathcal{L}^\text{big}_\tau$ are sometimes called list~3 and list~4,
and denoted by $\mathcal{L}^\text{(3)}_\tau$ and $\mathcal{L}^\text{(4)}_\tau$,
respectively \cite{martinsson2019fast}.

It is not difficult to show that the interactions between any two leaf cells
on a non-uniform quadtree are accounted for by the
\textit{outgoing-from-source} to \textit{incoming-from-outgoing} to
\textit{target-from-incoming} calculation, if and only if those two leaf
cells are in each other's \textit{separated lists} (see
Fig.~\ref{fig:exp_relation}). If two leaf cells $\Omega_\alpha$ and
$\Omega_\beta$ are not mutually well-separated, then one of the following
three cases must occur: $\Omega_\alpha$ and $\Omega_\beta$ are touching;
$\Omega_\beta$ (or one of its parent cells on a sufficiently high level) is
in the \textit{smaller separated list} of $\Omega_\alpha$; or $\Omega_\beta$ is
in the \textit{bigger separated list} of $\Omega_\alpha$ (or one of its
parent cells on a sufficiently high level).

When two cells are not mutually well-separated and are not touching,
it is necessary to skip part of \textit{incoming-from-outgoing} and  
directly proceed \textit{target-from-outgoing} or \textit{incoming-from-source}
%
 %\vspace{-5pt}
\begin{figure}[h]
 \includegraphics[width=\linewidth]{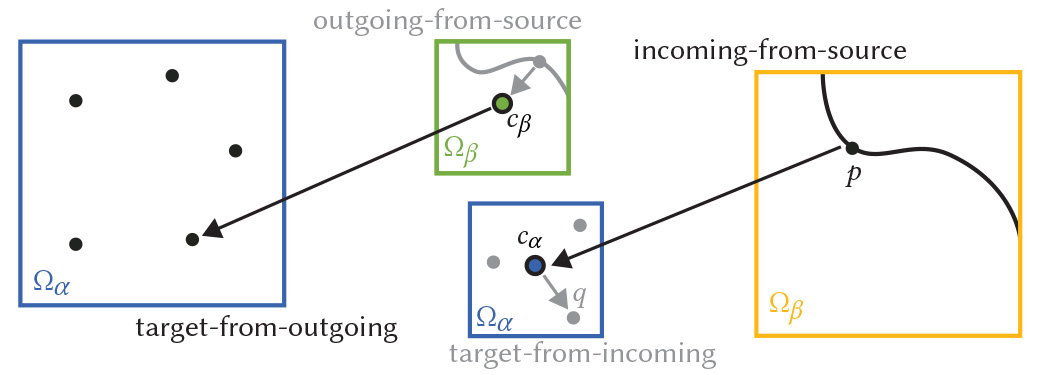}
%\caption{\emph{incoming-from-source} , \emph{target-from-outgoing}}
\label{fig:nonuni_exp}
\end{figure}
 %\vspace{-15pt}
%

\subsubsection{Target-from-Outgoing}

Suppose that $\Omega_\beta$ is in the \textit{smaller separated list}
$\mathcal{L}_\alpha^\text{small}$ of the leaf cell $\Omega_\alpha$.  In this
case, we cannot use the incoming expansion on $\Omega_\beta$, but the
outgoing expansion on $\Omega_\alpha$ is still valid.  Hence, we use the
\emph{target-from-outgoing} formula Eq.~\ref{eq:t-f-o}, skipping the incoming
expansion step (see above Figure and Fig.~\ref{fig:exp_relation}).

\subsubsection{Incoming-from-Source}

Suppose that the leaf cell $\Omega_\beta$ is in the \textit{bigger separated
list} $\mathcal{L}_\alpha^\text{big}$ of $\Omega_\alpha$.  In this case, we
cannot use the outgoing expansion on $\Omega_\beta$, but the incoming
expansion on $\Omega_\alpha$ is still valid.  Hence, we use the
\emph{incoming-from-source} formula Eq.~\ref{eq:i-f-s}, skipping the
outgoing expansion step (see above Figure and Fig.~\ref{fig:exp_relation}).

\subsection{Interaction between Cells}

Through Sec.~\ref{sec:multipole_exp} to Sec.~\ref{sec:non_uniform_quadtree},
we have established various routes of integration depending on the cell relations.
In summary, there are 4 paths of integration between \textit{source} and \textit{target}.

\begin{enumerate}
\item \vspace{1ex}\begin{rounded1} \textit{target-from-source} (direct
integration) \end{rounded1}

\item \begin{rounded2} \textit{outgoing-from-source} \textrightarrow
\textit{incoming-from-outgoing} \textrightarrow
\textit{target-from-outgoing} \end{rounded2}

\item \begin{rounded3} \textit{outgoing-from-source} \textrightarrow
\textit{target-from-outgoing} \end{rounded3}

\item \begin{rounded4} \textit{incoming-from-source} \textrightarrow
\textit{target-from-incoming} \end{rounded4}
\end{enumerate}

\begin{figure}[h]
  \includegraphics[width=\linewidth]{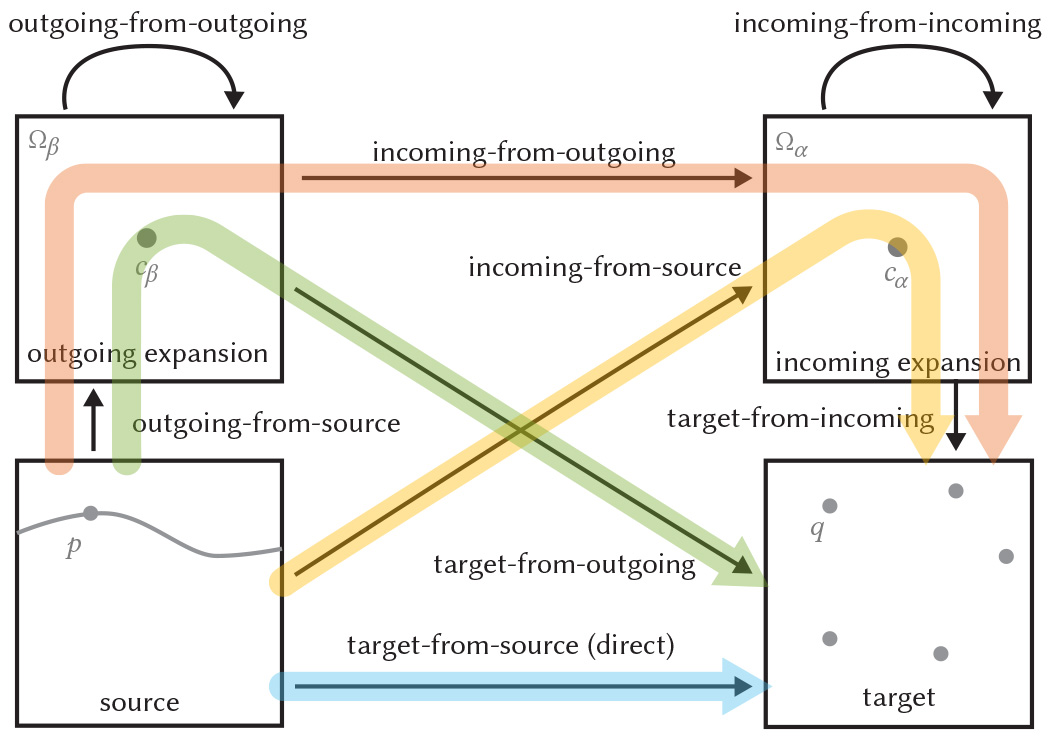}
  \caption{Diagram describing the flow of computations in the FMM
  by various routes. }
\label{fig:exp_relation}
\end{figure}

Fig.~\ref{fig:exp_relation} visualizes each path of integration 
between cells.

\subsection{The Fast Multipole Method}
\label{sec:fmm}

Suppose that the computational domain contains a collection of diffusion
curves which can be discretized using $M$ degrees of freedom. Suppose
further that we would like to evaluate the potential induced by these curves
at $N$ points. Suppose finally that we construct our quadtree by subdividing
until fewer than $b$ degrees of freedom are contained in each leaf box,
allowing leaf boxes of different sizes.  The Fast Multipole
Method, on a nonuniform quadtree, proceeds as follows. 
First, all outgoing expansions are formed on all of the
\text{leaf} cells on the lowest level in the hierarchy using the
\textit{outgoing-from-source} formula Eq.~\ref{eq:o-f-s}.
Next, the outgoing expansions are
formed on coarser grids, going from finest to coarsest, by merging outgoing
expansions using the \textit{outgoing-from-outgoing} formula Eq.~\ref{eq:o-f-o}. 

Once all outgoing expansions on all levels are formed, incoming expansions
are constructed, going from the coarsest grid to the finest. 
In particular,
incoming expansions are first formed on all cells on the coarsest grid by
applying the \textit{incoming-from-outgoing} formula Eq.~\ref{eq:i-f-o}
to all outgoing
expansions in the well-separated list $\mathcal{L}^\text{sep}_\alpha$ of each
cell $\Omega_\alpha$. 
Then, on each successive level, the incoming expansion of each
cell $\Omega_\alpha$ is formed from the incoming
expansion of the parent cell $\Omega_{\alpha'}$, by
transferring the expansion
to $\Omega_\alpha$ using the \textit{incoming-from-incoming}
formula Eq.~\ref{eq:i-f-i}. These incoming expansions are added to
the incoming expansions constructed from all
of the cells in the \textit{interaction list}
$\mathcal{L}^\text{int}_\alpha$ using the \textit{incoming-from-outgoing}
formula Eq.~\ref{eq:i-f-o},  and from all cells in the \textit{bigger
separated list} $\mathcal{L}^\text{big}_\alpha$ using the
\textit{incoming-from-source} formula Eq.~\ref{eq:i-f-s}.

For each leaf cell $\Omega_\alpha$, the potential is evaluated at the target
points using the \textit{target-from-incoming} formula Eq.~\ref{eq:t-f-i}. 
This potential
is then added to the potential produced by every cell $\Omega_\beta$ is
the \textit{smaller separated list} $\mathcal{L}^\text{small}_\alpha$ using
the \textit{target-from-outgoing} formula Eq.~\ref{eq:t-f-o}.
Finally, these values are added to the potentials from the source cells in
the \textit{neighbor list} $\mathcal{L}_\alpha^\text{nei}$, which are
evaluated directly using the \textit{target-from-source}
formula Eq.~\ref{eq:t-f-s}.

\begin{figure}
  \includegraphics[width=\linewidth]{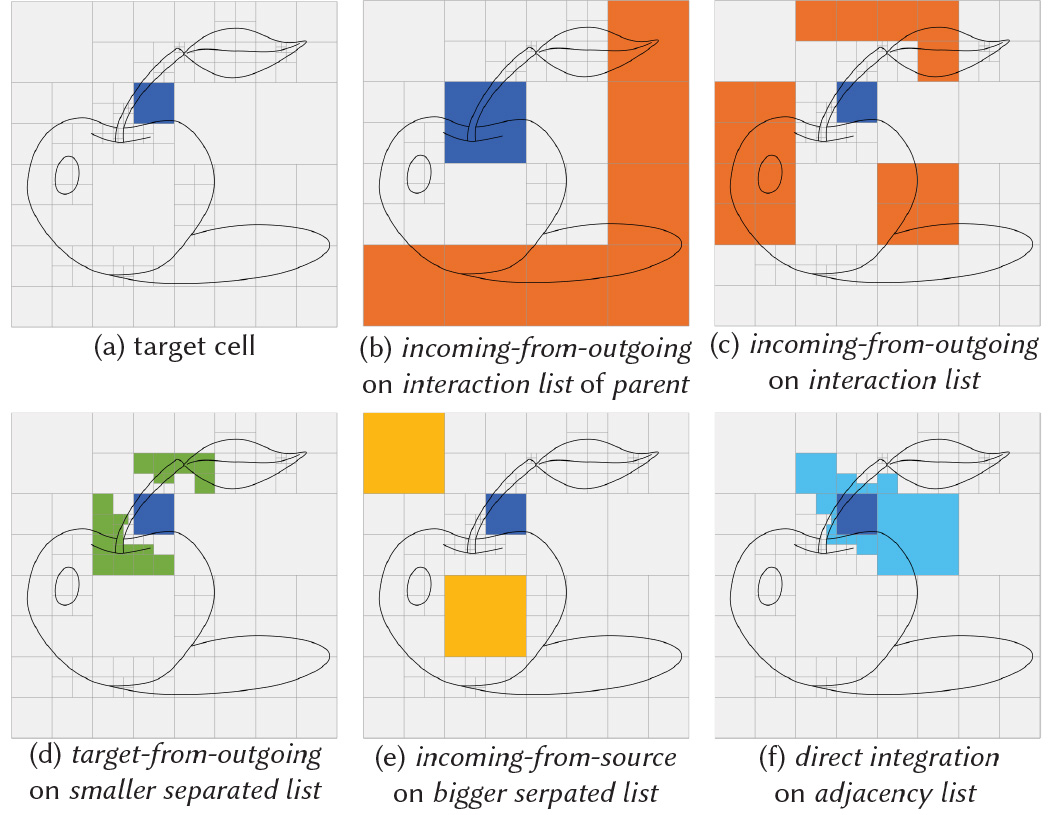}
  \caption{If a target point is contained in the blue
  cell (a), assuming the outgoing expansions are computed on every cell,
  \emph{incoming-from-outgoing} expansions are computed on its
  \emph{interaction list} from coarsest to finest (b), and transferred to
  its children. At its finest level, potentials are computed by summing over
  \emph{incoming-from-outgoing} on its \emph{interaction list} (c),
  \emph{incoming-from-incoming} on its parent, \emph{target-from-outgoing}
  on its \emph{smaller separated list} (d), \emph{incoming-from-source} on
  its \emph{bigger separated list} (e), \emph{target-from-incoming} from its
  own cell, and \emph{target-from-source} (direct integration) on its
  \emph{adjacency list} (f). }
\label{fig:quad_list}
\end{figure}

The overall algorithm in pseudo code is described in Table~\ref{alg:fmm} and 
the cell relations with an example set of input curves is visualized in
Fig.~\ref{fig:quad_list}.  Recalling that the incoming and outgoing
expansions contain only $K$ terms,
if we set $b = K$, then we have,  by an essentially identical argument to the
one provided in Appendix~\ref{app:fmm_cost},
that the total cost is $\mathcal{O}(KM+KN)$, so it is
linear in both the number of degrees of freedom in the source curves $M$,
and the number of targets $N$.

\begin{algorithm}
\caption{Fast Multipole Method}
\begin{flushleft}
 \hspace*{\algorithmicindent} \textbf{Inputs:} source curves $S$, 
 density values $\sigma$, target
 points $q$ \\
 \hspace*{\algorithmicindent} \textbf{Outputs:} target values $u^q$ 
 \end{flushleft}
\begin{algorithmic}[1]
\algloop{upward}
\algloop{downward}
\algloop{from}
	\Loop \, over every leaf cell $\Omega_\beta$:
		\State $\hat{\bm{\sigma}}^\beta = \text{outgoing-from-source}(S_\beta,\sigma)$
	\EndLoop
	\upward \, (from finest to coarsest level):
		\Loop \, over every cell $\Omega_{\beta'}$ in level:
      \State Initialize $\hat{\bm{\sigma}}^{\beta'} = \bm{0}$
      \Loop \, over every child cell $\Omega_\beta \in
      \mathcal{L}_{\beta'}^\text{child}$: 
				\State $\hat{\bm{\sigma}}^{\beta'} = \hat{\bm{\sigma}}^{\beta'} +
        \text{outgoing-from-outgoing}(\hat{\bm{\sigma}}^{\beta})$
      \EndLoop
		\EndLoop
	\downward \, (from coarsest to finest level):
		\Loop \, over every cell $\Omega_\alpha$ in level:
      \State Initialize $\hat{\bm{u}}^\alpha = \bm{0}$
			\Loop \, over every cell $\Omega_\beta$ in interaction list $\mathcal{L}_\alpha^\text{int}$
      (or separated list $\mathcal{L}_\alpha^\text{sep}$ when at the top level):
        \State $\hat{\bm{u}}^\alpha = \hat{\bm{u}}^\alpha +
        \text{incoming-from-outgoing}(\hat{\bm{\sigma}}_{\beta})$
      \EndLoop
      \Loop \, over every cell $\Omega_\beta$ in bigger separated list
      $\mathcal{L}_\alpha^\text{big}$: 
        \State $\hat{\bm{u}}^\alpha = \hat{\bm{u}}^\alpha +
        \text{incoming-from-source}(S_\beta,\sigma)$
      \EndLoop
			\from \, the parent cell $\Omega_{\alpha'}$, when it exists:
        \State $\hat{\bm{u}}^\alpha = \hat{\bm{u}}^\alpha +
        \text{incoming-from-incoming}(\hat{\bm{u}}^{\alpha'})$
		\EndLoop
	\Loop \, over every target $q$ in every leaf cell $\Omega_\alpha$:
    \State Initialize $u^q = 0$
		\from \, own cell:
			\State $u^q = u^q + \text{target-from-incoming}(\hat{\bm{u}}^\alpha)$
		\Loop \, over every 
    cell $\Omega_\beta$ in smaller separated list
    $\mathcal{L}_\alpha^\text{small}$:
      \State $u^q = u^q +
      \text{target-from-outgoing}(\hat{\bm{\sigma}}^\beta)$
    \EndLoop
    \Loop \, over every leaf cell $\Omega_\beta$ in neighbor list
    $\mathcal{L}_\alpha^\text{nei}$:
			\State $u^q = u^q + \text{target-from-source}(S_\beta,\sigma)$
    \EndLoop
	\EndLoop
\end{algorithmic}
\label{alg:fmm}
\end{algorithm}

The final FMM Algorithm \ref{alg:fmm} is thus an algorithm 
for evaluating the single-layer potential integral operator Eq.~\ref{eq:int_V}
efficiently, in linear time in both the number of degress of freedom in the sources, 
and in the number of targets.  We denote the  evaluation of the 
single-layer integral operator using the FMM by:
  \begin{align}
u^q = \mathbf{FMM}_{G}(S,\sigma,q),
\label{eq:fmm_eval_app}
  \end{align}
where $S$ is a collection of source curves, $\sigma$ is a density, and $q$ is
a collection of target points.

\subsection{Integration of Outgoing Expansion of Green's function}
\label{app:int_out_exp_g}

Here, we describe the analytic integration of Eq.~\ref{eq:o-f-s}. 
Suppose the we have a line element starting at point $p_a$ and ending at $p_b$.
The integral can be expressed as: 
\begin{equation}
\hat{\bm{\sigma}}^\beta_k = \sigma \int_{p_a}^{p_b} I_k(p-c_{\beta}) dS(p) = \sigma \bar{w} [I_{k+1}(p_b-c_{\beta}) - I_{k+1}(p_a-c_{\beta})],
\label{eq:o-f-s_anlytic}
\end{equation}
where $w$ is the complex tangential vector along the boundary $\overline{p_a
p_b}$, and $\bar{w}$ its complex conjugate.

\subsection{Integration of Outgoing Expansion of Normal Derivative of
Green's function}
\label{app:int_out_exp_f}

Similarly, we describe the analytic integration of the outgoing
expansion of the normal derivative of the Green's function:
\begin{equation}
\hat{\bm{\mu}}^\beta_k = \mu n \int_{p_a}^{p_b} I_{k-1}(p-c_{\beta}) dS(p) = \mu n \bar{w} [I_{k}(p_b-c_{\beta}) - I_{k}(p_a-c_{\beta})].
\label{eq:o-f-s_anlytic_f}
\end{equation}

\subsection{Integration of Incoming Expansion of Green's function}
\label{app:int_inc_exp_g}

Here we describe the analytic integration of Eq.~\ref{eq:i-f-s}:
\begin{equation}
\hat{\mathbf{u}}^\alpha_k= \sigma \int_{p_a}^{p_b} O_k(p-c_{\alpha})dS(p) = \sigma \bar{w} [O_{k-1}(p_a-c_{\alpha}) - O_{k-1}(p_b-c_{\alpha})].
\label{eq:i-f-s_anlytic_g}
\end{equation}

\subsection{Integration of Incoming Expansion of Normal Derivative of
Green's function}
\label{app:int_inc_exp_f}
Similarly, we describe analytic integration of the incoming expansion of
the normal derivative of the Green's function:
\begin{equation}
\hat{\mathbf{v}}^\alpha_k= \mu n \int_{p_a}^{p_b} O_{k+1}(p-c_{\alpha})dS(p) = \mu n \bar{w} [O_{k}(p_b-c_{\alpha}) - O_{k}(p_a-c_{\alpha})].
\label{eq:i-f-s_anlytic_f}
\end{equation}

\section{Cost of the Single-Level Method}
\label{app:single_level}

Here, we estimate the computational cost of the single-level method
described in Section~\ref{sec:single_level}.  The total cost of computing
the potentials can be estimated as follows, recalling that all incoming and
outgoing expansions contain only $K$ terms:
  \begin{itemize}

\item The cost of constructing all outgoing expansions is $\O(MK)$, since
each degree of freedom in the discretization of the diffusion curves
contributes to $K$ outgoing expansion coefficients of the box containing
that degree of freedom.

\item The cost of computing all incoming expansions from all outgoing
expansions is $\O(K^2m^2)$, since each incoming expansion coefficient of each
box receives a contribution from each of the outgoing expansion coefficients
of each well-separated box.

\item The cost of evaluating the potential at all $N$ points using the
incoming expansions is $\O(KN + NM/m)$, since every point receives a contribution
from $K$ incoming expansion coefficients of the box containing that point, and
also receives a contribution from every degree of freedom in the $\O(1)$ boxes
in the \textit{neighbor} list, each of which contains $\O(M/m)$ degrees of freedom.

  \end{itemize}
The total cost is thus $\O(MK + K^2m^2 + KN + NM/m)$. It is not too hard to
see that this cost is minimized when the number of boxes $m$ is set to
$m=(NM)^{1/3}$, which gives an overall cost of $\O((NM)^{2/3})$.  This can be
substantially better than the $\O(NM)$ cost of computing the potential at all
$N$ points naively.

\section{Cost of the FMM on a Uniform Quadtree}
\label{app:fmm_cost}

Here, we estimate the computational cost of the FMM on a uniform quadtree
described in Section~\ref{sec:fmm_uni_quadtree}.  The cost can be estimated
as follows, once again recalling that all incoming and outgoing expansions
contain only $K$ terms.
\begin{itemize}

\item The cost of constructing all outgoing expansion from the sources in
the leaf boxes is $\O(KM)$, since each source point contributes to $K$
outgoing expansion coefficents in the containing cell.

\item The cost of transferring outgoing expansion from the child boxes to
their parents is $\O(K^2 M/b(1 + 1/4 + 1/16 + \cdots)) = \O(K^2M/b)$, since
there are $M/b$ boxes at the lowest level, $M/b\cdot 1/4$ at the next level,
and so on. The exact same reasoning shows that the cost of transferring the
outgoing exapnsions to incoming expansions is also $\O(K^2 M/b)$, as is the
cost of transferring the incoming expansions from parent boxes to child
boxes.

\item The cost of evaluating the potential at all $N$ target points using
the incoming expansions is $\O(KN + Nb)$, since every point receievs a
contribution from $K$ incoming expansion coefficients of the cell containing
that point, and also receives a contribution from every degree of freedom in
the $\O(1)$ boxes in the \textit{neighbor list}, each of which contains
$\O(b)$ degrees of freedom.

\end{itemize}

The total cost is thus $\O(KM + K^2M/b + KN + Nb)$. By choosing $b = K$,
we have that the total cost is $\O(KM + KN)$, so it is linear in both the
number of degrees of freedom in the sources $M$ and the number of targets
$N$.

\section{Precomputations for the FMM}
\label{app:pre_fmm}

When the Fast Multipole Method is used to evaluate the potential produced by
several different density functions $\sigma$ over a single set of
discretized curves $\bar{S}$ and target points $q$, a large number of
quantities that are independent of the density function can be precomputed.
In particular, many terms appearing in the various operators used by the FMM
can be precomputed. Below, we describe the precise quantities contained in
the precomputatons $\mathcal{P}_G$ and $\mathcal{P}_F$ for the single- and
double-layer respectively (omitting some of the quantities associated only
with the FMM for the double-layer):
\begin{equation}
\hat{\mathbf{u}}^\alpha_k = \underbrace{\int_{S}
O_k(p-c_{\alpha})\sigma(p)dS(p)}_{\textrm{incoming-from-source}}
= \sum_{i=1}^{M_\alpha} \bar{\sigma}_i \underbrace{\int_{(\bar{S}_\alpha)_i}
O_k(p-c_{\alpha})dS(p)}_{\in \mathcal{P}_G},
\label{eq:pre-i-f-s}
\end{equation}
\begin{equation}
\hat{\bm{\sigma}}^\beta_k
= \underbrace{\int_{S_{\beta}} I_k(p-c_{\beta})\sigma(p)
dS(p)}_{\textrm{outgoing-from-source}}
= \sum_{i=1}^{M_\alpha} \bar{\sigma}_i \underbrace{\int_{(\bar{S}_{\beta})_i} 
I_k(p-c_{\beta}) dS(p)}_{\in \mathcal{P}_G},
\label{eq:pre-o-f-s}
\end{equation}
\begin{equation}
\hat{\bm{\sigma}}^{\beta'}_k = \underbrace{\sum_{l=0}^{k}
I_{k-l}(c_{\beta}-c_{\beta^{'}}) \hat{\bm{\sigma}}^{\beta}_l
}_{\textrm{outgoing-from-outgoing}}
= \sum_{l=0}^{k}
\underbrace{I_{k-l}(c_{\beta}-c_{\beta^{'}})}_{\in \mathcal{P}_G, \mathcal{P}_F}
\hat{\bm{\sigma}}^{\beta}_l,
\label{eq:pre-o-f-o}
\end{equation}
\begin{equation}
\hat{\mathbf{u}}^\alpha_l = \underbrace{(-1)^l \sum_{k=0}^{K} O_{l+k}
(c_{\alpha}-c_{\beta})
\hat{\bm{\sigma}}^\beta_k}_{\textrm{incoming-from-outgoing}}
= (-1)^l \sum_{k=0}^{K} \underbrace{O_{l+k}
(c_{\alpha}-c_{\beta})}_{\in \mathcal{P}_G, \mathcal{P}_F}
\hat{\bm{\sigma}}^\beta_k,
\label{eq:pre-i-f-o}
\end{equation}
\begin{equation}
\hat{\bm{u}}^\alpha_l = \underbrace{\sum_{k=l}^{K}
I_{k-l}(c_{\alpha}-c_{\alpha^{'}})
\hat{\bm{u}}_k^{\alpha'}}_{\textrm{incoming-from-incoming}}
= \sum_{k=l}^{K}
\underbrace{I_{k-l}(c_{\alpha}-c_{\alpha^{'}})
}_{\in \mathcal{P}_G, \mathcal{P}_F}\hat{\bm{u}}_k^{\alpha'},
\label{eq:pre-i-f-i}
\end{equation}
\begin{align}
u^q = \underbrace{\sum_{k=0}^{K} I_k(q-c_{\alpha})
\hat{\mathbf{u}}^\alpha_k}_{\textrm{target-from-incoming}} = \sum_{k=0}^{K}
\underbrace{I_k(q-c_{\alpha})}_{\in \mathcal{P}_G, \mathcal{P}_F}
\hat{\mathbf{u}}^\alpha_k.
\label{eq:pre-t-f-i}
\end{align}
\begin{align}
u^q = 
\underbrace{\frac{1}{2\pi}  \sum_{k=0}^{K} O_k(q-c_{\beta})
\hat{\bm{\sigma}}^\beta_k}_{\textrm{target-from-outgoing}}
= 
\frac{1}{2\pi}  \sum_{k=0}^{K}
\underbrace{O_k(q-c_{\beta})}_{\in \mathcal{P}_G, \mathcal{P}_F}
\hat{\bm{\sigma}}^\beta_k.
\label{eq:pre-t-f-o}
\end{align}

\section{Pre-processing of Diffusion Curves}

We pre-process input curves to lessen the burden on our adaptive subdivision
algorithm.  We split the curves that have intersections, and also remove
(redundant) overlapping curves.

\end{appendices}

\end{document}